\documentclass[12pt]{book}

\usepackage[vmargin={3cm, 3.3cm},hmargin=3cm]{geometry}

\date{}

\usepackage[pdfpagelabels]{hyperref}

\usepackage{latexsym} 
\usepackage{amsmath}
\usepackage{amssymb}
\usepackage{amsthm}
\usepackage{mathrsfs}
\usepackage{amsfonts}
\usepackage{graphicx}

\usepackage{color}
\usepackage{bbm,verbatim}

\usepackage{minitoc} 

\input labelfig.tex

\usepackage{geometry}
\usepackage{multirow} 
\usepackage{array}
\newcolumntype{M}[1]{>{\centering}m{#1}}

\newcommand{\clearemptydoublepage}{%
\newpage{\pagestyle{empty}\cleardoublepage}}

\theoremstyle{plain}
\newtheorem{theorem}{Theorem}
\numberwithin{theorem}{chapter}
\newtheorem{thm}[theorem]{Theorem} 
\newtheorem{corollary}[theorem]{Corollary}
\newtheorem{lemma}[theorem]{Lemma}
\newtheorem{proposition}[theorem]{Proposition}
\newtheorem{prop}[theorem]{Proposition}

\newtheorem{conjecture}{Conjecture}[chapter] 
\newtheorem{question}[conjecture]{Question}

\newtheorem{definition}{Definition}[chapter]
\newtheorem{defn}[definition]{Definition} 

\theoremstyle{definition}
\newtheorem{example}{Example}

\newtheorem{exercise}{Exercise}[chapter]    
\newtheorem{problem}[exercise]{Problem}

\theoremstyle{remark}
\newtheorem{remark}{Remark}[chapter]

\numberwithin{equation}{chapter}
\numberwithin{figure}{chapter}
\def\LyonsPeres{Lyons:book}

\def\BKS{MR2001m:60016}
\def\BKSfpp{MR2016607} 
\def\SS{SchrammSteif}

\def\RandomTurnHex{MR2309980}
\def\BalancedBoolean{MR2181623}

\def\SchrammSmirnovNoise{SSblacknoise}

\def\GPS{GPS}

\def\Nelson{MR0210416}

\def\Gross{MR0420249}
\def\HPS{MR1465800}

\def\KKL{KKL} 
\def\BKKKL{MR1194785}  
\def\FriedgutKalai{MR1371123} 

\def\Grimm{Grimmett:newbook}

\def\MargulisRusso{MR0472604}

\def\RussoApprox{MR671248}

\def\TalagrandRusso{MR1303654}
\def\TalaPositiveCorr{MR1401897}

\def\OdonnellThesis{ODonnellThesis}
\def\OdonnellBlog{OdonnellBlog}
\def\DecisionTrees{DecisionTrees}

\def\SurveySteif{SurveySteif}
\def\WWperc{arXiv:0710.0856}

\def\SmirnovIsing1{arXiv:0708.0039}
\def\SmirnovPerc{MR1851632}
\def\KestenScaling{MR88k:60174}
\def\SmirnovWerner{MR1879816}
\def\StrictInequalities{MR893131}

\def\JohanssonShapeFluct{MR1737991}

\def\BenaimRossignol{MR2451057}

\def\HaraSladeLace{MR1283177}
\def\MajorityStablest{MajorityStablest}

\def\SchSLE{MR1776084}

\def\LSWoneArm{MR2002k:60204}

\def\NolinKesten{NolinKesten}

\def\FriedgutLowInf{MR1645642}
\def\FriedgutRevisited{MR2034300}
\def\SurveyNoise{SurveyNoise}
\def\NewmanPiza{MR1349159}
\def\PemantlePeres{MR1283187}

\def\Russo81{MR618273}
\def\Kesten1/2{MR575895}
\def\Bollobas1/2{MR2239042}
\def\Triangle{MR1127713}

\def\OdonnellServedio{MR2341918}
\def\FKW{MR2047023}
\def\AS{MR2003f:60003}

\let\qqed=\qed
\def\QED{\qqed\medskip}
\let\qed=\QED
\newcommand{\Prob} {{\mathbb P}}
\newcommand{\R}{\mathbb{R}}
\newcommand{\Q}{\mathbb{Q}}
\newcommand{\C}{\mathbb{C}}
\newcommand{\Z}{\mathbb{Z}}

\def\T{\mathbb T}

\def\Piv{\calP}  
\def\H{\mathbb{H}}

\def\var{\mathop{\mathrm{var}}}

\def\dist{\mathop{\mathrm{dist}}}

\def\SLEkk#1/{$\mathrm{SLE}(#1)$}
\def\SLEr#1/{$\mathrm{SLE(\kappa;#1)}$}
\def\SLEkr#1;#2/{$\mathrm{SLE(#1;#2)}$}
\def\SLEk/{\SLEkk{\kappa}/}
\def\SLEtwo/{\SLEkk2/}

\def\SLEab/{\SLEkr 4; {a/\hco-1}, {b/\hco-1}/}
\def\Var{\mathrm{Var}}
\def\Cov{\mathrm{Cov}}
\def\Ito/{It\^o}
\def \eps {\epsilon}
\def \P {\Prob}
\def\md{\mid}
\def\Bb#1#2{{\def\md{\bigm| }#1\bigl[#2\bigr]}}
\def\BB#1#2{{\def\md{\Bigm| }#1\Bigl[#2\Bigr]}}
\def\Bs#1#2{{\def\md{\mid}#1[#2]}}

\def\BUFF#1#2#3{{#1_{#2} \bigl[ #3 \bigr]}}

\def\SPf{\BUFF{\hat \P}}
\def\SQf{\BUFF{\hat \Q}}

\def\SQ{\hat \Q}
\def\SP{\hat \P}  
\def\SQb{\Bb{\hat \Q}}
\def\SPb{\Bb{\hat \P}}
\def\SEb{\Bb{\hat \E}}

\def\Pb{\Bb\P}
\def\Eb{\Bb\E}

\def\EB{\BB\E}

\def\Es{\Bs\E}

\def \p {{\partial}}
\def \E {{\mathbb E}}

\def\ev#1{{\mathcal{#1}}}
\def \proof {{ \medbreak \noindent {\bf Proof.} }}
\def \proofoutline {{ \medbreak \noindent {\bf Outline of Proof.} }}

\def\sign{\mathrm{sign}}

\def\Spec{\mathscr{S}}
\def\Tree{\mathcal{T}}
\def\Ann{\mathcal{A}}

\def\1{1}

\def\lala(#1,#2){\lambda_{#1,#2}}

\def\nn{[n]}

\def\Ess_#1{\E |\Spec_{f_{#1}}|}

\def\Inf{{\bf I}}
\def\InfV{\mathrm{\bf Inf}}
\def\II{{\bf H}}
\def\Circ{\mathrm{Circ}}
\def\Exc{\mathcal{E}}
\def\reveal{\boldsymbol{\delta}}

\def\calP{{\mathcal P}}

\def\Stab{\mathbb{S}}

\def\D{\mathbb{D}}
\def\MAJ{\boldsymbol{\mathrm{MAJ}}}
\def\DICT{\boldsymbol{\mathrm{DICT}}}
\def\PAR{\boldsymbol{\mathrm{PAR}}}
\def\IT-3-MAJ{\boldsymbol{\mathrm{IT-3-MAJ}}}
\def\CLIQ{\boldsymbol{\mathrm{CLIQ}}}

\def\ni{\noindent}
\def\bi{\begin{itemize}}
\def\ei{\end{itemize}}

\newcommand{\margin}[1]{\marginpar{ \vskip -1cm \textcolor{magenta} {\it #1 }  }}
\newcommand{\note}[2]{ \hskip 2cm  \textcolor{blue}{\large \bf #1 :}   \vline\,\vline \hskip 0.5 cm \parbox{10 cm}{ #2}  }
\newcommand{\CommentBook}[1]{}

\newcommand{\BookChapter}[1]{}

\renewcommand{\margin}[1]{}
\renewcommand{\note}[2]{}{}

\begin{document}


\frontmatter

\thispagestyle{empty}

\begin{center}
\vskip 2.5 cm
\Huge{ \bf Noise sensitivity of Boolean functions and percolation}
\end{center}

\vskip 1 cm
\begin{center}
\Large{Christophe Garban\footnote{ENS Lyon, CNRS} \hskip 1 cm Jeffrey E. Steif\footnote{Chalmers University}}
\end{center}

\vskip 2 cm

\begin{figure}[!htp]
\begin{center}
\includegraphics[width=\textwidth]{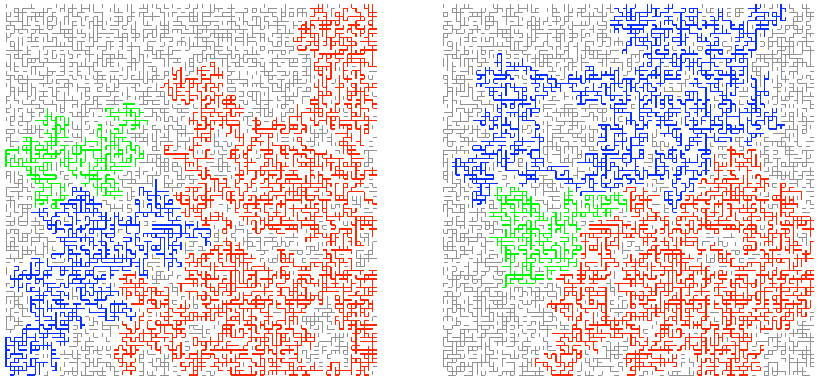}
\end{center}
\end{figure}

\clearemptydoublepage

\mainmatter

\tableofcontents

\newpage
\thispagestyle{empty}

\chapter*{Overview}
\addcontentsline{toc}{chapter}{Overview}

\vskip -0.2 cm
The goal of this set of lectures is to combine two seemingly unrelated topics:
\bi
\item The study of {\bf Boolean functions}, a field particularly active in computer science
\item Some models in statistical physics, mostly {\bf percolation}
\ei

\ni
The link between these two fields can be loosely explained as follows: a percolation configuration
is built out of a collection of i.i.d.\ ``bits'' which determines whether the corresponding edges, sites, or blocks
are present or absent. In that respect, any event concerning percolation can be seen as a Boolean function whose
input is precisely these ``bits''.

Over the last 20 years, mainly thanks to the computer science community,
a very rich structure has emerged concerning the properties of Boolean functions.
The first part of this course will be devoted to a 
description of some of the main achievements in this field.

 In some sense one can say, although this is an exaggeration, that computer scientists are mostly interested in the
{\it stability} or {\it robustness} of Boolean functions. As we will see later in this course, the Boolean functions
which ``encode'' large scale properties of {\bf critical} percolation will turn out to be very {\it sensitive}
to small perturbations. This phenomenon corresponds to what we will call {\bf noise sensitivity}. Hence, the
Boolean functions one wishes to describe here are in some sense {\it orthogonal} to the Boolean functions one encounters, ideally, in computer science.  Remarkably,
it turns out that the tools developed by the computer science community
to capture the properties and stability of Boolean functions are also suitable for the study of noise sensitive functions.
This is why it is worth us first spending some time on the general properties of Boolean functions.

One of the main tools needed to understand properties of Boolean functions is Fourier analysis on the
hypercube. Noise sensitivity will correspond to our Boolean function being of  
``high frequency''  while
stability will correspond to our Boolean function being of  ``low frequency''.
We will apply these ideas to some other models from statistical mechanics as well; namely, first passage percolation and dynamical percolation.

Some of the different topics here can be found (in a more condensed form) in \cite{\SurveyNoise}.

\note{Previously:}{Also, updated versions of these lecture notes will be posted on \url{http://www.umpa.ens-lyon.fr/~cgarban}}

\newpage
\thispagestyle{empty}

\vskip 1 cm

\ni
{\bf \Large Acknowledgements}
\vskip 1 cm

We wish to warmly thank the organizers David Ellwood, Charles Newman, Vladas 
Sidoravicius and Wendelin Werner for inviting us to give  this course at the 
Clay summer school 2010 in Buzios. It was a wonderful experience for us to give 
this set of lectures. We also wish to thank Ragnar Freij who served as a 
very good teaching assistant for this course and for various comments on the
manuscript.

\newpage
\thispagestyle{empty}

\vskip 1 cm

\ni
{\bf \Large Some standard notations}
\vskip 1.5 cm

\ni
In the following table, $f(n)$ and $g(n)$ are any sequences of positive real numbers.
\vskip 1 cm

\begin{center}
\begin{tabular}{| m{3cm} |   m{10 cm} |  }
\hline
 $f(n)\asymp g(n)$
&   
there exists some constant $C>0$ such that 
\[
C^{-1} \le \frac{ f(n)}{g(n)}  \le C\,, \;  \forall n\geq 1 
\] \\ \hline

$f(n) \leq O(g(n))$ 
& there exists some constant $C>0$ such that 
\[
f(n)\leq C g(n)\,, \; \forall n\geq 1
\] \\
\hline

$f(n) \geq \Omega(g(n))$
& there exists some constant $C>0$ such that 
\[
f(n)\geq C g(n)\,, \; \forall n\geq 1
\] \\ \hline

$f(n) = o(g(n))$
& 
\[
\lim_{n\to \infty} \frac {f(n)}{g(n)} =0 
\] \\

\hline

\end{tabular}
\end{center}


\CommentBook{
\chapter*{Introduction/Motivation(s)}
\addcontentsline{toc}{chapter}{Introduction/Motivation}

/////////// 30 minutes ////////////

{\bf Summarize in some sentences the content of the 30 minutes slides talk}

In this short and informal introduction, one shall present some overview of what the course will be about
and what are the main goals one wishes to achieve.

\begin{itemize}
\item Boolean functions naturally arise in Computer science
\item Voting schemes and ``Majority is Stablest'' Theorem stated in an ``intuitive way''.
\item Somehow, we will be interested in functions which are ``orthogonal'' to these stable functions: they are called {\it noise sensitive functions}.
They somehow have a more chaotic behavior.
These functions occur naturally in statistical physics in systems near their {\it critical temperature}.
They are the natural candidates for the description of systems for which {\bf small microscopic perturbations have big macroscopic effects}.

\item {\bf Dynamical percolation movie} (Vincent helped me with doing one). Start explaining how Boolean functions might eventually help explaining these types of 
phenomenon.

\item First Passage Percolation Small Fluctuations: seems to describe a smoothing effect, but in fact is due to a very sensitive dependency on the 
geodesics structure.

\end{itemize}
}



\chapter{Boolean functions and key concepts}\label{ch.BF}

\note{Timing}{//////////// 1 hour ////////////}

\section{Boolean functions}

\begin{definition}
A {\bf Boolean function} is a function from the hypercube 
$\Omega_n:= \{-1,1\}^n$ into either $\{-1 ,1\}$ or  $\{0,1\}$.
\end{definition}

$\Omega_n$ will be endowed with the uniform measure 
$\P=\P^n = (\frac 1 2 \delta_{-1} + \frac 1 2 \delta_1)^{\otimes n}$ and
$\E$ will denote the corresponding expectation.
At various times, $\Omega_n$ will be endowed with the general product measure
$\P_p=\P_p^n = ((1-p) \delta_{-1} + p \delta_1)^{\otimes n}$ but in such
cases the $p$ will be explicit.
$\E_p$ will then denote the corresponding expectations.

An element of $\Omega_n$ will be denoted by either
$\omega$ or $\omega_n$ and its $n$ bits by
$x_1,\ldots, x_n$ so that $\omega= (x_1,\ldots,x_n)$.

Depending on the context, concerning the range,
it might be more pleasant to work 
with one of $\{-1 ,1\}$ or $\{0,1\}$ rather than the other and
at some specific places in these lectures, we will even 
relax the Boolean constraint (i.e.\ taking only two possible values).
In these cases (which will be clearly mentioned), we will 
consider instead real-valued functions $f : \Omega_n \to \R$.
\vskip 0.2 cm

A Boolean function $f$ is canonically 
identified with a subset $A_f$ of $\Omega_n$ via 
$A_f:=\{\omega: f(\omega)=1\}$.

\begin{remark}
Often, Boolean functions are defined on $\{0,1\}^n$ rather than
$\Omega_n= \{-1,1\}^n$. This does not make any fundamental difference 
at all but, as we will see later,
the choice of $\{-1,1\}^n$ turns out to be more
convenient when one wishes to apply Fourier analysis on the hypercube.
\end{remark}

\section{Some Examples}

We begin with a few examples of Boolean functions. Others will appear
throughout this chapter.

\begin{example}[Dictatorship]\label{ex.dict}

\[
\DICT_n(x_1,\ldots,x_n):= x_1
\]

The first bit determines what the outcome is.
\end{example}

\begin{example}[Parity]\label{ex.par}
\[
\PAR_n(x_1,\ldots,x_n):= \prod_{i=1}^n  x_i
\]
This Boolean function tells whether the number of $-1$'s is even or odd.
\end{example}

These two examples are in some sense trivial, 
but they are good to keep in mind since in many 
cases they turn out to be the ``extreme cases'' for properties 
concerning Boolean functions.

The next rather simple Boolean function is of interest in
social choice theory.
\begin{example}[Majority function]\label{ex.maj}
Let $n$ be odd and define
\[
\MAJ_n(x_1,\ldots,x_n):= \sign(\sum_{i=1}^n x_i)\,.
\]
\end{example}

Following are two further examples which will also arise in our
discussions.

\begin{example}[Iterated 3-Majority function]\label{ex.it3maj}
Let $n=3^k$ for some integer $k$. The bits are
indexed by the leaves of a rooted 3-ary tree (so the root has degree
3, the leaves have degree 1 and all others have degree 4) with depth
$k$. One iteratively applies the previous example (with $n=3$) to
obtain values at the vertices at level $k-1$, then level
$k-2$, etc. until the root is assigned a value. The root's value is
then the output of $f$. For example when $k=2$, $f(-1,1,1;1,-1,-1;-1,1,-1)=-1$.
The recursive structure of this Boolean function will enable explicit computations
for various properties of interest.
\end{example}

\begin{example}[Clique containment]\label{ex.clique}
If $r=\binom{n}{2}$ for some
integer $n$, then $\Omega_r$ can be identified with the set of labelled graphs
on $n$ vertices. ($x_i$ is 1 iff the $i$th edge is present.) Recall that a {\bf clique} of size $k$ of a graph $G=(V,E)$ is 
a complete graph on $k$ vertices embedded in $G$. 

Now for any $1\leq k \le \binom{n}{2}=r$, let $\CLIQ_n^k$ be the indicator function of the event that 
the random graph 
$G_\omega$ defined by $\omega\in \Omega_r$ contains a clique of size $k$.
Choosing $k=k_n$ so that this Boolean function is non-degenerate turns out to be a rather delicate issue. The interesting regime is 
near $k_n \approx 2 \log_2(n)$. See the exercises for this 
``tuning'' of  $k=k_n$. 
It turns out that for most values of $n$, the Boolean function $\CLIQ_n^k$ is 
degenerate (i.e.\ has small variance) for all values of $k$. However, there is a 
sequence of $n$ for which there is some $k=k_n$ for which
$\CLIQ_n^k$ is nondegerate.
\end{example}

\section{Pivotality and Influence}

This section contains our first fundamental concepts.
We will abbreviate $\{1,\ldots,n\}$ by $[n]$.

\begin{definition}
Given a Boolean function  $f$ from
$\Omega_n$ into either $\{-1 ,1\}$ or  $\{0,1\}$
and a variable $i\in [n]$, we say that  
{\bf $i$ is pivotal for $f$ for $\omega$}
if $\{f(\omega)\neq f(\omega^i)\}$ where $\omega^i$ is
$\omega$ but flipped in the $i$th coordinate. 
Note that this event is measurable with respect to
$\{x_j\}_{j\neq i}$.
\end{definition}

\begin{definition}
The {\bf pivotal set}, 
$\calP$, for $f$ is the random set of $[n]$ given by
$$
\calP(\omega)=\calP_f(\omega):= \{i\in [n]: \mbox{ $i$ is pivotal for
  $f$ for $\omega$} \}.
$$
\end{definition}

In words, it is the (random) set of bits with the property that if you
flip the bit, then the function output changes.

\begin{definition}
The {\bf influence} of the $i$th bit, $\Inf_i(f)$, is defined by 
$$
\Inf_i(f):=\P(\mbox{ $i$ is pivotal for $f$ }) = \P(i\in \calP).
$$ 
Let also the {\bf influence vector}, $\InfV(f)$, be the collection of all the influences: i.e. $\{\Inf_i(f)\}_{i\in [n]}$.
\end{definition}
In words, the influence of the $i$th bit, $\Inf_i(f)$,
is the probability that, on flipping this bit, the function output changes.

\begin{definition}
The {\bf total influence}, $\Inf(f)$, is defined by
$$
\Inf(f):=\sum_i\Inf_i(f) \, = \| \InfV(f) \|_1 \,\, (=\E(|\calP|)).
$$
\end{definition}

It would now be instructive to go and compute these quantities for
examples \ref{ex.dict}--\ref{ex.maj}. See the exercises.

Later, we will need the last two concepts in the context when our
probability measure is $\P_p$ instead. We give the corresponding
definitions.

\begin{definition}
The {\bf influence vector at level $p$}, $\{\Inf^p_i(f)\}_{i\in [n]}$, is
defined by 
$$
\Inf^p_i(f):=\P_p(\mbox{ $i$ is pivotal for $f$ })= \P_p(i\in \calP).
$$
\end{definition}

\begin{definition}
The {\bf total influence at level $p$}, $\Inf^p(f)$, is defined by
$$
\Inf^p(f):=\sum_i\Inf^p_i(f)  \,\, (=\E_p(|\calP|)).
$$
\end{definition}

It turns out that the total influence has a geometric-combinatorial
interpretation as the size of the so-called edge-boundary of the
corresponding subset of the hypercube. See the exercises.

\begin{remark}
Aside from its natural definition as well as its geometric
interpretation as measuring the edge-boundary of the corresponding subset
of the hypercube (see the exercises), 
the notion of {\it total influence} arises very naturally when 
one studies {\bf sharp thresholds} for {\it monotone functions}
(to be defined in Chapter \ref{ch.ST}). 
Roughly speaking, as we will see in detail in Chapter \ref{ch.ST}, 
for a monotone event $A$, one has that $d\P_p\bigl[ A \bigr]/dp$ 
is the total influence at level $p$ (this is the Margulis-Russo formula). 
This tells us that the speed at which one changes from the event 
$A$  ``almost surely'' {\em not} occurring 
to the case where it ``almost surely'' {\em does} occur 
is very sudden if the Boolean function happens to have a large total influence.
\end{remark}

\section{The Kahn, Kalai, Linial theorem}
This section addresses the following question.
Does there always exist some variable $i$ with (reasonably)
large influence? In other words, for large $n$, what is the 
smallest value (as we vary over Boolean functions) that the largest influence
(as we vary over the different variables) can take on?

Since for the constant function all influences are 
0, and the function which is 1 only if all the bits are 1 has all
influences $1/2^{n-1}$, clearly one wants to deal with functions which
are reasonably balanced (meaning having variances not so close to 0)
or alternatively, obtain lower bounds on
the maximal influence in terms of the variance of the Boolean
function.

The first result in this direction is the following result. A sketch
of the proof is given in the exercises.

\begin{theorem}[Discrete Poincar\'e] \label{th:Poincare}
If $f$ is a Boolean function mapping $\Omega_n$ into $\{-1,1\}$, then
$$
\Var(f) \le \sum_i\Inf_i(f).
$$
It follows that there exists some $i$ such that
$$
\Inf_i(f) \ge \Var(f)/n.
$$
\end{theorem}

This gives a first answer to our question. For reasonably balanced functions,
there is some variable whose influence is at least of order $1/n$. 
{\em Can we find a better ``universal'' lower bound on the maximal influence?}
Note that for Example \ref{ex.maj} all the influences are of order
$1/\sqrt{n}$ (and the variance is 1). 
In terms of our question, this universal lower bound one is looking for should lie somewhere between 
$1/n$ and $1/\sqrt{n}$. The following celebrated result improves by a 
logarithmic factor on the above $\Omega(1/n)$ bound.

\begin{theorem}[\cite{\KKL}] \label{th:KKL1}

There exists a universal $c>0$ such that
if $f$ is a Boolean function mapping $\Omega_n$ into $\{0,1\}$, 
then there exists some $i$ such that
$$
\Inf_i(f) \ge c\Var(f)(\log n)/n.
$$
\end{theorem}

What is remarkable about this theorem is that this 
``logarithmic'' lower bound on the maximal influence turns out to be
{\em sharp}! This is shown by the following example by Ben-Or and Linial.

\begin{example}[Tribes]\label{ex.tribes}
Partition $[n]$ into subsequent blocks of length 
$\log_2(n)-\log_2(\log_2(n))$ with perhaps some leftover
debris. Define $f=f_n$ to be 1 if there exists at least one 
block which contains all 1's,  and 0 otherwise.
\end{example}

It turns out that one can check that the sequence of variances 
stays bounded away from 0 and that all the influences
(including of course those belonging to the debris which are equal to 0) 
are smaller than $c(\log n)/n$ for some $c<\infty$.  See the
exercises for this. Hence the above theorem is indeed sharp.
\note{For the book}{Add reference to Ben Or-Linial paper.}
\medbreak

Our next result tells us that if all the influences are ``small'',
then the total influence is large.

\begin{theorem}[\cite{\KKL}] \label{th:KKL2}
There exists a $c>0$ such that
if $f$ is a Boolean function mapping $\Omega_n$ into $\{0,1\}$
and $\delta:=\max_i \Inf_i(f)$ then
$$
\Inf(f) \ge c \, \Var(f)\log(1/\delta).
$$
Or equivalently,
\[
\| \InfV(f)\|_1 \ge c\, \Var(f) \log{\frac 1 {\| \InfV(f) \|_\infty }}.
\]
\end{theorem}

\medbreak

One can in fact talk about the influence of a set of variables
rather than the influence of a single variable.

\begin{definition} \label{def:setinfluence}
Given $S\subseteq [n]$, the {\bf influence of $S$}, $\Inf_S(f)$, is
defined by
$$
\Inf_S(f):= \P(\mbox{ $f$ is not determined by the bits in $S^c$}).
$$
\end{definition}

It is easy to see that when $S$ is a single bit, this corresponds to
our previous definition.
The following is also proved in \cite{KKL}. We will not indicate the
proof of this result in these lecture notes.

\begin{theorem}[\cite{\KKL}] \label{th:KKL3}
Given a sequence $f_n$ of Boolean functions mapping $\Omega_n$
into $\{0,1\}$
such that $0< \inf_n \E_n(f)\le \sup_n \E_n(f)< 1$ and any sequence
$a_n$ going to $\infty$ arbitrarily slowly, 
then there exists a sequence of sets
$S_n\subseteq [n]$ such that $|S_n|\le a_n n/\log n$ and
$\Inf_{S_n}(f_n)\to 1$ as $n\to\infty$.
\end{theorem}

Theorems \ref{th:KKL1} and \ref{th:KKL2} will be proved in
Chapter \ref{ch.hyper}.

\section{Noise sensitivity and noise stability}

This section introduces our second set of fundamental concepts.

Let $\omega$ be uniformly chosen from $\Omega_n$ and let 
$\omega_\epsilon$ be $\omega$ but with each bit independently
``rerandomized'' with probability $\epsilon$. This
means that each bit, independently of everything
else, rechooses whether it is $1$ or $-1$, each with probability $1/2$.
Note that $\omega_\epsilon$ then has the same distribution as $\omega$.

The following definition is {\it central} for these lecture notes.
Let $m_n$ be an increasing sequence of integers and
let $f_n:\Omega_{m_n} \rightarrow \{\pm 1\}$ or $\{0,1\}$.

\begin{defn}\label{d.NS}
The sequence $\{f_n\}$ is {\bf noise sensitive} if for every
$\epsilon> 0$,
\begin{equation}\label{e.NS}
\lim_{n\to\infty} \E[f_n(\omega)f_n(\omega_\epsilon)]- \E[f_n(\omega)]^2 =0.
\end{equation}
\end{defn}

Since $f_n$ just takes 2 values, this says that 
the random variables $f_n(\omega)$ and $f_n(\omega_\epsilon)$ are 
asymptotically independent for $\epsilon>0$ fixed and $n$ large.
We will see later that (\ref{e.NS}) holds for one value of
$\epsilon\in (0,1)$ if and only if it holds for all such $\eps$.
The following notion captures the opposite situation where
the two events above are close to being the same event 
if $\epsilon$ is small, uniformly in $n$.

\begin{defn}\label{d.STAB}
The sequence $\{f_n\}$ is {\bf noise stable} if 
$$
\lim_{\epsilon\to 0} 
\sup_n \P(f_n(\omega)\neq f_n(\omega_\epsilon))=0.
$$
\end{defn}

It is an easy exercise to check that a sequence
$\{f_n\}$ is both noise sensitive and noise stable 
if and only it is degenerate in the sense that
the sequence of variances $\{\Var(f_n)\}$ goes to 0.
Note also that a sequence of Boolean functions could be neither 
noise sensitive nor noise stable (see the exercises).

It is also an easy exercise to check that 
Example \ref{ex.dict} (dictator) is noise stable and 
Example \ref{ex.par} (parity) is noise sensitive. 
We will see later, when Fourier analysis is brought into
the picture, that these examples are the two opposite extreme
cases.
For the other examples, it turns out that Example \ref{ex.maj} (Majority) is
noise stable, while Examples \ref{ex.it3maj}--\ref{ex.tribes}
are all noise sensitive. See the exercises. In fact, there is a deep
theorem (see \cite{\MajorityStablest}) which says in some sense that, among
all low influence Boolean functions, Example \ref{ex.maj} (Majority) is the
stablest. 

\begin{figure}[!h]
\hskip  -1 cm \includegraphics[width = 1.12\textwidth]{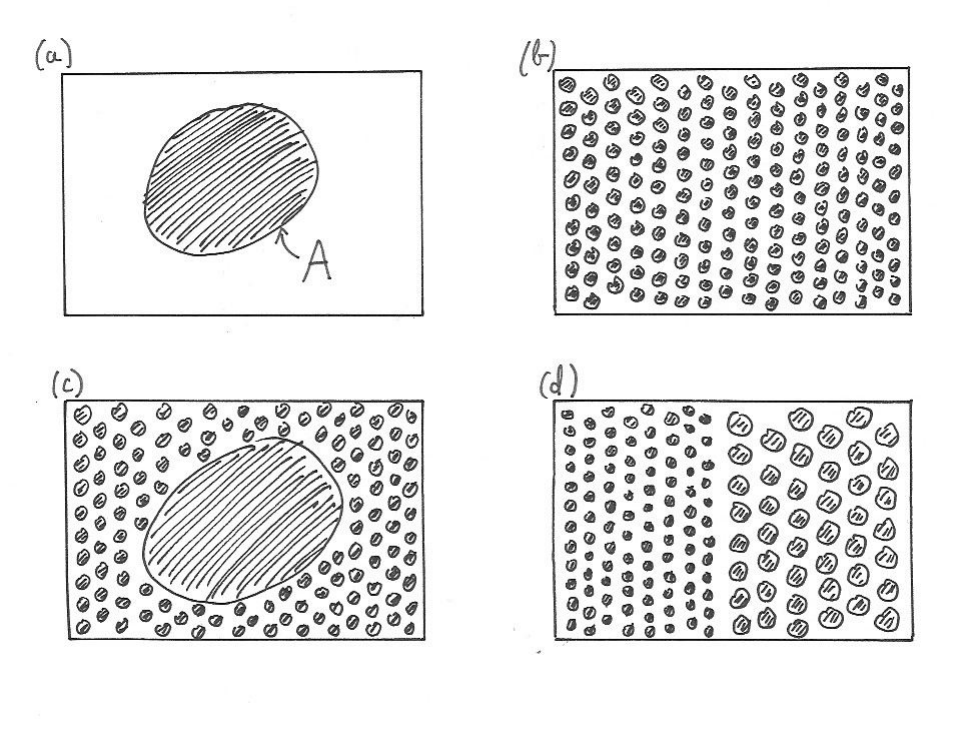}
\caption{Let us consider the following {\em ``experiment''}: take a bounded domain in the plane, say a rectangle, and consider a measurable subset $A$ 
of this domain. What would be an analogue of the above definitions of being {\it noise sensitive} or {\it noise stable} in this case?
Start by sampling a point $x$ uniformly in the domain according to Lebesgue measure. Then let us apply some noise to this position $x$ so that we end up with a new position $x_\epsilon$. One can think of many natural ``noising'' procedures here. For example, let $x_\epsilon$ be a uniform point in the ball of radius $\epsilon$ around $x$, conditioned to remain in the domain. (This is not quite perfect yet since this procedure does not exactly preserve Lebesgue measure, but let's not worry about this.)
The natural analogue of the above definitions is to ask whether $1_A(x)$ and $1_A(x_\epsilon)$ are decorrelated or not.\newline
{\em Question:} According to this analogy, discuss the stability versus 
sensitivity of the sets $A$ sketched in pictures (a) to (d) ? Note that in order to match with definitions \ref{d.NS} and \ref{d.STAB}, one should consider sequences of subsets $\{A_n\}$
instead, since noise sensitivity is an asymptotic notion.
}\label{f.impressionist}
\end{figure}

In Figure \ref{f.impressionist}, we give a slightly 
{\it impressionistic} view of what ``noise sensitivity'' is.

\section{The Benjamini, Kalai and Schramm noise sensitivity theorem}

The following is the main theorem concerning noise sensitivity.

\begin{theorem}[\cite{\BKS}] \label{th:NSmainresult}
If
$$
\lim_n\sum_k \Inf_k(f_n)^2 = 0,
$$
then $\{f_n\}$ is noise sensitive.
\end{theorem}

\begin{remark}
The converse is clearly false as shown by Example \ref{ex.par}. 
However, it turns
out that the converse is true for so-called {\bf monotone functions} (see
the next chapter for the definition of this) as we will see in 
Chapter \ref{ch.FA}.
\end{remark}

This theorem will allow us to conclude noise sensitivity of many of
the examples we have introduced in this first chapter. See the exercises.
This theorem will also be proved in Chapter \ref{ch.hyper}.

\section{Percolation crossings: our final and most important example}

We have saved our most important example to the end. This set of notes
would not be being written if it were not for this example and for
the results that have been proved for it.

Let us consider percolation on $\Z^2$ at the critical point 
$p_c(\Z^2)=1/2$.
(See Chapter \ref{ch.perc} for a fast review on the model.) At this critical point, 
there is no infinite cluster, but somehow clusters are `large' (there are clusters at all scales).
This can be seen 
using duality or with the RSW Theorem \ref{th.RSW}. In order to understand the geometry 
of the critical picture, the following large-scale {\it observables} turn out to be very useful:
Let $\Omega$ be a piecewise smooth domain with two disjoint open arcs $\p_1$ and $\p_2$
on its boundary $\p \Omega$. For each $n\geq 1$, we consider the scaled domain $n \Omega$.
Let $A_n$ be the event that there is an open path in $\omega$ from $n\p_1$ to $n\p_2$ which stays inside 
$n \Omega$. Such events are called {\bf crossing events}. They 
are naturally associated with Boolean functions
whose entries are indexed by the set of edges inside $n \Omega$ (there are $O(n^2)$ such variables).

  For simplicity, let us consider the particular case of rectangle crossings:

\begin{example}[Percolation crossings]\label{ex.perc}
\-
\vskip 0.2 cm

\ni
\begin{minipage}{0.55 \textwidth}
\includegraphics[width=\textwidth]{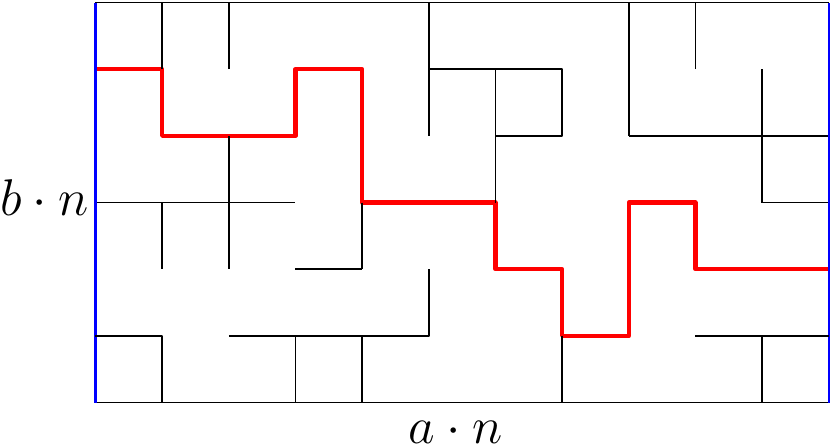}
\end{minipage}
\hspace{0.02\textwidth}
\begin{minipage}{0.46 \textwidth}
Let $a,b>0$ and let us consider the rectangle $[0,a\cdot n]\times[0,b\cdot n]$.
The left to right crossing event corresponds to the Boolean function
 $f_n : \{-1,1\}^{O(1)n^2} \to \{0, 1\}$ defined as follows:
\[
f_n(\omega):= \left\lbrace \begin{array}{ll} 1 & \,   \begin{array}{l} \text{if there is a left-} \\ \text{right crossing} \end{array}  \\ 0 &\, \text{ otherwise}\end{array} \right. 
\]
\end{minipage}

\end{example}

We will later prove that this sequence of Boolean functions $\{f_n\}$ is noise sensitive.
This means that if a percolation configuration $\omega \sim \P_{p_c=1/2}$ is given to us,
one cannot predict anything about the large scale clusters of the slightly perturbed percolation 
configuration $\omega_\eps$ (where only an $\eps$-fraction of the edges have been resampled).

\begin{remark}
The same statement holds for the above more general crossing events (i.e. in $(n\Omega, n\p_1, n\p_2)$).
\end{remark}


\chapter*{Exercise sheet of Chapter \ref{ch.BF}}

\setcounter{exercise}{0}

\begin{exercise}
Determine the pivotal set, the influence vector and the total
influence for Examples \ref{ex.dict}--\ref{ex.maj}.
\end{exercise}

\begin{exercise}
Determine the influence vector for iterated 3-majority and tribes.
\end{exercise}

\begin{exercise}
Show that in Example 6 (tribes) the variances stay bounded away from
0. If the blocks are taken to be of size 
$\log_2 n$ instead, show that the influences would all be of order $1/n$.
Why does this not contradict the KKL Theorem?
\end{exercise}

\begin{exercise}
\label{exercise:edgeboundary}
$\Omega_n$ has a graph structure where two elements are neighbors
if they differ in exactly one location.  The {\bf edge boundary} of 
a subset $A\subseteq \Omega_n$, denoted by $\partial_E(A)$,
is the set of edges where exactly one of
the endpoints is in $A$.  

Show that for any Boolean function,
$\Inf(f)=|\partial_E(A_f)|/2^{n-1}$.
\end{exercise}

\begin{exercise}
Prove Theorem \ref{th:Poincare}.
This is a type of Poincar\'e inequality.
Hint: use the fact that $\Var(f)$ can be written 
$2 \Pb{f(\omega) \neq f(\widetilde{\omega})}$, where 
$\omega,\widetilde{\omega}$ are independent and try to 
``interpolate'' from $\omega$ to $\widetilde{\omega}$.

\note{Further hint}{Let $\omega$ be chosen randomly. Let $\omega_1$ be 
$\omega$ except with the first bit replaced by an independent random
bit. Let $\omega_2$ be $\omega_1$ except with the second bit replaced 
by an independent random bit. Continue until $\omega_n$ is defined.
Consider $\P(f(\omega)\neq f(\omega_n))$ and relate this to the two
sides of the inequality. Note that $\omega$ and $\omega_n$ are independent.}

\end{exercise}

\begin{exercise}
Show that Example \ref{ex.maj} (Majority) is noise stable.
\end{exercise}

\begin{exercise}
Prove that Example \ref{ex.it3maj}
(iterated 3-majority) is noise sensitive 
directly without relying on Theorem \ref{th:NSmainresult}. Hint: use 
the recursive structure of this example in order to show that the criterion of noise sensitivity is satisfied.

\note{Further hint} {Let $a_k=\P(f_{3^k}(\omega)=1=f_{3^k}(\omega_\epsilon))$ and
describe $a_{k+1}$ exactly in terms of $a_k$. (While these terms
depend of course on $\epsilon$, the formula of
$a_{k+1}$ in terms of $a_k$ does not.)
Then use a 'graphical fixed
point argument' to show that $a_k\to 1/4$.}
\end{exercise}

\begin{exercise}\label{ex:noBKS}
Prove that Example \ref{ex.tribes} (tribes) is noise sensitive directly
without using Theorem \ref{th:NSmainresult}. Here there is no recursive 
structure, so a more ``probabilistic'' argument is needed.

\note{A way to do it}{Let $A$ denote the event there at least a tribe, then it's
enough to show that knowing that $A$ occurs, it does 
not help so much for $A^\eps$ to occur. Indeed conditioned on $A$
occurs, let's $T_1, T_2 ,\ldots, T_k$ be the concerned tribes.
Easy to show that $k$ is of order 1. All of these are whp killed. 
The good thing is that the $n/\log{n}- k$ remaining tribes have some 
small amount of negative information, which thus does not help $A^\eps$.}

\note{Another way}{ Use the natural ordering of the blocks. Condition on the first tribe which implies $A$. This particular one is killed whp.
There is negative information for the smaller blocks. The rest of the blocks are untouched, hence unbiased.}
\end{exercise}

\begin{problem}\label{exer.clique}
Recall Example \ref{ex.clique} (clique containment).

\bi 

\item[(a)] Prove that when $k_n=o(n^{1/2})$, $\CLIQ_n^{k_n}$ is asymptotically noise sensitive. Hint: start by obtaining an upper bound on 
the influences (which are identical for each edge) using Exercise
\ref{exercise:edgeboundary}. Conclude by using Theorem \ref{th:NSmainresult}.

\margin{Think of the Total influence as the edge boundary. This gives that individual influences are 
 $O(1) k_n^2/n^2$.}

\item[(b)] \underline{\it Open exercise:} Find a more direct proof of this fact (in the spirit of exercise \ref{ex:noBKS}) which would avoid using Theorem 
\ref{th:NSmainresult}.
\ei

{\it As pointed out after Example \ref{ex.clique}, for most values of $k=k_n$, 
the Boolean function $\CLIQ_n^{k_n}$ becomes degenerate. 
The purpose of the rest of this problem is to determine what the 
interesting regime is where $\CLIQ_n^{k_n}$ has a chance of being 
non-degenerate
(i.e. variance bounded away from 0). The rest of this exercise is somewhat 
tangential to the course.}

\bi

\item[(c)] If $1\leq k \leq \binom{n}{2}=r$, what is the expected number of cliques in $G_\omega$, $\omega\in \Omega_r$ ?

\item[(d)] Explain why there should be at most one choice of $k=k_n$ such that the variance of  $\CLIQ_n^{k_n}$ remains bounded away from 0 ?
(No rigorous proof required.) Describe this choice of $k_n$. Check that it is indeed in the regime $2 \log_2(n)$.

\item[(e)] Note retrospectively that in fact, for any choice of $k=k_n$, $\CLIQ_n^{k_n}$ is noise sensitive.

\ei

\end{problem}

\begin{exercise}
Deduce from Theorem \ref{th:NSmainresult} that both Example 
\ref{ex.it3maj} (iterated 3-majority) and Example 
\ref{ex.tribes} (tribes) are noise sensitive.
\end{exercise}

\begin{exercise}
\medbreak
Give a sequence of Boolean functions which is neither 
noise sensitive nor noise stable.
\end{exercise}

\begin{exercise}
In the sense of Definition \ref{def:setinfluence},
show that for the majority function and for fixed $\eps$, 
any set of size $n^{1/2+\epsilon}$
has influence approaching 1 while any set of size $n^{1/2-\epsilon}$
has influence approaching 0.
\end{exercise}

\begin{exercise}
Show that there exists $c>0$ such that for any Boolean function
$$
\Inf^2_i(f) \ge c\Var^2(f)(\log^2 n)/n
$$
and show that this is sharp up to a constant.
This result is also contained in \cite{\KKL}.
\end{exercise}

\begin{problem}\label{pr:generic}
\medbreak
Do you think a {\it ``generic''} Boolean function would be stable or sensitive? 
Justify your intuition. Show that if $f_n$ was a ``randomly'' chosen
function, then a.s.\  $\{f_n\}$ is noise sensitive.
\end{problem}

\note{Remark}
{JEFF: Here are some relevant notes on this.

One first proves (type of reverse Markov)
that if $X$ is between 0 and $1/2$ with expected value
at least $(1-\delta)/2$ then
$$
P(X\le (1-\ell\delta)/2)\le 1/\ell
$$
for every $\ell$. 

Then one defines (for fixed epsilon)
a rv on Boolean functions (the latter chosen completely
randomly) $X(f)$ which is the probability that $f$ changes when we 
epsilon rerandomize. The mean is easily computed and then one applies the previous
lemma with $\ell=n^2$.  Then Borel Cantelli.

One can also do it via spectrum, described below which is the approach one should
use when the reader is asked to redo this is Chapter 4 I think it is.

It seems as you said that things look very simple.
There are two ways to see that sensitivity is generic
(I did not check the computations but it should roughly be correct)

1/ for any $k\in[1,n]$,
the expectation of $P[|S_{f_\sigma}|= k]$ is exactly $\binom(n,k)/2^n$

(notice the sum over k gives 1 which is good).

this computation shows that things should be very concentrated around
$n/2$ which leads us to the second way

$2/ E_\sigma(E[|S|]) = n/2$

and $E_\sigma (E[|S|^2]) = n^2/4 +n/4$

which gives in expectation $Var(E|S|)= n/4 << (E|S|)^2$
so it is (with high probability) noise sensitive.
Some more details needed but not hard. 
}

\chapter{Percolation in a nutshell}\label{ch.perc}

In order to make these lecture notes as self-contained as possible,
we review various aspects of the percolation model
and give a short summary of the main useful results.
  
For a complete account of percolation, see \cite{\Grimm} and for
a study of the 2-dimensional case, which we are concentrating on here,
see the lecture notes \cite{\WWperc}. 

\section{The model}
Let us briefly start by introducing the model itself. 

We will be concerned mainly with two-dimensional percolation and we will 
focus on two lattices: $\Z^2$ and the triangular lattice $\T$.
(All the results stated for $\Z^2$ in these lecture notes are also valid for 
percolations on ``reasonable'' 2-d translation invariant graphs for which the
RSW Theorem (see the next section) is known to hold at the corresponding 
critical point.)

Let us describe the model on the graph $\Z^2$ which has
$\Z^2$ as its vertex set and edges between vertices having Euclidean distance
1. Let $\E^2$ denote the set of edges of the graph $\Z^2$.
For any $p\in[0,1]$ we define a random subgraph of $\Z^2$ as follows: independently for each edge $e\in \E^2$,
we keep this edge with probability $p$ and remove it with probability $1-p$. Equivalently, this corresponds to defining 
a random configuration $\omega\in \{-1,1\}^{\E^2}$ where, independently for each edge $e\in \E^2$,  we declare the
edge to be open ($\omega(e)=1$) with probability $p$ or  closed $(\omega(e)=-1)$ with probability $1-p$.
The law of the so-defined random subgraph (or configuration) is denoted by $\P_p$.

Percolation is defined similarly on the triangular grid $\T$, except
that on this lattice we will instead  consider {\em site} percolation
(i.e. here we keep each site with probability $p$). 
The sites are the points $\Z+e^{i\pi/3}\Z$ so that neighboring
sites have distance one from each other in the complex plane.

\begin{figure}[!h]
\begin{minipage}[t]{0.48\textwidth}
\centering
\includegraphics[width=\textwidth]{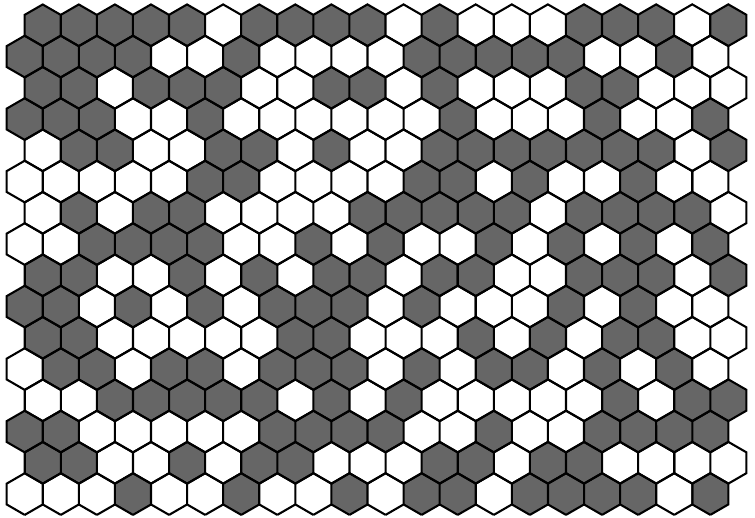}
\end{minipage}
\hspace{0.03\textwidth}
\begin{minipage}[t]{0.48\textwidth}
\centering
\includegraphics[width=\textwidth]{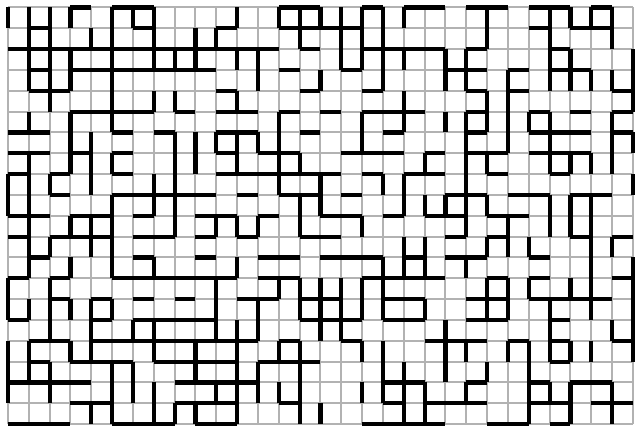}
\end{minipage}
\caption{Pictures (by Oded Schramm) representing two percolation configurations respectively on $\T$ and on $\Z^2$
(both at $p=1/2$). The sites of the triangular grid 
are represented by hexagons.}\label{f.perco}
\end{figure}

\section{Russo-Seymour-Welsh}
We will often rely on the following celebrated result known as the {\bf RSW Theorem}.

\begin{theorem}[RSW]\label{th.RSW} (see \cite{\Grimm}) 
For percolation on $\Z^2$ at $p=1/2$, one has the following property concerning the 
crossing events.
Let $a,b>0$. There exists a constant $c=c(a,b)>0$, such that for any 
$n\geq 1$, if $A_n$ denotes the event that there is a left to right crossing in the rectangle 
$\left([0,a\cdot n]\times [0,b\cdot n]\right) \cap \Z^2$, then 
\[
c< \P_{1/2}\bigl[A_n\bigr]<1-c\,.
\]
In other words, this says that the Boolean functions $f_n$ defined in Example \ref{ex.perc} of Chapter \ref{ch.BF}
are non-degenerate.
\end{theorem}

The same result holds also in the case of site-percolation on $\T$ (also at $p=1/2$). 

\ni
\begin{minipage}{0.55 \textwidth}
The parameter $p=1/2$ plays a very special role for the two models under consideration.
Indeed, there is a natural way to associate to each percolation configuration 
$\omega_p \sim \P_p$ a dual configuration $\omega_{p^*}$ on the so-called dual graph.
In the case of $\Z^2$, its dual graph can be realized as $\Z^2+(\frac 1 2, \frac 1 2)$.
In the case of the triangular lattice, $\T^* = \T$. The figure on the right illustrates 
this {\em duality}
for percolation on $\Z^2$. It is easy to see that in both cases $p^* = 1-p$. 
Hence, at $p=1/2$, our two models 
happen to be {\em self-dual}.
\end{minipage}
\hspace{0.05\textwidth}
\begin{minipage}{0.4\textwidth}
\includegraphics[width=\textwidth]{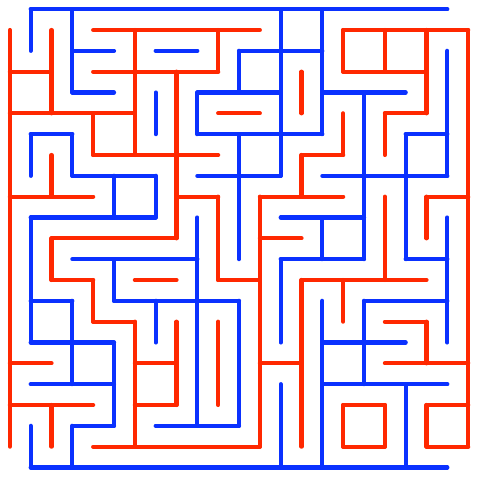}
\end{minipage}

This duality has the following very important consequence. For a domain in 
$\T$ with two specified boundary arcs, there is a 'left-right'
crossing of white hexagons if and only if there is no 'top-bottom'
crossing of black hexagons.

\section{Phase transition}
In percolation theory, one is interested in large scale connectivity properties of the random configuration $\omega=\omega_p$. 
In particular, as one raises the level $p$ above a certain critical parameter $p_c(\Z^2)$, an infinite cluster (almost surely) emerges.
This corresponds to the well-known {\it phase transition} of percolation. By a famous theorem of Kesten this transition takes place at 
$p_c(\Z^2)=\frac 1 2$. On the triangular grid, one also has 
$p_c(\T)=1/2$. 
The event $\{0\overset{\omega}{\longleftrightarrow} \infty\}$
denotes the event that there exists a self-avoiding path from 0
to $\infty$ consisting of open edges.
\vskip 0.2 cm

\noindent\begin{minipage}{0.52 \textwidth}
This phase transition can be measured with the {\it density function} $\theta_{\Z^2}(p): = \P_p (0\overset{\omega}{\longleftrightarrow} \infty)$
which encodes important properties of the 
large scale connectivities of the random configuration $\omega$: it corresponds to the density averaged over the space $\Z^2$ of 
the (almost surely unique) infinite cluster. The shape of the function $\theta_{\Z^2}$ is pictured on the right (notice the infinite derivative at $p_c$).
\end{minipage}
\hspace{0.06\textwidth}
\begin{minipage}{0.4 \textwidth}
\centering
\includegraphics[width=\textwidth]{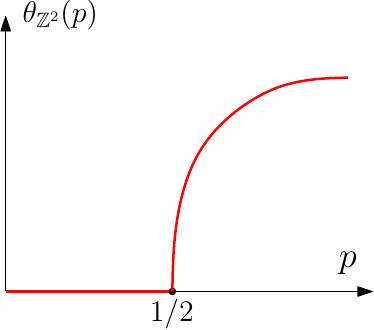}
\end{minipage}
\vskip 0.2 cm

\section{Conformal invariance at criticality and SLE processes}

It has been conjectured for a long time that percolation should be {\it asymptotically} conformally invariant at the critical point.
This should be understood in the same way as the fact that a Brownian
motion (ignoring its time-parametrization) is a 
conformally invariant probabilistic object. One way to picture this conformal invariance
is as follows: consider the `largest' cluster $C_\delta$ surrounding $0$ in 
$\delta \Z^2 \cap \D$ and such that $C_\delta \cap \p \D =\emptyset$. Now consider some other simply connected domain 
$\Omega$ containing 0. Let $\hat C_\delta$ be the largest cluster surrounding $0$ in a critical configuration in $\delta \Z^2 \cap \Omega$
and such that $\hat C_\delta \cap \p \Omega = \emptyset$. Now let $\phi$ be the conformal map from $\D$ to $\Omega$ such that 
$\phi(0)=0$ and $\phi'(0)>0$. Even though the random sets $\phi(C_\delta)$ and $\hat C_\delta$ do not lie on the same lattice, the conformal 
invariance principle claims that when $\delta=o(1)$, these two random clusters are very close in law.
\vskip 0.2 cm

Over the last decade, two major breakthroughs have enabled a much better understanding of the critical regime of percolation:
\bi
\item The invention of the $\mathrm{SLE}$ processes by Oded Schramm(\cite{\SchSLE}).
\item The proof of conformal invariance on $\T$ by Stanislav Smirnov (\cite{\SmirnovPerc}).
\ei

The simplest precise statement concerning conformal invariance is the following. 
Let $\Omega$ be a bounded simply connected domain of the plane and let
$A,B,C$ and $D$ be 4 points on the boundary of $\Omega$ in clockwise order.
Scale the hexagonal lattice $T$ by $1/n$ and perform critical percolation on this 
scaled lattice.
Let $\P(\Omega, A,B,C,D,n)$ denote the probability that in the
$1/n$ scaled hexagonal lattice there is an open path of hexagons in
$\Omega$ going from the boundary of $\Omega$ between $A$ and $B$
to the boundary of $\Omega$ between $C$ and $D$. 

\begin{theorem}(Smirnov, \cite{\SmirnovPerc})

\ni
(i) For all $\Omega$ and $A,B,C$ and $D$ as above,
$$
\P(\Omega, A,B,C,D,\infty):=\lim_{n\to\infty}\P(\Omega, A,B,C,D,n)
$$
exists and is conformally invariant in the sense that
if $f$ is a conformal mapping, then
$\P(\Omega, A,B,C,D,\infty)= \P(f(\Omega), f(A),f(B),f(C),f(D),\infty)$. \\
(ii) If $\Omega$ is an equilateral triangle (with side lengths 1), $A,B$ and $C$ 
the three corner points and $D$ on the line between $C$ and $A$
having distance $x$ from $C$, then the above limiting probability is
$x$.  (Observe, by conformal invariance, that this gives the limiting probability 
for all domains and 4 points.)
\end{theorem}

The first half was conjectured by M.\ Aizenman while J.\ Cardy conjectured the 
limit for the case of rectangles using the four corners. In this case, the formula 
is quite complicated involving hypergeometric functions but
Lennart Carleson realized that this is then equivalent to the simpler formula 
given above in the case of triangles.

Note that, on $\Z^2$ at $p_c=1/2$, proving the conformal invariance is still a challenging  
open problem.

We will not define the $\mathrm{SLE}$ processes in these notes. 
See the lecture notes by Vincent Beffara and references therein.
The illustration below explains how $\mathrm{SLE}$ curves arise naturally 
in the percolation picture.
\vskip 0.2 cm

\noindent\begin{minipage}{0.35 \textwidth}
This celebrated picture (by Oded Schramm)
represents an {\bf exploration path} on the triangular lattice.
This exploration path, which turns right when encountering black hexagons
and left when encountering white ones, asymptotically converges towards $\mathrm{SLE}_6$ (as the mesh size goes to 0).
\end{minipage}
\hspace{0.04\textwidth}
\begin{minipage}{0.6 \textwidth}
\centering
\includegraphics[width=\textwidth]{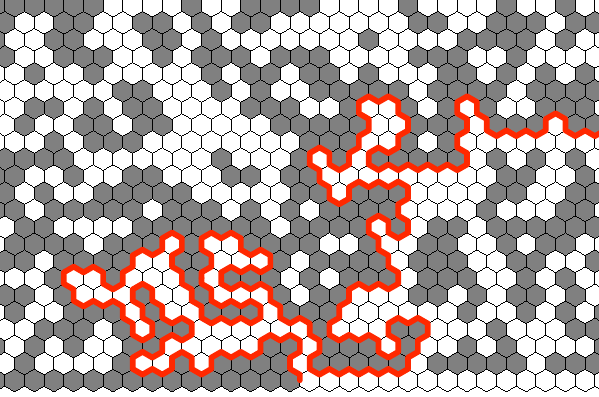}
\end{minipage}
\vskip 0.2 cm

\section{Critical exponents}

The proof of conformal invariance combined with the detailed information given by the $\mathrm{SLE}_6$ process 
enables one to obtain very precise information on the critical and {\it near-critical} behavior of site percolation on $\T$. For instance,
it is known that on the triangular lattice the density function $\theta_{\T}(p)$  has the following behavior near $p_c=1/2$:
\begin{equation}\label{e.5/36}
\theta(p) = (p-1/2)^{5/36 + o(1)}\,,\nonumber
\end{equation}
when $p\to 1/2+$ (see \cite{\WWperc}).
\vskip 0.2 cm

\begin{figure}[!h]
\begin{minipage}[t]{0.48\textwidth}
\centering
\includegraphics[width=\textwidth]{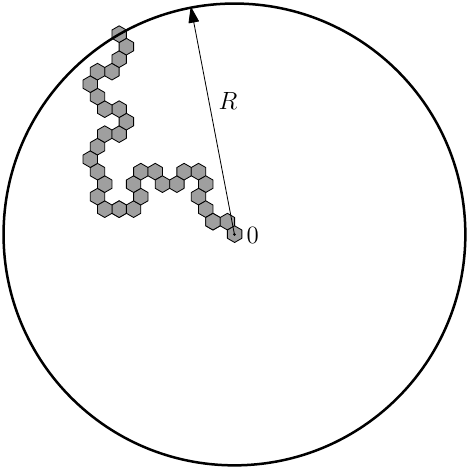}
\end{minipage}
\hspace{0.03\textwidth}
\begin{minipage}[t]{0.48\textwidth}
\centering
\includegraphics[width=\textwidth]{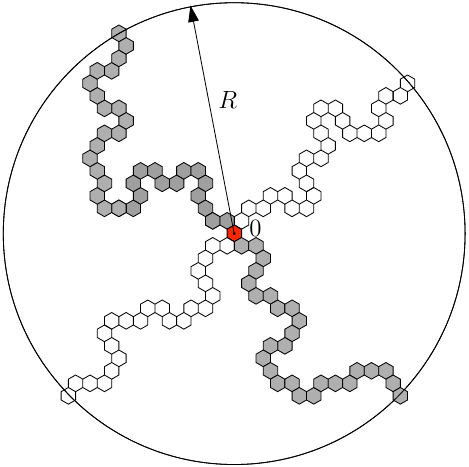}
\end{minipage}
\caption{A realization of the {\it one-arm} event is pictured on the left; the {\it four-arm} event is pictured on the right.}\label{f.armsevents}
\end{figure}

In the rest of these lectures, we will often rely on three 
types of percolation events: namely
the {\it one-arm}, {\it two-arm} and {\it four-arm}
events. They are defined as follows: for any radius $R>1$, let $A_R^1$
be the event that the site 0 is connected to distance $R$ by some open 
path (one-arm). Next, let $A_R^2$ be the event that there are two ``arms'' of
different colors from the site 0 (which itself can be of either color)
to distance $R$ away. Finally, let $A_R^4$ be the event that there are 
four ``arms'' of alternating color from the site 0 (which itself can be of 
either color) to distance $R$ away
(i.e. there are four connected paths, two open, two closed from 0 to
radius $R$ and the closed paths lie  between the open paths). 
See Figure \ref{f.armsevents} for a realization of two of these events.

It was proved in  \cite{\LSWoneArm} that the probability of the one-arm event decays as follows:
\[
\Pb{A_R^1}:= \alpha_1(R) = R^{-\frac 5 {48} +o(1)}\,.
\]
For the two-arms and four-arms events, 
it was proved by Smirnov and Werner in \cite{\SmirnovWerner} that
these probabilities decay as follows:
\[
\Pb{A_R^2}:= \alpha_2(R) = R^{-\frac 1 4 +o(1) }
\]
and 
\[
\Pb{A_R^4}:= \alpha_4(R) = R^{-\frac 5 4 +o(1) }\,.
\]

\begin{remark}
Note the $o(1)$ terms in the above statements 
(which means of course goes to zero as $R\to \infty$). 
Its presence reveals that the above critical exponents are known so far only up to `logarithmic' corrections.
It is conjectured that there are no such `logarithmic' corrections, 
but at the moment one has to deal with their possible existence.
More specifically, it is believed that for the one-arm event,
$$
\alpha_1(R) \asymp R^{-\frac 5 {48}}\,
$$
where $\asymp$ means that the ratio of the two sides is bounded
away from 0 and $\infty$ uniformly in $R$; similarly for the other
arm events.
\end{remark}

The four exponents we encountered concerning $\theta_{\T}$,
$\alpha_1$, $\alpha_2$ and $\alpha_4$ (i.e. $\frac 5 {36}$, 
$\frac 5 {48}$, $\frac 1 4$ and $\frac 5 4$)
are known as {\it critical exponents}.

The {\it four-arm} event is clearly of particular relevance to us in these lectures. 
Indeed, if a point $x$ is in the `bulk' of a domain $(n\Omega, n\p_1, n \p_2)$, 
the event that this point is pivotal for the Left-Right crossing event 
$A_n$ is intimately related to the four-arm event. See Chapter \ref{ch.FE} 
for more details.

\section{Quasi-multiplicativity}\label{s.quasi}

Finally, let us end this overview by a type of scale invariance property of these arm events.
More precisely, it is often convenient to ``divide'' these arm events into different scales. For this purpose, we introduce 
$\alpha_4(r,R)$ (with $r\le R$) to be the probability that the
four-arm event holds from radius $r$ to radius $R$
($\alpha_1(r,R)$, $\alpha_2(r,R)$ and $\alpha_3(r,R)$ are defined analogously).
By independence on disjoint sets, it is clear that if $r_1\le r_2\le r_3$ then 
one has $\alpha_4(r_1,r_3) \le \alpha_4(r_1,r_2) \,
\alpha_4(r_2,r_3)$. A very useful property known as {\bf quasi-multiplicativity}
claims that up to constants, these two expressions are the same (this makes the division into several scales practical).
This property can be stated as follows.
\begin{proposition}[{\bf quasi-multiplicativity}, \cite{\KestenScaling}]\label{pr.quasi}
For any $r_1\le r_2\le r_3$, one has (both for $\Z^2$ and $\T$ percolations)
\begin{equation}
\alpha_4(r_1,r_3) \asymp \alpha_4(r_1,r_2) \,\alpha_4(r_2,r_3)\,.\nonumber
\end{equation}
\end{proposition}
See \cite{\WWperc, \NolinKesten, \SS} 
for more details. Note also that the same property holds for the
one-arm event. However, this is much easier to prove: it is an easy 
consequence of the RSW Theorem \ref{th.RSW} and the so-called
FKG inequality which says that increasing events are positively correlated. 
The reader might consider doing this as an exercise.


\chapter[\;\;\; 
Sharp thresholds and the critical point]
{Sharp thresholds and 
the critical point for 2-d percolation}\label{ch.ST}

\note{Timing}{/////// 45 Minutes ///////}

\section{Monotone functions and the Margulis-Russo formula}

The class of so-called monotone functions plays a very central role in this
subject.

\begin{definition}
A function $f$ is {\bf monotone} if $x\le y$ (meaning
$x_i\le y_i$ for each $i$) implies that $f(x)\le f(y)$.
An event is monotone if its indicator function is monotone.
\end{definition}

Recall that when the underlying variables are independent with 1 having
probability $p$, we let $\P_p$ and $\E_p$ 
denote probabilities and expectations.

It is fairly obvious that for $f$ monotone,
$\E_p(f)$ should be increasing in $p$. The 
Margulis-Russo formula gives us an explicit formula for
this (nonnegative) derivative.

\medskip\noindent
\begin{thm} \label{th:russo} 
Let $A$ be an increasing event in $\Omega_n$. Then
$$
d(\P_p(A))/dp = \sum_{i} \Inf^p_i(A).
$$
\end{thm}

\proof
Let us allow each variable $x_i$ to have its own parameter $p_i$ and
let $\P_{p_1,\ldots,p_n}$ and $\E_{p_1,\ldots,p_n}$ be the corresponding probability measure and expectation.
It suffices to show that
$$
\partial(\P_{(p_1,\ldots,p_n)}(A))/\partial p_i = \Inf^{(p_1,\ldots,p_n)}_i(A)
$$
where the definition of this latter term is clear. WLOG, take $i=1$.
Now
$$
\P_{p_1,\ldots,p_n}(A)=
\P_{p_1,\ldots,p_n}(A\backslash \{1 \in \calP_A\})+
\P_{p_1,\ldots,p_n}(A\cap \{1 \in \calP_A\}).
$$
The event in the first term is measurable with respect to the other variables
and hence the first term does not depend on $p_1$ while the second term is
$$
p_1\P_{p_2,\ldots,p_n}(\{1 \in \calP_A\})
$$
since $A\cap \{1 \in \calP_A\}$ is the event
$\{x_1=1\}\cap \{1 \in \calP_A\}$. 
\qed

\section{KKL away from the uniform measure case}

Recall now Theorem \ref{th:KKL1}. 
For sharp threshold results, one needs lower bounds on the total
influence not just at the special parameter $1/2$ but at all $p$.

The following are the two main results concerning the KKL result for general
$p$ that we will want to have at our disposal. The proofs of these theorems
will be outlined in the exercises in Chapter \ref{ch.hyper}.

\begin{theorem}[\cite{\BKKKL}] \label{th:KKL1generalp}
There exists a universal $c>0$ such that 
for any Boolean function $f$  mapping $\Omega_n$ into $\{0,1\}$
and, for any $p$, there exists some $i$ such that
$$
\Inf_i^p(f) \ge c\Var_p(f)(\log n)/n
$$
\end{theorem}

\begin{theorem}[\cite{\BKKKL}] \label{th:KKL2generalp}
There exists a universal $c>0$ such that
for any Boolean function $f$  mapping $\Omega_n$ into $\{0,1\}$
and for any $p$, 
$$
\Inf^p(f) \ge c\Var_p(f)\log(1/\delta_p)
$$
where $\delta_p:=\max_i \Inf_i^p(f)$.
\end{theorem}

\section{Sharp thresholds in general : the Friedgut-Kalai Theorem}

\begin{theorem}[\cite{\FriedgutKalai}] \label{th:FG}
There exists a $c_1<\infty$ such that for any monotone event $A$ on $n$
variables where all the influences are the same,
if $\P_{p_1}(A)> \epsilon$, then
$$
\P_{p_1+\frac{c_1\log(1/(2\epsilon))}{\log n}}(A)> 1-\epsilon.
$$
\end{theorem}

\begin{remark}
This says that for fixed $\eps$, the probability of $A$ moves from below
$\eps$ to above $1-\eps$ in an interval of $p$ of length of order at
most $1/\log(n)$. The assumption of equal influences holds for example
if the event is invariant under some transitive action, which is often
the case. For example, it holds for Example \ref{ex.it3maj} 
(iterated 3-majority) as well as for any graph property in the context of
the random graphs $G(n,p)$.
\end{remark}

\proof
Theorem \ref{th:KKL1generalp} and all the 
influences being the same tell us that 
$$
\Inf^p(A) \ge c\min\{\P_p(A),1-\P_p(A)\}\log n 
$$
for some $c>0$. Hence Theorem \ref{th:russo} yields
$$
d(\log(\P_p(A)))/dp \ge c\log n 
$$
if $\P_p(A)\le 1/2$. Letting
$p*:=p_1+\frac{\log(1/2\epsilon)}{c\log n}$,
an easy computation (using the fundamental theorem of calculus) yields
$$
\log(P_{p^*}(A))\ge  \log(1/2).
$$
Next, if $\P_p(A)\ge 1/2$, then
$$
d(\log(1-\P_p(A)))/dp \le -c\log n 
$$
from which another application of the fundamental theorem yields
$$
\log(1-P_{p^{**}}(A))\le  -\log(1/\eps)
$$
where $p^{**}:=p^*+\frac{\log(1/2\epsilon)}{c\log n}$. Letting
$c_1=2/c$ gives the result.
\qed

\section{The critical point for percolation for $\Z^2$ and
$\T$ is $\frac 1 2$}

\begin{theorem}[\cite{\Kesten1/2}] \label{th:crit1/2}
$$
p_c(\Z^2) = p_c(\T) = \frac 1 2.
$$
\end{theorem}

\proof
We first show that $\theta(1/2)=0$. 
Let ${\rm Ann}(\ell) := [-3\ell,3\ell] \backslash [-\ell,\ell]$ and
$C_k$ be the event that there is a circuit in ${\rm Ann}(4^k) +1/2$
in the dual lattice around the origin consisting of closed edges.
The $C_k$'s are independent and RSW and FKG show that
for some $c >0$, $\P_{1/2}(C_k)\ge c$ for all $k$.
This gives that 
$\P_{1/2}(C_k \mbox{ infinitely often}) =1$ and hence
$\theta(1/2)=0$.

The next key step is a {\em finite size criterion} which implies
percolation and which is interesting in itself. We outline its proof 
afterwards.

\begin{prop} \label{prop:finitesize}(Finite size criterion)
Let $J_n$ be the event that there is a crossing of a $2n\times (n-2)$ box.
For any $p$, if there exists an $n$ such that
$$
\P_p(J_n) \ge .98,
$$
then a.s.\ there exists an infinite cluster.
\end{prop}

Assume now that $p_c=1/2+\delta$ with $\delta>0$.
Let $I=[1/2,1/2+\delta/2]$. Since $\theta(1/2+\delta/2)=0$, 
it is easy to see that the maximum influence over all variables
and over all $p\in I$ goes to 0 with $n$ since being pivotal implies
the existence of an open path from a neighbor of the given edge to
distance $n/2$ away. Next, by RSW, 
$\inf_n \P_{1/2}(J_n) > 0$. If for all $n$,
$\P_{1/2+\delta/2}(J_n) < .98$, then
Theorems \ref{th:russo} and \ref{th:KKL2generalp} 
would allow us to
conclude that the derivative of $\P_p(J_n)$ goes to $\infty$
uniformly on $I$ as $n\to\infty$, giving a contradiction. Hence 
$\P_{1/2+\delta/2}(J_n) \ge .98$ for some $n$ implying, by 
Proposition \ref{prop:finitesize}, that $\theta(1/2+\delta/2)>0$, a 
contradiction.
\qed

{ \medbreak \noindent {\bf Outline of proof of Proposition \ref{prop:finitesize}.} }

\begin{figure}[!htp]
\begin{center}
\includegraphics[width=0.5\textwidth]{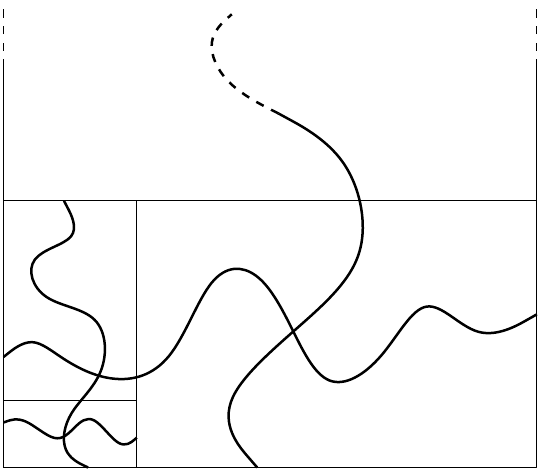}
\end{center}
\caption{}\label{f.finitecriterion}
\end{figure}

The first step is to show that 
for any $p$ and for any $\epsilon \le .02$,
if $\P_p(J_n) \ge 1-\epsilon$, then 
$\P_p(J_{2n}) \ge  1-\epsilon/2$. The idea is that by FKG 
and ``glueing'' one can show
that one can cross a $4n\times (n-2)$ box with probability at least
$1-5\epsilon$ and hence one obtains
that $\P_p(J_{2n}) \ge  1-\epsilon/2$ since, for this event to fail,
it must fail in both the top and bottom halves of the box. It follows
that if we place down a sequence of
(possibly rotated and translated) boxes of sizes
$2^{n+1}\times 2^n$ anywhere, then with probability 1, all but finitely 
many are
crossed. Finally, one can place these boxes down in an intelligent way
such that crossing all but finitely many of them necessarily entails
the existence of an infinite cluster (see Figure \ref{f.finitecriterion}).
\qed

\section{Further discussion}


The Margulis-Russo formula is due independently to 
Margulis \cite{\MargulisRusso} and Russo \cite{\Russo81}.

The idea to use the results from KKL to show that $p_c=1/2$ is due to
Bollob\'as and Riordan (see \cite{\Bollobas1/2}). 
It was understood much earlier
that obtaining a sharp threshold was the key step.
Kesten (see \cite{\Kesten1/2}) showed the necessary sharp threshold 
by obtaining a lower bound on the expected number of pivotals 
in a hands on fashion.
Russo (see \cite{\RussoApprox}) had developed an earlier weaker,
more qualitative, version of KKL
and showed how it also sufficed to show that $p_c=1/2$.

\chapter*{Exercise sheet of Chapter \ref{ch.ST}}

\setcounter{exercise}{0}

\begin{exercise}
Develop an alternative proof of
the Margulis-Russo formula using classical couplings.
\end{exercise}

\begin{exercise}
Study, as best as you can, what the ``threshold windows'' are
(i.e. where and how long does it take to go from a probability of 
order $\eps$ to a probability of order $1-\eps$) in the following examples:
\bi

\item[(a)] for $\DICT_n$
\item[(b)] for $\MAJ_n$
\item[(c)] for the tribes example
\item[(d)] for the iterated majority example.
\ei

Do not rely on \cite{\KKL} type of results, 
but instead do hands-on computations specific to each case.
\end{exercise}

\begin{exercise}
Write out the details of the proof of Proposition \ref{prop:finitesize}.
\end{exercise}

\begin{problem}[{\it What is the ``sharpest'' monotone event ?}]
\label{prob:MajIsSharpest}
Show that among all monotone Boolean functions on $\Omega_n$, $\MAJ_n$
is the one with largest total influence (at $p=1/2$). \\
\ni
Hint: Use the Margulis-Russo formula.

\note{Further hints}{ This comes from Lemma 6.1 in \cite{\FriedgutKalai}. See also Chayes,Chayes, Fisher, Spencer 86.

Let $A$ be some (not necessarily even) increasing event of $\Omega_n$. A good way to look at our problem is 
to think of the total influence via Russo.

\[
\P_p(A) = \sum_k \sum_{|S|=k} p^k q^{n-k} \1_A(S)
\]

\begin{eqnarray}
d(\P_p(A))/dp &= &\sum_{k=1}^n \sum_{|S|=k} p^{k-1}q^{n-k-1} \1_A(S) (k - np) \nonumber \\
&\le & \sum_{k>np} \sum_{|S|=k} ... \1_A(S) (k-np) \nonumber \\
& \le & \sum_{k>np} \sum_{|S|=k} (k- np)  \,,\nonumber
\end{eqnarray}

which is achieved by the ``Majority'' event $B:=\{ S, |S|> np \} $.
\-
Finally, check that for $p=1/2$, this indeed gives asymptotically $\sqrt{\frac 2 {\pi} n} $.
}
\end{problem}

\begin{exercise}
A consequence of Problem \ref{prob:MajIsSharpest} is
that the total influence at $p=1/2$ of any monotone function is at most
$O(\sqrt{n})$. A similar argument shows that for any $p$, there is a constant
$C_p$ so that the total influence at level $p$ of any monotone function is at 
most $C_p\sqrt{n}$. Prove nonetheless that there exists $c>0$ such for
for any $n$, there exists a monotone function $f=f_n$ and a $p=p_n$ so that
the total influence of $f$ at level $p$ is at least $cn$.
\end{exercise}

\begin{exercise}
Find a monotone function $f : \Omega_n \to \{0,1\}$ such that $d(\E_p(f))/dp$ 
is very large at $p=1/2$, but nevertheless there is no sharp threshold for $f$ 
(this means that a large total influence at some value of $p$ is not in general 
a sufficient condition for sharp threshold).
\end{exercise}

\chapter[Fourier analysis of Boolean functions]{Fourier analysis of Boolean functions (first facts)}\label{ch.FA}

\note{Timing}{/////// 45 Minutes ///////}

\section{Discrete Fourier analysis and the energy spectrum}

It turns out that in order to understand and analyze
the concepts previously introduced,
which are in some sense purely probabilistic, a critical tool is 
Fourier analysis on the hypercube.

Recall that we consider our Boolean functions as
functions from the hypercube $\Omega_n:= \{-1,1\}^n$ 
into $\{-1,1\}$ or $\{0,1\}$ where $\Omega_n$ is endowed with the
uniform measure $\P=\P^n = (\frac 1 2 \delta_{-1} + \frac 1 2 \delta_1)^{\otimes n}$.

In order to apply Fourier analysis, the natural setup is to enlarge
our discrete space of Boolean functions and to consider instead the
larger space 
$L^2(\{-1,1\}^n)$ of real-valued functions on $\Omega_n$ endowed with
the inner product:
\begin{eqnarray}
\langle f,g\rangle & :=& \sum_{x_1,\ldots, x_n} 2^{-n} f(x_1,\ldots,x_n)g(x_1,\ldots,x_n) \nonumber\\
&=&\Eb{f g}\; \text{ for all } f,g\, \in L^2(\Omega_n)\,, \nonumber
\end{eqnarray}
where $\E$ denotes expectation with respect to the uniform measure $\P$ on $\Omega_n$.

For any subset $S\subseteq \{1,2\ldots,n\}$, let $\chi_S$ be the
function on $\{-1,1\}^n$ defined for any $x=(x_1,\ldots,x_n)$ by
\begin{equation}
\chi_S(x): = \prod_{i\in S} x_i\,.
\end{equation}
(So $\chi_\emptyset \equiv 1$.)
It is straightforward (check this!) 
to see that this family of $2^n$ functions forms an orthonormal
basis of $L^2(\{-1,1\}^n)$. Thus, any function $f$ on $\Omega_n$ (and a fortiori any Boolean function $f$) can be decomposed as
\[
f = \sum_{S\subseteq \{1,\ldots,n\} } \hat f(S)\, \chi_S,
\]
where $\{\hat f(S)\}_{S\subseteq [n]}$ are the so-called Fourier coefficients of 
$f$. They are also sometimes called the {\bf Fourier-Walsh} coefficients
of $f$ and they satisfy 
\begin{equation}
\hat f(S) := \langle f, \chi_S \rangle = \Eb{ f \chi_S} .\nonumber
\end{equation}
Note that $\hat f(\emptyset)$ is the average $\Eb{f}$ and that we have Parseval's formula which 
states that
\[
\E(f^2) = \sum_{S\subseteq \{1,\ldots,n\} } \hat f^2(S).
\]
As in classical Fourier analysis, if $f$ is some Boolean function, 
its Fourier(-Walsh) coefficients provide information on the
``regularity'' of $f$.  We will sometimes use the term {\it spectrum}
when referring to the set of Fourier coefficients.

\vskip 0.2 cm

Of course one may find many other orthonormal bases for
$L^2(\{-1,1\}^n)$, but there are many situations
for which this particular set of functions 
$\{\chi_S\}_{S\subseteq \{1,\ldots,n \} }$ arises naturally. First of all there is a well-known theory 
of Fourier analysis on groups, a theory which is particularly simple and elegant on Abelian groups (thus including
our special case of $\{-1, 1\}^n$, but also $\R\slash \Z$, $\R$ and so
on).
For Abelian groups, what turns out to be relevant for doing harmonic analysis
is the set $\hat G$ of {\bf characters} of $G$ (i.e. the group homomorphisms from $G$ to $\C^*$). In our case of $G=\{-1,1\}^n$,
the characters are precisely our functions $\chi_S$ indexed by $S\subseteq \{1,\ldots,n\}$ since they satisfy $\chi_S(x\cdot y) = \chi_S(x) \chi_S(y)$.
This background is not however needed and we won't talk in these terms.

These functions also arise naturally if one performs simple random walk on the hypercube (equipped with the Hamming graph structure), since they 
are the eigenfunctions of the corresponding Markov chain
(heat kernel) on $\{-1,1\}^n$.
Last but not least, we will see later in
this chapter that the basis $\{\chi_S\}$ turns out to be particularly 
well adapted to our study of noise sensitivity.

\vskip 0.2 cm
We introduce one more concept here without motivation; it will be very
well motivated later on in the chapter.
\begin{defn}\label{d.EnergySpectrum}
For any real-valued function $f : \Omega_n \to \R$, the
{\bf energy spectrum} $E_f$ is defined by 
\[
E_f(m):= \sum_{|S|=m} \hat f(S)^2, \; \forall \,m \in \{1,\ldots, n\}\,.
\]
\end{defn}

\section{Examples}

First note that, from the Fourier point of view, Dictator and Parity have
simple representations since they are $\chi_1$ and $\chi_{[n]}$
respectively.  Each of the two corresponding energy 
spectra are trivially concentrated on 1 point, namely 1 and $n$.

For Example \ref{ex.maj}, the Majority function, Bernasconi 
explicitly computed the Fourier coefficients 
and when $n$ goes to infinity, one ends up with the following asymptotic
formula for the energy spectrum:
\[
E_{\MAJ_n}(m) = \sum_{|S|=m} \widehat{\MAJ_n}(S)^2 = \left\lbrace \begin{array}{ll} \frac{4}{ \pi \,m 2^m} \binom{m-1}{\frac{m-1}{2}}+O(m/n) & \text{if $m$ is odd}\,, \\ 0 & \text{if $m$ is even}\,. \end{array} \right. 
\]
(The reader may think about why the ``even'' coefficients are 0.)
See \cite{\OdonnellThesis} for a nice overview and references therein
concerning the spectral behavior of the majority function.

\begin{figure}[!h]
\centering
\includegraphics[width=0.5\textwidth]{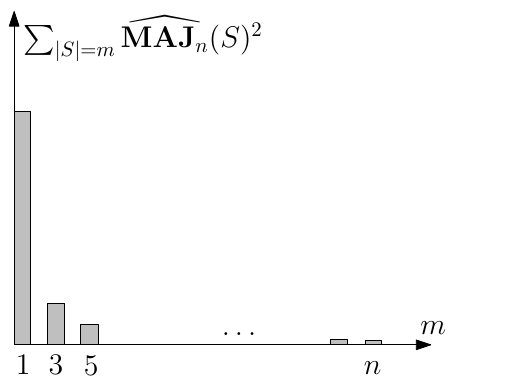}
\caption{Shape of the energy spectrum for the Majority function}\label{f.MAJ}
\end{figure}

Picture \ref{f.MAJ} represents the shape of the energy spectrum of $\MAJ_n$: its spectrum is 
concentrated on low frequencies, which is typical of stable functions.

\section{Noise sensitivity and stability in terms of the energy spectrum}\label{s.NSandEnergySpectrum}

In this section, we describe the concepts of noise sensitivity and
noise stability in terms of the energy spectrum. 

The first step is to note that, given any real-valued function
$f : \Omega_n \to \R$, the
correlation between $f(\omega)$ and $f(\omega_\eps)$
is nicely expressed in terms of the Fourier coefficients of $f$ as 
follows:

\begin{eqnarray}\label{e.correlationFourier}
\Eb{f(\omega) f(\omega_{\eps})} & =& \Eb{\Bigl( \sum_{S_1} \hat f(S_1) \chi_{S_1}(\omega)\Bigr) \Bigl( \sum_{S_2} \hat f(S_2) \chi_{S_2}(\omega_\eps)\Bigr)} \nonumber \\
&=& \sum_{S} \hat f(S)^2\Eb{\chi_S(\omega) \,\chi_S(\omega_\eps)} \nonumber\\
&=& \sum_S \hat f(S)^2 (1-\eps)^{|S|}\,. 
\end{eqnarray}

Moreover, we immediately obtain
\begin{equation}\label{e.mainequation}
\mathrm{Cov}(f(\omega), f(\omega_\eps)) = 
\sum_{m=1}^n E_f(m) (1-\eps)^m.
\end{equation}

Note that either of the last two expressions tell us that
$\mathrm{Cov}(f(\omega), f(\omega_\eps))$ is nonnegative and decreasing in
$\eps$.
Also, we see that the ``level of noise sensitivity'' of a Boolean function is
naturally encoded in its energy spectrum.
It is now an an easy exercise to prove the following proposition.

\begin{proposition}[\cite{\BKS}]\label{pr.NSeq}
A sequence of Boolean functions $f_{n} : \{-1,1\}^{m_n} \to \{0,1\}$
is noise sensitive if and only if, for any $k\geq 1$,
\begin{equation}
\sum_{m=1}^k \, \sum_{|S|=m} \hat f_n(S)^2  = \sum_{m=1}^k E_{f_n}(m) \underset{n\to \infty}{\longrightarrow} 0\,. \nonumber
\end{equation}
Moreover, (\ref{e.NS}) holding does not depend on the value of
$\epsilon\in (0,1)$ chosen.
\end{proposition}

There is a similar spectral description of noise stability which,
given (\ref{e.correlationFourier}), is an easy exercise.

\begin{proposition}[\cite{\BKS}]\label{pr.NStableeq}
A sequence of Boolean functions $f_{n} : \{-1,1\}^{m_n} \to \{0,1\}$
is noise stable if and only if, for any
$\eps>0$, there exists $k$ such that for all $n$,
\begin{equation}
\sum_{m=k}^\infty \, \sum_{|S|=m} \hat f_n(S)^2  = 
\sum_{m=k}^\infty E_{f_n}(m) <\eps. \nonumber
\end{equation}
\end{proposition}

So, as argued in the introduction, a function of ``high frequency'' will 
be sensitive to noise while a function of ``low frequency'' will be
stable. 

\section{Link between the spectrum and influence}\label{s.InfandSpec}

In this section, we relate the notion of influence with that of the
spectrum.

\begin{proposition}\label{pr.inf&spectrum}
If $f:\Omega_n \rightarrow \{0,1\}$, then for all $k$,
\[
\Inf_k(f)= 4 \sum_{S:k\in S} \hat f(S)^2
\]
and
\[
\Inf(f)= 4 \sum_{S} |S| \hat f(S)^2.
\]
\end{proposition}

\proof
If $f:\Omega_n\to \R$, we introduce the functions
\[
\nabla_k f : \left\lbrace
\begin{array}{ccl}
 \Omega_n & \to & \R  \\
\omega & \mapsto & f(\omega) -f(\sigma_k(\omega))
\end{array} \right. \; \text{for all } k\in [n]\,,
\]
where $\sigma_k$ acts on $\Omega_n$ by flipping the $k^{\mathrm{th}}$ bit
(thus $\nabla_k f$ corresponds to a discrete derivative along the
$k^{\mathrm{th}}$ bit).

Observe that
$$
\nabla_k f (\omega)=
\sum_{S\subseteq \{1,\ldots,n\} } \hat f(S)\, 
[\chi_S(\omega)-\chi_S(\sigma_k(\omega))]=
\sum_{S\subseteq \{1,\ldots,n\}, k\in S} 2\hat f(S)\,\chi_S(\omega),
$$
from which it follows that for any $S\subseteq [n]$,
\begin{equation} \label{e.easyfourier}
\widehat{\nabla_k f} (S) = \left\lbrace \begin{array}{ll} 
2 \hat f (S) & \text{if } k\in S \\ 0 & \text{otherwise} \end{array}
\right.
\end{equation}

Clearly, if $f$ maps into $\{0,1\}$, then
$\Inf_k(f) := \| \nabla_k f \|_1$ and since
$\nabla_k f$ takes values in $\{-1,0,1\}$ in this case, we have 
$\| \nabla_k f \|_1= \| \nabla_k f \|_2^2$. Applying Parseval's formula to
$\nabla_k f$ and using (\ref{e.easyfourier}), one obtains the first
statement of the proposition. The second is obtained by summing over 
$k$ and exchanging the order of summation.
\qed

\begin{remark}
If  $f$ maps into $\{-1,1\}$ instead,
then one can easily check that
$\Inf_k(f)=  \sum_{S:k\in S} \hat f(S)^2$
and $\Inf(f)=  \sum_{S} |S| \hat f(S)^2$.
\end{remark}

\section{Monotone functions and their spectrum}

It turns out that for monotone functions, there is an alternative
useful spectral description of the influences.

\begin{proposition}\label{pr.inf&spectrum.monotone}
If $f:\Omega_n \rightarrow \{0,1\}$ is monotone, then for all $k$
\[
\Inf_k(f)= 2\hat f(\{k\})
\] 
If $f$ maps into $\{-1,1\}$ instead,
then one has that $\Inf_k(f)= \hat f(\{k\})$. 
(Observe that Parity shows that the assumption of monotonicity 
is needed here; note also that the proof shows that the weaker result
with $=$ replaced by $\ge$ holds in general.)
\end{proposition}

\proof
We prove only the first statement; the second is proved in
the same way.
$$
\hat f(\{ k\}):= \Eb{f \chi_{\{k\}}} = 
\Eb{f \chi_{\{k\}}I_{\{k\not\in \calP\}}}+
\Eb{f \chi_{\{k\}}I_{\{k\in \calP\}}} 
$$
It is easily seen that the first term is 0 (independent
of whether $f$ is monotone or not)
and the second term is $\frac{\Inf_k(f)}{2}$ due to
monotonicity.
\qed

\begin{remark} \label{rem:BKScomments}
This tells us that, for monotone functions
mapping into $\{-1,1\}$,
the sum in Theorem \ref{th:NSmainresult} is exactly
the total weight of the level 1 Fourier coefficients, that is,
the energy spectrum at 1, $E_f(1)$. (If we map into 
$\{0,1\}$ instead, there is simply an extra irrelevant factor of
4.) So Theorem \ref{th:NSmainresult} 
and Propositions \ref{pr.NSeq} and \ref{pr.inf&spectrum.monotone}
imply that 
for monotone functions, if the energy spectrum at 1 goes to 0,
then this is true for any fixed level. 
In addition, Propositions \ref{pr.NSeq} (with $k=1$) 
and \ref{pr.inf&spectrum.monotone} easily imply that for
monotone functions the converse of Theorem \ref{th:NSmainresult}
holds. 
\end{remark}

Another application of Proposition \ref{pr.inf&spectrum.monotone} gives
a general upper bound for the total influence for monotone functions.

\begin{proposition}\label{pr.monsqrtn}
If $f:\Omega_n \rightarrow \{-1,1\}$ or $\{0,1\}$ 
is monotone, then 
\[
\Inf(f)\le\sqrt{n}.
\] 
\end{proposition}

\proof 
If the image is $\{-1,1\}$, then by Proposition 
\ref{pr.inf&spectrum.monotone}, we have
$$
\Inf(f)=\sum_{k=1}^n \Inf_k(f)=\sum_{k=1}^n \hat f(\{k\}).
$$
By the Cauchy-Schwarz inequality, this is at most 
$(\sum_{k=1}^n \hat f^2(\{k\}))^{1/2}\sqrt{n}$. By Parseval's
formula, the first term is at most 1 and we are done.
If the image is $\{0,1\}$, the above proof can easily be modified
or one can deduce it from the first case
since the total influence of the corresponding
$\pm 1$-valued function is the same.
\qed

\begin{remark}
The above result with some universal $c$ on the right hand side
follows (for odd $n$) from an earlier exercise showing that Majority has the
largest influence together with the known influences for
Majority. However, the above argument yields a more direct proof of the 
$\sqrt{n}$ bound.
\end{remark}

\chapter*{Exercise sheet of chapter \ref{ch.FA}}
\setcounter{exercise}{0}

\begin{exercise}
Prove the discrete Poincar\'e inequality, Theorem \ref{th:Poincare},
using the spectrum.
\end{exercise}

\begin{exercise}
Compute the Fourier coefficients for the indicator function that there
are all 1's.
\end{exercise}

\begin{exercise}
Show that all even size Fourier coefficients for the Majority function
are 0. Can you extend this result to a broader class of Boolean functions? 
\end{exercise}

\begin{exercise}\label{ex.k3}
For the Majority function $\MAJ_n$, find  
the limit (as the number of voters $n$ goes to 
infinity) of the following quantity (total weight of the level-3 Fourier coefficients)
\[
E_{\MAJ_n}(3):= \sum_{|S|=3} \widehat \MAJ_n(S)^2\,.
\]
\end{exercise}

\begin{exercise}
Let $\{f_n\}$ be a sequence of Boolean functions which is noise sensitive
and $\{g_n\}$ be a sequence of Boolean functions which is noise stable.
Show that $\{f_n\}$ and $\{g_n\}$ are asymptotically uncorrelated.
\end{exercise}

\begin{exercise}[Another equivalent definition of noise sensitivity]
 Assume that $\{A_n\}$ is a noise sensitive sequence.
(This of course means that the indicator functions of these events
is a noise sensitive sequence.)

\begin{itemize}
\item[(a)] Show for each $\epsilon >0$, we have
that $\Pb{ \omega_\epsilon\in A_n \bigm| \omega} -\Pb{ A_n} $ approaches 0
in probability. \\
Hint: use the Fourier representation.
\item[(b)] Can you show the above implication without
using the Fourier representation? 
\item[(c)] Discuss if this implication is surprising.
\item[(d)] Show that the condition in part (a) implies that the sequence is
noise sensitive directly without the Fourier representation.
\end{itemize}

\end{exercise}

\begin{exercise}
How does the spectrum of a generic Boolean function look? Use this to 
give an alternative answer to the question asked in problem \ref{pr:generic} 
of Chapter \ref{ch.BF}.
\end{exercise}

\begin{exercise} {\it (Open exercise)}.
For Boolean functions, can one have ANY (reasonable) shape of the
energy spectrum or are there restrictions?
\end{exercise}

For the next exercises, we introduce the following functional  which measures the stability of Boolean functions.
For any Boolean function $f: \Omega_n \to \{-1,1\}$, let 

\[
\mathbb{S}_f \,: \, \eps  \mapsto \Pb{f(\omega)\neq f(\omega_{\eps})} \,.
\]
Obviously, the smaller $\Stab_f$ is, the more stable $f$ is.

\begin{exercise}
Express the functional $\Stab_f$ in terms of the Fourier expansion of $f$.
\end{exercise}

By a {\bf balanced} Boolean function, we mean one which takes its two
possible values each with probability $1/2$.

\begin{exercise}
Among balanced Boolean functions,  
does there exist some function $f^*$ which is ``stablest''  in the sense that for any balanced Boolean function $f$ and any $\eps>0$,
\[
\mathbb{S}_{f^*}(\eps) \le \mathbb{S}_f(\eps)\, ?
\]
If yes, describe the set of these extremal functions and 
prove that these are the only ones.
\end{exercise}

\begin{problem}
In this problem, we wish to understand the \underline{asymptotic shape} of the 
energy spectrum for $\MAJ_n$.

\bi
\item[(a)] Show that for all $\eps\ge 0$, 
$$
\lim_{n\to \infty} \Stab_{\MAJ_n}(\eps)=
\frac 1 2 - \frac{\arcsin(1-\epsilon)}{\pi}=
 \frac{\arccos(1-\epsilon)}{\pi}.
$$
Hint: The relevant limit is easily expressed as the probability
that a certain 2-dimensional Gaussian variable (with a particular correlation
structure) falls in a certain area of the plane. One can write down
the corresponding density function and this probability as an explicit 
integral but
this integral does not seem so easy to evaluate. However, this Gaussian
probability can be computed directly by representing the joint distribution
in terms of two independent Gaussians.

Note that the above limit immediately implies that for $f_n=\MAJ_n$,
$$
\lim_{n\to \infty} \E(f_n(\omega)f_n(\omega_{\eps}))=
\frac{2\arcsin(1-\epsilon)}{\pi}.
$$

\item[(b)] Deduce from (a) and the Taylor expansion for $\arcsin(x)$
the limiting value, as $n\to\infty$ of
$E_{\MAJ_n}(k)= \sum_{|S|=k} \widehat{\MAJ_n}(S)^2$
for all $k\geq 1$. Check that the answer is consistent 
with the values obtained earlier for $k=1$ and $k=3$ (Exercise \ref{ex.k3}).
\ei

\note{A way to do it}{ One needs that if $f_n$ is a sequence of polynomials converging
to a function $f$ such that all the cooeffients of all the 
$f_n$'s are bounded in absolute value by 1 
and $f$ has a power series representation with radius of convergence
at least 1, then for all $k$, the $k$th coefficient of $f_n$ converges to the
$k$th coefficient of $f$. To see this, assume this is not true for the 
$k$th coefficient. Then, since all coefficients are bounded in absolute value by 1,
we can choose a subsequence of the original sequence so that all the coefficients converge.
It is easy to see that this subsequence converges uniformly on $[0,1/2]$ to the function
which has a power series representation with coefficients consisting of the
respective limits of the coefficients and that this has radius of convergence at least 1.
But this limiting function must be $f$ and hence we have two different power series 
representations of $f$ both with radius at least one.  But we know this can't be.}

\end{problem}

\chapter{Hypercontractivity and its applications}\label{ch.hyper}

\note{Timing}{///////// $1.5$ hour  ///////////}

In this lecture, we will prove the main theorems about influences stated in 
Chapter \ref{ch.BF}. As we will see, these proofs rely on techniques imported
from harmonic analysis, in particular {\it hypercontractivity}. 
As we will see later in this chapter and in Chapter \ref{ch.AF}, these types 
of proofs extend to other contexts which will be of interest to us:
noise sensitivity and sub-Gaussian fluctuations.

\note{Remark}{ Try to localize/ understand along the proof why the log term that one finds happens to be OPTIMAL !! (Tribes) }

\section{Heuristics of proofs}
All the subsequent proofs which will be based on {\it hypercontractivity} will have more or less the same flavor. Let us now explain 
in the particular case of Theorem \ref{th:KKL1} what the overall
scheme of the proof is.

\vskip 0.2 cm
Recall that we want to prove that there exists a universal constant $c>0$ such that for any function $f : \Omega_n \to \{0,1\}$,
one of its variables has influence at least  $c \frac {\log n\, \Var(f)}{n}$. 

\vskip 0.2 cm

Let $f$ be a Boolean function. 
Suppose all its influences $\Inf_k(f)$ are ``small'' (this would need to be made quantitative). 
This means that $\nabla_k f$ must have small support.
Using the intuition coming from the Weyl-Heisenberg uncertainty, $\widehat{\nabla_k f}$ 
should then be quite spread out in the sense that most of its spectral mass should be 
concentrated on high frequencies.

This intuition, which is still vague at this point, says that having small influences 
pushes the spectrum of $\nabla_k f$ towards high frequencies. Now, summing up
as we did in Section \ref{s.InfandSpec} of Chapter \ref{ch.FA}, but restricting ourselves 
only to frequencies $S$ of size smaller than some large (well-chosen) $1\ll M \ll n$, one 
easily obtains 

\begin{eqnarray}\label{e.estimate}
\sum_{0<|S|< M} \hat f(S)^2 & \le & 4 \sum_{0<|S| < M} |S| \hat f(S)^2  \nonumber \\
&=& \sum_k \sum_{0<|S|<M} \widehat{\nabla_k f}(S)^2 \nonumber \\
& " \ll " & \sum_k \| \widehat{\nabla_k f} \|_2^2 \nonumber \\
&=& \Inf(f) \,,
\end{eqnarray}
where, in the third line, we used the informal statement that $\widehat{\nabla_k f}$ should be supported on high frequencies if 
$f$ has small influences. Now recall (or observe) that 
\[
\sum_{|S|>0} \hat f(S)^2 = \Var(f)\,.
\]
Therefore, in the above equation~\eqref{e.estimate}, if we are in the case where
a positive fraction of the Fourier mass of $f$ is concentrated below $M$, 
then~\eqref{e.estimate} says that $\Inf(f)$ is much larger than $\Var(f)$. 
In particular, at least one of the influences has to be ``large''.
If, on the other hand, we are in the case where most of the spectral mass of $f$ is supported on frequencies of 
size higher than $M$, then we also obtain that $\Inf(f)$ is large by using the formula: 
\[
\Inf(f)= 4 \sum_S |S| \hat f(S)^2\,.
\]



\begin{remark}
Note that these heuristics suggest that 
there is a subtle balance between $\sum_k \Inf_k(f) = \Inf(f)$ and 
$\sup_k \Inf_k(f)$. Namely, if influences are all small
(i.e.\ $\| \cdot \|_\infty$ is small), then their sum on the other hand has to 
be ``large''. The right balance is exactly quantified by Theorem \ref{th:KKL2}.
\end{remark}

Of course it now remains to convert the above sketch into a proof. The main 
difficulty in the above program is to obtain quantitative spectral information
on functions with values in $\{-1,0,1\}$  knowing that they have small support.
This is done (\cite{\KKL}) using techniques imported from harmonic analysis, namely 
hypercontractivity.

\section{About hypercontractivity}

First, let us state what hypercontractivity corresponds to.
Let $(K_t)_{t\ge 0}$ be the heat kernel on $\R^n$. Hypercontractivity is a statement which quantifies how functions are regularized under 
the heat flow. The statement, which goes back to 
a number of authors,
can be simply stated as follows:

\begin{theorem}[Hypercontractivity]\label{th.hyp}
Consider $\R^n$ with standard Gaussian measure.
If $1<q<2$, there is some $t=t(q)>0$ (which does not depend on the dimension $n$) such that for any $f\in L^q(\R^n)$,
\[
\| K_t \ast f \|_2 \le \| f \|_q \,.
\]
\end{theorem}
The dependence $t=t(q)$ is explicit but will not concern us in the Gaussian case. 
Hypercontractivity is thus a regularization statement: if one starts with some initial 
``rough'' $L^q$ function $f$ outside of $L^2$ and waits long enough ($t(q)$) 
under the heat flow, then we end up being in $L^2$ with a good control on its $L^2$ norm.
We will not prove nor use Theorem \ref{th.hyp}.

This concept has an interesting history as is nicely explained in O'Donnell's
lecture notes (see \cite{\OdonnellBlog}). 
It was originally invented by Nelson 
in \cite{\Nelson} where he needed regularization estimates on Free fields 
(which are the building blocks of quantum field theory) in order to apply these 
in ``constructive field theory''. It was then generalized by Gross in his 
elaboration of logarithmic Sobolev inequalities (\cite{\Gross}),
which is an important tool in analysis. Hypercontractivity is intimately 
related to these Log-Sobolev inequalities and thus has many applications in 
the theory of Semigroups, mixing of Markov chains and other topics.

We now state the result in the case which concerns us, namely the hypercube.
For any $\rho\in[0,1]$, let $T_\rho$ be the following {\bf noise
operator} on the set of functions on the hypercube:
recall from Chapter \ref{ch.BF} that if $\omega\in \Omega_n$, we denote by $\omega_\eps$ an $\eps$-noised configuration of $\omega$.
For any $f : \Omega_n \to \R$, we define $T_{\rho} f: \omega \mapsto  \Eb{f(\omega_{1-\rho}) \bigm| \omega}$. This noise operator
acts in a very simple way on the Fourier coefficients, as the reader can check:
\[
T_\rho : f=\sum_S \hat f(S) \,\chi_S \mapsto \sum_S \rho^{|S|} \hat f(S)\, \chi_S\,.
\]
  
We have the following analogue of Theorem \ref{th.hyp}
\begin{theorem}[Bonami-Gross-Beckner]\label{th.BB}
For any $f : \Omega_n \to \R$ and any $\rho\in [0,1]$,
\[
\| T_\rho f \|_{2} \le \| f\|_{1+\rho^2} \,.
\]
\end{theorem}

The analogy with the classical result \ref{th.hyp} is clear: the heat flow is 
replaced here by the random walk on the hypercube. You can find the proof of 
Theorem \ref{th.BB} in the appendix attached to the present chapter.

\begin{remark}
The term {\it hypercontractive} refers here to the fact that one has
an operator which maps $L^q$ into $L^2$ ($q<2$), which is a contraction.
\end{remark}

\vskip 0.2 cm
\begin{center}
--------------------
\end{center}

Before going into the detailed proof of Theorem \ref{th:KKL1}, let us see why 
Theorem \ref{th.BB} provides us with the type of spectral information we need.
In the above sketch, we assumed that all influences were small. This 
can be written as 
\[
\Inf_k(f) = \| \nabla_k f\|_1 = \|\nabla_k f \|_2^2 \ll 1\,,
\]
for any $k\in [n]$. Now if one applies the hypercontractive estimate to these 
functions $\nabla_k f$ for some fixed $0<\rho<1$, we obtain that

\begin{equation}\label{e.heuristic}
\| T_\rho (\nabla_k f) \|_2 \le \| \nabla_k f \|_{1+\rho^2} = \|\nabla_k f \|_2^{2/(1+\rho^2)} \ll \| \nabla_k f\|_2
\end{equation}
where, for the equality, we used once again that $\nabla_k f \in \{-1,0,1\}$.
After squaring, this gives on the Fourier side,

\[
\sum_{S}  \rho^{2|S|} \widehat{\nabla_k f}(S)^2 \ll \sum_S \widehat{\nabla_k f}(S)^2\,.
\]
This shows (under the assumption that $\Inf_k(f)$ is small)
that the spectrum of $\nabla_k f$ is indeed mostly concentrated on high frequencies.

\begin{remark}
We point out that Theorem \ref{th.BB} in fact tells us that any function with small 
support has its frequencies concentrated on large sets as follows. It is easy to see that given 
any $p< 2$, if a function $h$ on a probability space has very small support, then its 
$L_p$ norm is much smaller than its $L_2$ norm. Using Theorem \ref{th.BB},
we would then have for such a function that
\begin{equation}
\| T_\rho (h) \|_2 \le \| h \|_{1+\rho^2} \ll \| h\|_2\,,\nonumber
\end{equation}
yielding that 
\[
\sum_{S}  \rho^{2|S|} \widehat{h}(S)^2 \ll \sum_S \widehat{h}(S)^2
\]
which can only occur if $h$ has its frequencies concentrated on large sets. 
From this point of view, one also sees that under the small influence assumption,
one did not actually need the third term in (\ref{e.heuristic}) in the above outline.
\end{remark}

\note{For the book}{Maybe one can try to be clearer here: The Boolean assumption here i.e. into $\{ 0, 1\}$ is not 
so important after all, bounded gradient combined with small support is enough}

\section[Proof of the KKL theorems]{Proof of the KKL Theorems on the influences of Boolean functions}

We will start by proving Theorem \ref{th:KKL1}, and then Theorem \ref{th:KKL2}. 
In fact, it turns out that one can recover Theorem \ref{th:KKL1}
directly from Theorem \ref{th:KKL2}; see the exercises.
Nevertheless, since the proof of Theorem \ref{th:KKL1} is slightly
simpler, we start with this one.

\subsection{Proof of Theorem \ref{th:KKL1}}
Let $f : \Omega_n \to \{0,1\}$. Recall that we want to show that there is some $k\in[n]$ such that 

\begin{equation}\label{e.KKL1}
\Inf_k(f) \ge c \Var(f) \frac{\log n}{n}\,,
\end{equation}
for some universal constant $c>0$.

We divide the analysis into the following two cases.

\ni
{\bf Case 1:}

Suppose that there is some $k\in[n]$ such that $\Inf_k(f) \ge n^{-3/4}\, \Var(f)$. Then 
the bound \ref{e.KKL1} is clearly satisfied for a small enough $c>0$.

\ni
{\bf Case 2:}

Now, if $f$ does not belong to the first case, this means that for all $k \in [n]$,
\begin{equation}\label{e.assump}
\Inf_k(f) = \| \nabla_k f \|_2^2 \le \Var(f) n^{-3/4}\,.
\end{equation}

Following the above heuristics, we will show that under this assumption, most of the Fourier spectrum of $f$ is supported on high frequencies.
Let $M\geq 1$, whose value will be chosen later.
We wish to bound from above the bottom part (up to $M$) of the Fourier spectrum of $f$.

\begin{eqnarray}
\sum_{1\le |S| \le M} \hat f(S)^2 &\le &  \sum_{1\le |S| \le M} |S| \hat f(S)^2 \nonumber \\
& \le &  \, 2^{2 M} \sum_{|S| \ge 1} (1/2)^{2|S|} |S| \hat f(S)^2 \nonumber\\
& = & \frac 1 4 \, 2^{2M} \sum_{k} \| T_{1/2}(\nabla_k f) \|_2^2 \,, \nonumber 
\end{eqnarray}
(see Section \ref{s.InfandSpec} of Chapter \ref{ch.FA}).
Now by applying hypercontractivity (Theorem \ref{th.BB}) with $\rho=1/2$ to the above sum, we obtain

\begin{eqnarray}
\sum_{1\le |S| \le M} \hat f(S)^2 &\le & \frac 1 4 \, 2^{2M} \sum_k \| \nabla_k f \|_{5/4}^2 \nonumber \\
& \le &  \, 2^{2M}  \sum_k \Inf_k(f)^{8/5} \nonumber \\
& \le &  \, 2^{2M}\,  n \, \Var(f)^{8/5} n^{\frac{-3} 4 \cdot \frac 8 5}  \nonumber \\
& \le &  \, 2^{2M}\, n^{-1/5} \, \Var(f) \,, \nonumber 
\end{eqnarray}
where we used the assumption \ref{e.assump} and the obvious fact that $\Var(f)^{8/5}\le \Var(f)$ (recall $\Var(f)\le 1$
since $f$ is Boolean). Now with $M:= \lfloor \frac 1 {20} \log_2 n \rfloor$, this gives 

\[
\sum_{1\le |S| \le \frac 1 {20} \log_2 n} \hat f(S)^2 \le  n^{1/10 - 1/5} \, \Var(f) = n^{-1/10}\, \Var(f)\,.
\]

This shows that under our above assumption, most of the Fourier spectrum is concentrated above $\Omega(\log n)$.
We are now ready to conclude:

\begin{eqnarray}
\sup_k \Inf_k(f) \ge \frac {\sum_k \Inf_k(f)} {n} & = & \frac {4 \sum_{|S|\ge 1} |S| \hat f(S)^2} n \nonumber \\
&\ge & \frac 1 n \bigl[  \sum_{|S|>M} |S| \hat f(S)^2  \bigr] \nonumber \\
&\ge & \frac M n   \bigl[ \sum_{|S|>M} \hat f(S)^2 \bigr] \nonumber \\
& = & \frac M n  \bigl[ \Var(f) - \sum_{1\le |S| \le M} \hat f(S)^2 \bigr] \nonumber \\
&\ge & \frac M n \Var(f) \bigl[ 1-  n^{-1/10}\bigr] \nonumber \\
&\ge & c_1 \, \Var(f) \frac{\log n}{n}\,, \nonumber
\end{eqnarray}
with $c_1 = \frac 1 {20 \log 2} (1-2^{-1/10})$. By combining with the constant given in case 1, this completes the proof. \QED

\begin{remark}
We did not try here to optimize the proof in order to find the best possible universal constant $c>0$.
Note though, that even without optimizing at all, the constant we obtain is not that bad.
\end{remark}

\subsection{Proof of Theorem \ref{th:KKL2}}
We now proceed to the proof of the stronger result, Theorem \ref{th:KKL2}, which states that there is a universal constant 
$c>0$ such that for any $f : \Omega_n \to \{0,1\}$,
\[
\| \Inf(f)\| = \| \InfV(f)\|_1 \ge c \, \Var(f) \log{\frac 1 {\| \InfV(f) \|_\infty}}\,.
\]

The strategy is very similar. 
Let $f : \Omega_n \to \{0,1\}$ and let $\delta:= \| \InfV(f) \|_\infty
= \sup_k \Inf_k(f)$. Assume for the moment that $\delta\le 1/1000$.
As in the above proof, we start by bounding the bottom part of the
spectrum up to some integer $M$ (whose 
value will be fixed later). Exactly in the same way as above, one has

\begin{eqnarray}
\sum_{1 \le |S| \le M} \hat f(S)^2  & \le & 2^{2M} \sum_{k} \Inf_k(f)^{8/5} \nonumber \\
& \le &  2^{2M} \delta^{3/5} \sum_k \Inf_k(f) = 2^{2M} \delta^{3/5} \, \Inf(f)\,. \nonumber 
\end{eqnarray}

Now, 
\begin{eqnarray}
\Var(f) = \sum_{|S|\geq 1} \hat f(S)^2 &\le & \sum_{1\le |S| \le M} \hat f(S)^2 + \frac 1 M \sum_{|S|>M} |S|  \hat f(S)^2 \nonumber \\
& \le & \bigl[ 2^{2M} \delta^{3/5} + \frac 1 M \bigr]\, \Inf(f) \,. \nonumber
 \end{eqnarray}

Choose $M:= \frac 3 {10} \log_2(\frac 1 \delta) - \frac 1 2 \log_2
\log_2 (\frac 1 \delta)$. Since $\delta<1/1000$, it is easy to check
that $M\ge \frac 1 {10} \log_2(1/\delta)$ which leads us to
\begin{eqnarray}
\Var(f) &\le&   \left[ \frac 1 {\log_2(1/\delta)} + \frac {10} {\log_2(1/\delta)} \right] \Inf(f) \,\nonumber \\
\end{eqnarray}
which gives 
\[
\Inf(f) = \| \InfV(f)\|_1 \geq \frac 1 {11 \log 2} \, \Var(f) \log \frac 1 {\| \InfV(f) \|_\infty}\,.
\]
This gives us the result for $\delta\le 1/1000$. 

Next the discrete Poincar\'e inequality, which says that $\Inf(f) \geq \Var(f)$,
tells us that the claim is true for $\delta\ge 1/1000$ if we take
$c$ to be $1/\log 1000$. Since this is larger than 
$\frac 1 {11 \log 2}$, we obtain the result with the constant
$c=\frac 1 {11 \log 2}$.
\QED

\section{KKL away from the uniform measure}

In Chapter \ref{ch.ST} (on sharp thresholds), we needed an extension of the above 
KKL Theorems to the $p$-biased measures 
$\P_p = (p \delta_1 + (1-p)\delta_{-1})^{\otimes n}$. These extensions are 
respectively Theorems \ref{th:KKL1generalp} and \ref{th:KKL2generalp}.

A first natural idea in order to extend the above proofs would be to extend the 
hypercontractive estimate (Theorem \ref{th.BB}) to these $p$-biased measures
$\P_p$. This extension of Bonami-Gross-Beckner is possible, but it turns out that 
the control it gives gets worse near the edges ($p$ close to 0 or 1). 
This is problematic since both in Theorems \ref{th:KKL1generalp} 
and \ref{th:KKL2generalp}, we need bounds which are uniform in $p\in [0,1]$.

Hence, one needs a different approach to extend the KKL Theorems.
A nice approach was provided in \cite{\BKKKL}, where they prove the following 
general theorem.

\begin{theorem}[\cite{\BKKKL}]\label{th.BKKKL}
There exists a universal $c>0$ such that for any 
measurable function $f : [0,1]^n \to \{0,1\}$, there exists a variable $k$ such that 
\[
\Inf_k(f) \geq c \, \Var(f) \frac {\log n}{ n}\,.
\]
Here the `continuous' hypercube is endowed with the uniform (Lebesgue) measure and for any $k\in[n]$,
$\Inf_k(f)$ denotes the probability that $f$ is not almost-surely constant on the fiber given by $(x_i)_{i\neq k}$.

In other words, 
\[
\Inf_k(f) = \Pb{ \Var\bigl( f(x_1,\ldots,x_n) \bigm| x_i, i\neq k \bigr) > 0  }\,.
\]
\end{theorem}

It is clear how to obtain Theorem \ref{th:KKL1generalp} from the above theorem. 
If $p\in[0,1]$ and $f : \Omega_n \to \{0,1\}$, consider $\bar f_p : [0,1]^n \to \{0,1\}$ defined by 
\[
\bar f_p(x_1,\ldots, x_n) = f( (1_{x_i<p}-1_{x_i\geq p})_{i\in [n]})\,.
\]

Friedgut noticed in \cite{\FriedgutRevisited} that one can recover 
Theorem \ref{th.BKKKL} from Theorem \ref{th:KKL1generalp}. The first idea is to 
use a symmetrization argument in such a way that the problem reduces to the 
case of monotone functions. Then, the main idea is the approximate the 
uniform measure on $[0,1]$ by the dyadic random variable
\[
X_M : (x_1,\ldots, x_M)\in \{-1,1\}^M \mapsto \sum_{m=1}^M \frac {x_m + 1}{2} 2^{-m}\,.
\]
One can then approximate $f : [0,1]^n \to \{0,1\}$ by the Boolean function $\hat f_M$ defined on $\{-1,1\}^{M\times n}$
by 
\[
\hat f_M(x_1^1,\ldots,x_M^1,\ldots, x_1^n,\ldots, x_M^n) := f(X_M^1,\ldots,X_M^n)\,.
\]

Still (as mentioned in the above heuristics) this proof requires two technical 
steps: a monotonization procedure and an ``approximation'' step (going from 
$f$ to $\hat f_M$). Since in our applications to sharp thresholds we used 
Theorems \ref{th:KKL1generalp} and \ref{th:KKL2generalp} 
only in the case of monotone functions, for the sake of simplicity we will not present the monotonization procedure in these notes.

Furthermore, it turns out that for our specific needs (the applications in 
Chapter \ref{ch.ST}), we do not need to deal with the approximation 
part either. The reason is that for any Boolean function $f$, the function $p\mapsto \Inf_k^p(f)$ is continuous. Hence it is enough 
to obtain uniform bounds on $\Inf_k^p(f)$ for dyadic values of $p$ (i.e. $p \in \{m 2^{-M}\}\cap [0,1]$).

See the exercises for the proof of Theorems \ref{th:KKL1generalp} and \ref{th:KKL2generalp} when $f$ is assumed to be monotone (problem \ref{ex:monotonecase}).

\begin{remark}
We mentioned above that generalizing hypercontractivity would not allow us to obtain uniform bounds (with $p$ taking any value in $[0,1]$)
on the influences. It should be noted though that Talagrand obtained (\cite{\TalagrandRusso}) results similar to Theorems \ref{th:KKL1generalp} and \ref{th:KKL2generalp} by somehow generalizing hypercontractivity, but along a different 
line. Finally, let us point out that 
both Talagrand (\cite{\TalagrandRusso}) and Friedgut 
and Kalai (\cite{\FriedgutKalai}) obtain sharper versions of 
Theorems \ref{th:KKL1generalp} and \ref{th:KKL2generalp} where the 
constant $c=c_p$ in fact improves (i.e. blows up) near the edges.
\end{remark}

\note{Comment}{
So far, we obtained lower bounds on influences at $p=1/2$. 
As we have seen in Chapter \ref{ch.ST}, we need to extend these results to the $p$-biased measures.

One natural approach would be to extend the hypercontractivity result to these measures.
This is not very hard (the tensorization works just as fine) but it turns out that hypercontractivity gets weaker and 
weaker when $p$ approaches the edges of $[0,1]$. See the appendix and the exercises.

There are two ways to generalize properly to the case $p\neq 1/2$.

The first way is through BKKKL
\cite{\BKKKL}. Using their result, one can ``easily'' obtain a uniform bound over the whole interval.

In fact, it seems from \cite{\FriedgutKalai}, section 3, that the proof of BKKKL gives a better estimate saying 
that things get better near the edges (in $(p \log(1/p))^{-1}$).

\ldots

In fact, surprisingly, the bounds are better near the edges, see Talagrand (\cite{\TalagrandRusso})(who extends hypercontractivity to $p \neq 1/2$ but along a different line 
as the naive one above) or \cite{\FriedgutKalai} (who obtain better bounds near the edge only relying on \cite{\BKKKL} !!).

\ldots

As opposed to what Friedgut claims in BKKKL revisited, what Talagrand obtains is in fact stronger.
(Be careful that his $\mu_p(A_i)$ corresponds to our $p\times I_i^p(A)$).
Hence he obtains the nice estimate that 

\[
\frac{\sum_{i} I_i^p(f)} {\log(1/ \| \InfV(f) \|_\infty)}  \gg  \sum_i \frac{I_i^p(f)}{ \log(1/ I_i^p(f))}  \geq c \Var(f) \frac 1 {p(1-p) \log(\frac {1}{p (1-p)})}
\]

}

\section{The noise sensitivity theorem}

In this section, we prove the milestone Theorem \ref{th:NSmainresult} from \cite{\BKS}.
Before recalling what the statement is, let us define the following functional on Boolean functions.
For any $f : \Omega_n \to \{0,1\}$, let 
\[
\II(f):= \sum_k \Inf_k(f)^2 = \| \InfV(f)\|_2^2\,.
\]

Recall the Benjamini-Kalai-Schramm Theorem.

\begin{theorem}[\cite{\BKS}]
Consider a sequence of Boolean functions $f_n : \Omega_{m_n} \to \{0,1\}$.
If
\[
\II(f_n) = \sum_{k=1}^{m_n} \Inf_k(f)^2  \to 0
\]
as $n\to \infty$, then $\{f_n\}_n$ is noise sensitive.
\end{theorem}

We will in fact prove this theorem under a stronger condition, namely that $\II(f_n) \le (m_n)^{-\delta}$ for some exponent 
$\delta>0$.
Without this assumption of ``polynomial decay'' on $\II(f_n)$, the proof is more technical and relies on estimates 
obtained by Talagrand. See the remark at the end of this proof. For our application to the {\it noise sensitivity of percolation}
(see Chapter \ref{ch.FE}),
this stronger assumption will be satisfied and 
hence we stick to this simpler case in these notes.

The assumption of polynomial decay in fact enables us to prove the following more 
quantitative result.
\begin{proposition}[\cite{\BKS}]\label{pr.BKSquantitative}
For any $\delta>0$, there exists a constant $M=M(\delta)>0$ such that if $f_n : \Omega_{m_n} \to \{0,1\}$ is any sequence of Boolean 
functions satisfying
\[
\II(f_n) \le (m_n)^{-\delta}\,,
\]
then 
\[
\sum_{1\le |S| \le M \log{(m_n)}} \widehat{f_n}(S)^2 \to 0 \,.
\]
\end{proposition}

Using Proposition \ref{pr.NSeq},
this proposition obviously implies Theorem \ref{th:NSmainresult} 
when $\II(f_n)$ decays as assumed. Furthermore, this gives a 
quantitative ``logarithmic''
control on the noise sensitivity of such functions.

\proof
The strategy will be very similar to the one used in the KKL Theorems (even 
though the goal is very different). The main difference here is that the 
regularization term $\rho$ used in the hypercontractive estimate must be chosen 
in a more delicate way than in the proofs of KKL results (where we simply took 
$\rho=1/2$).

Let $M>0$ be a constant whose value will be chosen later.

\begin{eqnarray}
\sum_{1\le |S| \le M \log(m_n)} \widehat{f_n}(S)^2 &\le &  4 \sum_{ 1\le |S| \le M \log(m_n)} |S| \widehat{f_n}(S)^2  
=  \sum_k  \sum_{ 1\le |S| \le  M \log(m_n)}  \widehat{\nabla_k f_n}(S)^2 \nonumber \\
& \le & \sum_k  (\frac 1 {\rho^2})^{M \log(m_n)} \| T_\rho(\nabla_k f_n)\|_2^2 \nonumber \\
&\le & \sum_k (\frac 1 {\rho^2})^{M \log(m_n)} \| \nabla_k f_n\|_{1+\rho^2}^2. \nonumber
\end{eqnarray}
by Theorem \ref{th.BB}.

Now, since $f_n$ is Boolean, one has $\|\nabla_k f_n\|_{1+\rho^2} = \|\nabla_k f_n \|_2^{2/(1+\rho^2)}$, hence
\begin{eqnarray}
\sum_{0<|S| < M \log(m_n)} \hat f_n(S)^2 & \le & \rho^{-2 M \log(m_n) }\sum_k \| \nabla_k f_n \|_2^{4/(1+\rho^2)} =  \rho^{-2 M \log(m_n) }\sum_k  \Inf_k(f_n)^{2/(1+\rho^2)}\nonumber\\ 
& \le & \rho^{-2 M \log(m_n)} (m_n)^{\rho^2/(1+\rho^2)} \Bigl( \sum_k \Inf_k(f_n)^2 \Bigr)^{\frac 1 {1+\rho^2}}\text{ (by H\"older)} \nonumber \\
& = & \rho^{-2 M \log(m_n)} (m_n)^{\rho^2/(1+\rho^2)}\, \II(f_n)^{\frac 1 {1+\rho^2}} \nonumber \\
&\le & \rho^{-2 M \log(m_n)} (m_n)^{\frac{\rho^2 - \delta}{1+\rho^2}}\,. \nonumber 
\end{eqnarray}

Now by choosing $\rho\in (0,1)$ close enough to 0, and then by choosing $M=M(\delta)$ small enough, we obtain the desired 
logarithmic noise sensitivity. \QED

\medskip
We now give some indications of the proof of Theorem \ref{th:NSmainresult} in the general case.

Recall that Theorem \ref{th:NSmainresult} is true independently of the speed of 
convergence of $\II(f_n) = \sum_k \Inf_k(f_n)^2$. The proof of this general result 
is a bit more involved than the one we gave here. The main lemma is as follows: 

\begin{lemma}[\cite{\BKS}]\label{lem:talagrandext}
There exist absolute constants $C_k$ 
such that for any monotone Boolean function $f$ and for any $k\geq 2$, 
one has
\[
\sum_{|S|= k} \hat f(S)^2 \le C_k \II(f)\, (-\log \II(f) )^{k-1}\,.
\]
\end{lemma}
This lemma ``mimics'' a result from Talagrand's \cite{\TalaPositiveCorr}. 
Indeed, Proposition 2.3 in \cite{\TalaPositiveCorr} can be translated as 
follows: for any monotone Boolean function $f$, its level-$2$ Fourier weight 
(i.e. $\sum_{|S|=2} \hat f(S)^2$) is bounded by $O(1) \II(f) \log(1/\II(f))$. 
Lemma \ref{lem:talagrandext} obviously implies Theorem \ref{th:NSmainresult} 
in the monotone case, while the general case can be deduced by a 
monotonization procedure. It is worth pointing out that hypercontractivity is 
used in the proof of this lemma.

\note{On influences}
{
After a discussion in the "Capuccino Caf\'e": some remarks/ examples 
came up (some of them being almost the same).

the topic of the discussion was the following: if $f: \Omega_n \to \R$, one can think of at least 3 definitions for the influence
of $k\in [n]$:
\begin{enumerate}
\item $\Inf_k(f) = \Pb{\nabla_k f \neq 0}$ \;  I.e. the size of the support
\item $\Inf_k(f) = \| \nabla_k f \|_1$  \; $L^1$ norm
\item $\Inf_k(f) = \| \nabla_k f \|_2^2$\; an $L^2$ definition of influence (probably related to the $\Eb{\Var( \cdot \md \{ k \} ^c) }$
\end{enumerate}

(Of course for Boolean functions $f: \Omega_n \to \{0,1\}$, they give more or less the same information.)
In some contexts, we've heard the $L^2$ one is more natural. One reason for this is that it makes some well-known inequalities 
HOMOGENOUS, for example {\bf Poincar\'e inequality} works great in this setting: it relates the Variance with some $L^2$ quantities.
{\bf BUT}, for what we are interested in, the definitions based on the size of the support (1 and somewhat 2) are what we want.
Indeed:

\bi
\item Example 1. Assume you want to prove a KKL theorem OR an analog of Lemma \ref{l.torusfpp} in a more general setting.
For example real-valued functions $f: \Omega_n \to [-1,1]$. Then it's easy to see that the $L^2$ definition fails ! Take $f:= (n^{-1/2} \sum x_i)$ truncated at $-1$ and $1$ 
(so that it fits in $[-1,1]$). 

\item In fact, at least in the case of Lemma \ref{l.torusfpp}, one can strengthen the hypothesis $f\to [-1,1]$ into $\nabla_k f \in [-1,1]$. 
which makes the above example even simpler, since no need to truncate anymore. 

\item The generalization of KKL $f: \Omega_n \to [-1,1]$ is not very surprising, it can probably be recovered from the Boolean case.
The generalization to bounded gradient (for example $\nabla_k f \in [-1,1]$) seems more interesting, but that's basically Lemma \ref{l.torusfpp}.
Note that in both cases, with the definition of the support, the inequality on one side remains unchanged, so KKL for example remains true if one scale DOWN the function.
But not if we scale Up. With an $L^1$ definition, it would maybe work slightly better, but not enough yet (not yet homogenous).

\ei
}

\chapter*{Appendix: proof of hypercontractivity}
\addcontentsline{toc}{chapter}{Appendix on Bonami-Gross-Beckner}
\setcounter{exercise}{0}
\setcounter{remark}{0}
\setcounter{section}{0}

The purpose of this appendix is to show that we are not using a 
giant ``hammer'' but rather that this needed inequality arising from 
Fourier analysis is understandable from first principles.
In fact, historically, the proof by Gross of the Gaussian case first 
looked at the case of the hypercube and 
so we have the tools to obtain the Gaussian case should we want to.
Before starting the proof, observe that for $\rho=0$ (where $0^0$ is
defined to be 1), this simply reduces to $|\int f|\le \int |f|$.

\bigskip

{\medbreak \noindent {\bf Proof of Theorem \ref{th.BB}.}}

\subsection{Tensorization}
In this first subsection, we show that it is sufficient,
via a tensorization procedure, that
the result holds for $n=1$ in order for us to be able to
conclude by induction the result for all $n$.

The key step of the argument is the following lemma.

\begin{lemma}\label{lem:BBlemma}
Let $q\ge p\ge 1$, $(\Omega_1,\mu_1),(\Omega_2,\mu_2)$ be two finite probability
spaces, $K_i:\Omega_i\times\Omega_i\to \R$ and assume that for $i=1,2$
$$
\|T_i(f)\|_{L_q(\Omega_i,\mu_i)}\le\|f\|_{L_p(\Omega_i,\mu_i)}
$$
where 
$T_i(f)(x):=\int_{\Omega_i} f(y)K_i(x,y) d\mu_i(y)$. Then
$$
\|T_1\otimes T_2(f)\|_{L_q((\Omega_1,\mu_1)\times (\Omega_2,\mu_2))}
\le\|f\|_{L_p((\Omega_1,\mu_1)\times (\Omega_2,\mu_2))}
$$
where
$T_1\otimes T_2(f)(x_1,x_2):=\int_{\Omega_1\times \Omega_2} 
f(y_1,y_2)K_1(x_1,y_1) K_2(x_2,y_2) 
d\mu_1(y_1)\times d\mu_2(y_2)$.
\end{lemma}

\proof
One first needs to recall Minkowski's inequality for integrals, which
states that, for $g\ge 0$ and $r\in [1,\infty)$, we have

$$
\left(
\int \left(
\int g(x,y) d\nu(y)
\right)^r d\mu(x)
\right)^{1/r} \le
\int \left(
\int g(x,y)^r d\mu(x)
\right)^{1/r}
d\nu(y).
$$
(Note that when $\nu$ consists of 2 point masses each of size 1, then this reduces
to the usual Minkowski inequality.)

One can think of $T_1$ acting on functions of both variables by
leaving the second variable untouched and analogously for $T_2$. It is
then easy to check that $T_1\otimes T_2=T_1\circ T_2$. By
thinking of $x_2$ as fixed, our assumption on $T_1$ yields
$$
\|T_1\otimes T_2(f)\|_{L_q((\Omega_1,\mu_1)\times
(\Omega_2,\mu_2))}^q\le 
\int_{\Omega_2} 
\left(\int_{\Omega_1}|T_2(f)|^p d\mu_1(x_1) \right)^{q/p} d\mu_2(x_2).
$$
(It might be helpful here to think of $T_2(f)(x_1,x_2)$ as 
a function $g^{x_2}(x_1)$ where $x_2$ is fixed).

Applying Minkowski's integral inequality to $|T_2(f)|^p$ with $r=q/p$,
this in turn is at most
$$
\left[
\int_{\Omega_1}
\left(\int_{\Omega_2}
|T_2(f)|^q 
d\mu_2(x_2)\right)^{p/q} 
d\mu_1(x_1) 
\right]^{q/p}.
$$
Fixing now the $x_1$ variable and applying our assumption on $T_2$
gives that this is at most 
$\|f\|_{L_p((\Omega_1,\mu_1)\times (\Omega_2,\mu_2))}^q$, as desired.
\qed

The next key observation, easily obtained by expanding and
interchanging of summation, is
that our operator $T_\rho$ acting on functions on $\Omega_n$
corresponds to an operator of the type dealt with in the previous
lemma with $K(x,y)$ being
$$
 \sum_{S\subseteq \{1,\ldots,n\} }\rho^{|S|} \chi_S(x)\chi_S(y).
$$
In addition, it is easily checked that the function $K$ for the
$\Omega_n$ is simply an $n$-fold product of the function for the
$n=1$ case.

Assuming the result for the case $n=1$, Lemma \ref{lem:BBlemma} 
and the above observations
allow us to conclude by induction the result for all $n$.

\subsection{The $n=1$ case}

We now establish the case $n=1$. We abbreviate $T_\rho$ by $T$.

Since $f(x)=(f(-1)+f(1))/2  + (f(1)-f(-1))/2\,\,\,x$, we have
$Tf(x)=(f(-1)+f(1))/2  + \rho(f(1)-f(-1))/2\,\,\,x$. Denoting 
$(f(-1)+f(1))/2$ by $a$ and $(f(1)-f(-1))/2$ by $b$, it suffices to
show that for all $a$ and $b$, we have
$$
(a^2+\rho^2b^2)^{(1+\rho^2)/2}\le
\frac{|a+b|^{1+\rho^2}+|a-b|^{1+\rho^2}}{2}.
$$

Using $\rho\in [0,1]$, the case $a=0$ is immediate. For the case,
$a\neq 0$, it is clear we can assume $a>0$. Dividing both sides by 
$a^{1+\rho^2}$, we need to show that
\begin{equation}\label{e.BB}
(1+\rho^2y^2)^{(1+\rho^2)/2}\le
\frac{|1+y|^{1+\rho^2}+|1-y|^{1+\rho^2}}{2}
\end{equation}
for all $y$ and clearly it suffices to assume $y\ge 0$.

We first do the case that $y\in [0,1)$. By the generalized Binomial
formula, the right hand side of (\ref{e.BB}) is
$$
\frac{1}{2}\left[\sum_{k=0}^\infty \binom{1+\rho^2}{k}y^k+
\sum_{k=0}^\infty \binom{1+\rho^2}{k}(-y)^k\right]=
\sum_{k=0}^\infty \binom{1+\rho^2}{2k}y^{2k}.
$$

For the left hand side of (\ref{e.BB}), we first note the following.
For $0<\lambda<1$, a simple calculation shows that the function
$g(x)=(1+x)^\lambda-1-\lambda x$ has a negative derivative on
$[0,\infty)$ and hence $g(x)\le 0$ on $[0,\infty)$. 

This yields that the left hand side of (\ref{e.BB}) is at most 
$$
1+\left(\frac{1+\rho^2}{2}\right)\rho^2y^2
$$
which is precisely the first two terms of the right hand side of (\ref{e.BB}).
On the other hand, the binomial coefficients appearing in the other
terms are nonnegative, since in the numerator there are an even number
of terms with the first two terms being positive and all the other
terms being negative. This verifies the desired inequality for 
$y\in [0,1)$.

The case $y=1$ for (\ref{e.BB}) follows by continuity.

For $y>1$, we let $z=1/y$ and note, by multiplying both sides
of (\ref{e.BB}) by $z^{1+\rho^2}$, we need to show
\begin{equation}\label{e.BBsecond}
(z^2+\rho^2)^{(1+\rho^2)/2}\le
\frac{|1+z|^{1+\rho^2}+|1-z|^{1+\rho^2}}{2}.
\end{equation}

Now, expanding $(1-z^2)(1-\rho^2)$, one sees that
$z^2+\rho^2\le 1+z^2\rho^2$ and hence the desired inequality follows
precisely from (\ref{e.BB}) for the case $y\in (0,1)$ 
already proved. This completes the $n=1$ case and thereby the proof.
\qed

\chapter*{Exercise sheet of Chapter \ref{ch.hyper}}
\setcounter{exercise}{0}

\begin{exercise}\label{ex.KKL21}
Find a direct proof that Theorem \ref{th:KKL2} implies Theorem \ref{th:KKL1}.
\end{exercise}

\begin{exercise}
Is it true that the smaller the influences are, the more noise sensitive the function is?
\end{exercise}

\note{Exercise ?}
{
Can we find an example in the Markov Chain setting 
where hypercontractivity is used to prove a mixing time. Fill/Aldous lecture notes ?
}

\note{Used to be an exercise}
{
Why hypercontractivity gets worse near the edges ? To see this, fix $q<2$ and show that in order the operator $T_\rho$
to be hypercontractive from $L^q(d\P_p)$ to $L^2(d\P_p)$, the noise parameter $\rho=\rho(p)$ needs to be chosen small enough.
In other words, near the edges, more regularization is needed. Hint: To show this, find an adequate test function $f$ on which you may test 
the hypercontractive estimate $\| T_\rho f \|_{L^2(p)} \le \| f \|_{L^q(p)}$.

\ni
The problem is that this natural operator $T_\rho$ which re-randomizes each bit with probability $(1-\rho)$ provides less regularization as well, since when $p$ is near the edges, when a bit is resampled, with high probability it remains unchanged.
}

\note{Remark}{The hypercontractive constants are known in this case thanks 
to a theorem by Oleszkiewicz. See Mossel lecture's notes Oct 6.}

\begin{exercise}
Prove that Theorem \ref{th.BKKKL} indeed implies Theorem \ref{th:KKL1generalp}.\\ 
Hint: use the natural projection.
\end{exercise}

\begin{problem}\label{ex:monotonecase}
In this problem, we prove Theorems \ref{th:KKL1generalp} and \ref{th:KKL2generalp} for the monotone case.
\vskip 0.3 cm

\begin{enumerate}
\item
Show that Theorem \ref{th:KKL2generalp}  implies \ref{th:KKL1generalp} and hence one needs to prove only
Theorem \ref{th:KKL2generalp} (This is the basically the same as Exercise \ref{ex.KKL21}).

\item 
Show that it suffices to prove the result when $p=k/2^\ell$
for integers $k$ and $\ell$.

\item 
Let $\Pi:\{0,1\}^\ell\to \{0,1/2^\ell, \ldots,
(2^\ell-1)/2^\ell\}$ by
$\Pi(x_1,\ldots, x_\ell)=\sum_{i=1}^\ell x_i/2^i$. Observe that
if $x$ is uniform, then $\Pi(x)$ is uniform on its range and
that $\P(\Pi(x)\ge i/2^\ell)=(2^\ell-i)/2^\ell$.

\item 
Define $g:\{0,1\}^\ell\to \{0,1\}$ by
$g(x_1,\ldots,x_\ell):=I_{\{\Pi(x)\ge 1-p\}}$. Note that $\P(g(x)=1)=p$.

\item
Define $\tilde{f}:\{0,1\}^{n\ell}\to \{0,1\}$ by
$$
\tilde{f}(x^1_1,\ldots,x^1_\ell,x^2_1,\ldots,x^2_\ell,\ldots,
x^n_1,\ldots,x^n_\ell)=
f(g(x^1_1,\ldots,x^1_\ell),g(x^2_1,\ldots,x^2_\ell),\ldots,
g(x^n_1,\ldots,x^n_\ell)).
$$
Observe that $\tilde{f}$ (defined on
$(\{0,1\}^{n\ell},\pi_{1/2})$) and
$f$ (defined on $(\{0,1\}^{n},\pi_{p})$) have the same distribution
and hence the same variance.

\item
Show (or observe) that
$\Inf_{(r,j)}(\tilde{f})\le \Inf^p_{r}(f)$ for each $r=1,\ldots,n$ and
$j=1,\ldots,\ell$. Deduce from
Theorem \ref{th:KKL2} that
$$
\sum_{r,j}\Inf_{(r,j)}(\tilde{f})\ge c\Var(f)\log(1/\delta_p)
$$
where $\delta_p:=\max_i \Inf^p_i(f)$
where $c$ comes from Theorem \ref{th:KKL2}.

\item
 (Key step). Show that for each $r=1,\ldots,n$ and
$j=1,\ldots,\ell$,
$$
\Inf_{(r,j)}(\tilde{f})\le \Inf^p_{r}(f)/2^{j-1}.
$$

\item  Combine parts 6 and 7 to complete the proof.

\end{enumerate}

\end{problem}

\chapter[\; First evidence of noise sensitivity of percolation]{First evidence of noise sensitivity of percolation}\label{ch.FE}

\note{Timing}{/////////// 45 minutes ($\approx$) //////////////// \vskip 1 cm}

In this lecture, our goal is to collect some of the facts and theorems we have seen so far in order 
to conclude that percolation crossings are indeed noise sensitive. 
Recall from the ``BKS'' Theorem (Theorem \ref{th:NSmainresult}) 
that it is enough for this purpose to prove that influences are ``small'' in the sense that $\sum_k \Inf_k(f_n)^2$ goes to zero.

In the first section, we will deal with a careful study of influences in the case of percolation crossings
on the triangular lattice.
Then, we will treat the case of $\Z^2$, where conformal invariance is not known. Finally, we will
speculate to what ``extent'' percolation is noise sensitive.
\vskip 0.3 cm

{\it 
This whole chapter should be considered somewhat of a ``pause'' in our program,
where we take the time to summarize what we have achieved so far in our understanding of the noise sensitivity of percolation, 
and what remains to be done if one wishes to prove things such as the existence of exceptional times in dynamical percolation. 
}

\section[Influences of crossing events]{Bounds on influences for crossing events in critical percolation on the triangular lattice}

\subsection{Setup}

Fix $a,b > 0$, let us consider some rectangle $[0,a\cdot n]\times [0, b\cdot n]$, 
and let $R_n$ be the set of of hexagons in $\T$ which 
intersect $[0, a\cdot n] \times [0, b\cdot n]$. Let $f_n$ be the event that there is 
a left to right crossing event in $R_n$. (This is the same event as in Example 
\ref{ex.perc} in chapter \ref{ch.BF}, but with $\Z^2$ replaced by $\T$).
By the RSW Theorem \ref{th.RSW}, we know that $\{ f_n \}$ is non-degenerate. 
Conformal invariance tells us that $\Eb{f_n}= \Pb{f_n =1}$
converges as $n\to \infty$. The limit is given by the so-called {\bf Cardy's formula}.
\vskip 0.3 cm

In order to prove that this sequence of Boolean functions $\{ f_n \}$ is noise 
sensitive, we wish to study its influence vector $\InfV(f_n)$ and we would like 
to prove that $\II(f_n) = \| \InfV(f_n) \|_2^2 = \sum \Inf_k(f_n)^2$ decays 
polynomially fast towards 0. (Recall that in these notes, we gave a 
complete proof of Theorem \ref{th:NSmainresult} only in the case where $\II(f_n)$ 
decreases as an inverse polynomial of the number of variables.)

\subsection{Study of the set of influences}

Let $x$ be a site (i.e. a hexagon) in the rectangle $R_n$. One needs to understand
\[
\Inf_x(f_n):= \Pb{ x \text{ is pivotal for }f_n}
\]

\ni
\begin{minipage}{0.5 \textwidth}
It is easy but crucial to note that if $x$ is at distance $d$ from the boundary of $R_n$, in order 
for $x$ to be pivotal,
the {\it four-arm} event described in Chapter \ref{ch.perc} (see Figure \ref{f.armsevents}) has to be satisfied in the ball 
$B(x,d)$ of radius $d$ around the hexagon $x$. See the figure on the right.
\end{minipage}
\hspace{0.02 \textwidth}
\begin{minipage}{0.45 \textwidth}
\includegraphics[width=\textwidth]{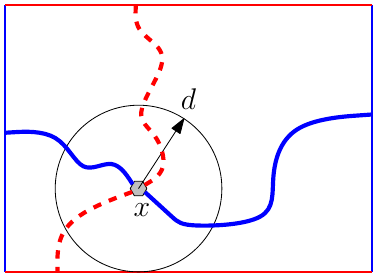}
\end{minipage}
\vskip 0.3 cm

In particular, this implies (still under the assumption that $\mathrm{dist}(x, \p R_n) =d$) that
\[
\Inf_x(f_n) \le \alpha_4(d) = d^{-\frac 5 4 +o(1)}\,,
\]
where $\alpha_4(d)$ denotes the probability of the four-arm event up to distance 
$d$. See Chapter \ref{ch.perc}. The statement
\[
\alpha_4(R) = R^{-5/4 +o(1)}
\]
implies that for any $\eps>0$, there exists a constant $C=C_\eps$, such that for all $R\geq 1$,
\[
\alpha_4(R) \le C \, R^{-5/4 + \eps}\,.
\]

The above bound gives us a very good control on the influences of the points in 
the {\it bulk} of the domain (i.e. the points far from the boundary).
Indeed, for any fixed $\delta>0$, let $\Delta_n^\delta$ be the set of hexagons 
in $R_n$ which are at distance at least $\delta n$ from $\p R_n$.
Most of the points in $R_n$ (except a proportion $O(\delta)$ of these) lie in 
$\Delta_n^\delta$, and for any such point $x\in \Delta_n^\delta$, 
one has by the above argument 
\begin{equation}
\Inf_x(f_n) \le \alpha_4(\delta n) \le  C\, (\delta n)^{-5/4 + \eps} \le  C \delta^{-5/4} n^{-5/4+\eps}\,.
\end{equation}

Therefore, the contribution of these points to $\II(f_n) = \sum_k \Inf_k(f_n)^2$ is bounded 
by $O(n^2) (C \delta^{-5/4} n^{-5/4+\eps})^2 = O(\delta^{-5/2} n^{-1/2 + 2\eps})$. 
As $n\to \infty$, this goes to zero polynomially fast. Since this estimate 
concerns ``almost'' all points in $R_n$, it seems we are close to proving the 
BKS criterion.

\subsection{Influence of the boundary}

Still, in order to complete the above analysis, one has to estimate what the influence of the points near the boundary is.
The main difficulty here is that if $x$ is close to the boundary, the probability for $x$ to be pivotal is not 
related anymore to the above {\it four-arm} event.  Think of the above 
figure when $d$ gets very small compared to $n$.
One has to distinguish two cases:
\bi 
\item $x$ is close to a {\em corner}. This will correspond to a {\it two-arm} event in a quarter-plane.
\item $x$ is close to an {\em edge}. This involves the {\it three-arm} event in the half-plane $\H$.
\ei

Before detailing how to estimate the influence of points near the boundary, let us start by 
giving the necessary background on the involved critical exponents.
\vskip 0.3 cm

\underline{The two-arm and three-arm events in $\H$}.
For these particular events, it turns out that the critical exponents are known 
to be {\it universal}: they are two of the very few critical exponents which 
are known also on the square lattice $\Z^2$. The derivations of these types of 
exponents do not rely on SLE technology but are ``elementary''. Therefore, in this 
discussion, we will consider both lattices $\T$ and $\Z^2$.
\vskip 0.3 cm

\ni
\begin{minipage}{0.4 \textwidth}
The {\it three-arm} event in $\H$ corresponds to the event that 
there are three arms (two open arms and one `closed' arm in the dual) going 
from 0 to distance $R$ and such that they remain in the upper half-plane.
See the figure for a self-explanatory definition. 
The {\it two-arm} event corresponds to just having one open and one closed
arm.
\end{minipage}
\hspace{0.03 \textwidth}
\begin{minipage}{0.54 \textwidth}
\includegraphics[width=\textwidth]{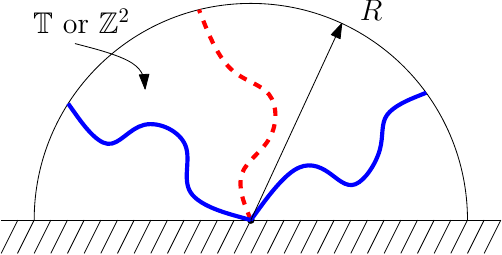}
\end{minipage}

\vskip 0.2 cm

Let $\alpha_2^+(R)$ and $\alpha_3^+(R)$ 
denote the probabilities of these events. As in chapter \ref{ch.perc},
let $\alpha_2^+(r,R)$ and $\alpha_3^+(r,R)$ be the natural extensions to the 
annulus case (i.e. the probability that these events are satisfied in the annulus 
between radii $r$ and $R$ in the upper half-plane).

We will rely on the following result, which goes back as far as we know to M. Aizenman.
See \cite{\WWperc} for a proof of this result.

\begin{proposition}\label{prop.aizenman}
Both on the triangular lattice $\T$ and on $\Z^2$, one has that 
\[
\alpha_2^+(r,R) \asymp (r/R)
\]
and
\[
\alpha_3^+(r,R) \asymp (r/R)^2\,.
\]
Note that, in these special cases, there are no $o(1)$ correction terms
in the exponent. The probabilities are in this case known up to constants.
\end{proposition}

\underline{The two-arm event in the quarter-plane}.
In this case, the corresponding exponent is unfortunately not known on $\Z^2$, so we will need to do some 
work here in the next section, where we will prove noise sensitivity of percolation crossings on $\Z^2$.

\vskip 0.3 cm

\ni
\begin{minipage}{0.55 \textwidth}
The {\it two-arm} event in a corner 
corresponds to the event illustrated 
on the following picture.

We will use the following proposition: 
\begin{proposition}[\cite{\SmirnovWerner}]
If $\alpha_2^{++}(R)$ denotes the probability of this event, then
\[
\alpha_2^{++}(R) = R^{-2+o(1)}\,,
\]
and with the obvious notations
\[
\alpha_2^{++}(r,R)=(r/R)^{2+o(1)}\,.
\]
\end{proposition}

\end{minipage}
\hspace{0.02 \textwidth}
\begin{minipage}{0.4 \textwidth}
\includegraphics[width=\textwidth]{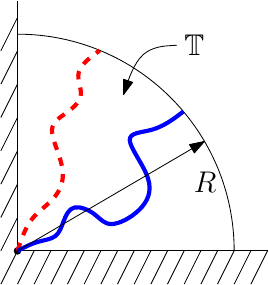}
\end{minipage}
\vskip 0.3 cm

Now, back to our study of influences, we are in good shape (at least for the triangular lattice) since the two critical exponents
arising from the boundary effects are larger than the {\it bulk} exponent $5/4$. This means that it is less likely for a point near the boundary 
to be pivotal than for a point in the bulk. Therefore in some sense the boundary helps us here.

More formally, summarizing the above facts, for any $\eps>0$, there is a constant $C=C(\eps)$ such 
that for any $1\leq r \leq R$,

\begin{equation}\label{e.BulkIsMaximal}
\max \{\alpha_4(r,R), \alpha_3^+(r,R), \alpha_2^{++}(r,R) \} \le C (r/R)^{\frac 5 4 -\eps}\,.
\end{equation}

Now, if $x$ is some hexagon in $R_n$, let $n_0$ be the distance to the 
closest edge of $\p R_n$ and let $x_0$ be the point on $\p R_n$ such 
that $\mathrm{dist}(x,x_0)=n_0$. Next, 
let $n_1\geq n_0$ be the distance from $x_0$
to the closest corner and let $x_1$ be this closest corner.
It is easy to see that for $x$ to be pivotal for $f_n$, the following events all have to be satisfied:
\bi 
\item The four-arm event in the ball of radius $n_0$ around $x$.
\item The $\H$-three-arm event in the annulus centered at $x_0$ of radii $2n_0$ and $n_1$. 
\item The corner-two-arm event in the annulus centered at $x_1$ of radii $2n_1$ and $n$.
\ei

By independence on disjoint sets, one thus concludes that 
\begin{eqnarray}
\Inf_x(f_n) &\le & \alpha_4(n_0) \, \alpha_3^+(2n_0, n_1) \, \alpha_2^{++}(2n_1, n) \nonumber \\
& \le & O(1) n^{-5/4 + \eps}\,. \nonumber
\end{eqnarray}

\subsection{Noise sensitivity of crossing events}

This uniform bound on the influences over the whole domain $R_n$ enables us to conclude that the BKS criterion 
is indeed verified. Indeed,
\begin{eqnarray}
\II(f_n) = \sum_{x\in R_n} \Inf_x(f_n)^2  \le  C n^2  (n^{-5/4+\eps})^2 = C n^{-1/2 + 2\eps}\,,
\end{eqnarray}
where $C=C(a,b,\eps)$ is a universal constant. By taking $\eps<1/4$,
this gives us the desired polynomial decay on $\II(f_n)$, which by 
Proposition \ref{pr.BKSquantitative}) implies 

\begin{theorem}[\cite{\BKS}]
The sequence of percolation crossing events $\{ f_n \}$ on $\T$ is noise sensitive.
\end{theorem}

We will give some other consequences (for example, to sharp thresholds)
of the above analysis on the influences of the crossing events in a later section.

\section{The case of $\Z^2$ percolation}
Let $R_n$ denote similarly the $\Z^2$ rectangle closest to 
$[0, a\cdot n] \times [0, b\cdot n]$ and let $f_n$ be the corresponding 
left-right crossing event (so here this corresponds exactly to example \ref{ex.perc}).
Here one has to face two main difficulties:
\bi
\item The main one is that due to the missing ingredient of {\it conformal invariance}, one does not have at our disposal the value 
of the {\it four-arm} critical exponent (which is of course believed to be $5/4$). In fact, even the {\it existence} of a critical exponent is an open problem.
\item  The second difficulty (also due to the lack of conformal invariance) is that it is now 
slightly harder to deal with boundary issues. Indeed, one can still use the above bounds on $\alpha_3^+$ which 
are {\it universal}, but the exponent $2$ for $\alpha_2^{++}$ is not known for $\Z^2$. So this requires some more analysis.
\ei
Let us start by taking care of the boundary effects.

\subsection{Handling the boundary effect}

What we need to do in order to carry through the above analysis 
for $\Z^2$ is to obtain a reasonable estimate on $\alpha_2^{++}$. 
Fortunately, the following bound, which follows immediately from 
Proposition \ref{prop.aizenman}, is sufficient.

\begin{equation}\label{e.cornerbound}
\alpha_2^{++}(r,R) \le O(1) \frac r R\,.
\end{equation}

Now let $e$ be an edge in $R_n$. We wish to bound from above $\Inf_e(f_n)$.
We will use the same notation as in the case of the triangular lattice: recall the definitions
of $n_0,x_0,n_1,x_1$ there.

We obtain in the same way
\begin{equation}\label{e.boundinf}
\Inf_e(f_n) \le   \alpha_4(n_0) \, \alpha_3^+(2n_0, n_1) \, \alpha_2^{++}(2n_1, n)\,.
\end{equation}

At this point, we need another {\it universal} exponent, which goes back also to M. Aizenman:
\begin{theorem}[M. Aizenman, see \cite{\WWperc}] \label{t.5armexponent}
Let $\alpha_5(r,R)$ denote the probability that there are 5 arms (with four of them being of `alternate colors').
Then there are some universal constants $c,C>0$ such that both for $\T$ and $\Z^2$, one has for all $1\le r \le R$,
\[
c\bigl( \frac r R\bigr)^2 \le \alpha_5(r,R) \le C \bigl( \frac r R\bigr)^2\,.
\]
\end{theorem}

This result allows us to get a lower bound on $\alpha_4(r,R)$. Indeed, it is clear that
\begin{equation}\label{e.comparison43}
\alpha_4(r,R) \ge \alpha_5(r,R) \ge \Omega(1) \alpha_3^+(r,R)\,.
\end{equation}

In fact, one can obtain the following better lower bound on $\alpha_4(r,R)$ which
we will need later.

\begin{lemma}\label{l.alpha4}
There exists some $\eps>0$ and some constant $c>0$ such that
for any $1\le r \le R$,
\[
\alpha_4(r,R) \ge c (r/R)^{2-\eps}\,.
\]
\end{lemma}

\proof
There are several ways to see why this holds, none of them being either very hard
or very easy. One of them is to use {\bf Reimer's inequality} (see 
\cite{\Grimm}) which in this case would imply that 
\begin{equation}\label{e.reimerapplication}
\alpha_5(r,R) \le \alpha_1(r,R) \alpha_4(r,R)\,.
\end{equation}
The RSW Theorem \ref{th.RSW} can be used to show that
\[
\alpha_1(r,R) \le (r/R)^\alpha
\]
for some positive $\alpha$. By Theorem \ref{t.5armexponent}, we are done.
[See [\cite{\GPS}, Section 2.2 as well as the appendix] for more on these 
bounds.] 
\qed

\ni
Combining  ~\eqref{e.boundinf} with~\eqref{e.comparison43}, one obtains
\begin{eqnarray}
\Inf_e(f_n) &\le & O(1)\alpha_4(n_0) \alpha_4(2n_0, n_1) \alpha_2^{++}(2n_1, n) \nonumber \\
& \le & O(1) \alpha_4(n_1) \,  \frac {n_1} n \,, \nonumber
\end{eqnarray}
where in the last inequality we used quasi-multiplicativity 
(Proposition \ref{pr.quasi}) as well as the bound given by ~\eqref{e.cornerbound}.
\vskip 0.3 cm

\ni
\begin{minipage}{0.4 \textwidth}
Recall that we want an upper bound on $\II(f_n)= \sum \Inf_e(f_n)^2$. 
In this sum over edges $e\in R_n$, let us divide the set of edges into dyadic annuli 
centered around the 4 corners as in the next picture.
\end{minipage}
\hspace{0.1\textwidth}
\begin{minipage}{0.4 \textwidth}
\includegraphics[width=\textwidth]{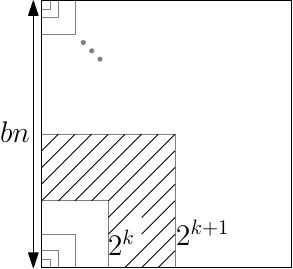}
\end{minipage}
\vskip 0.2 cm

Notice that there are $O(1) 2^{2k}$ edges in an annulus of radius $2^k$. This enables us to bound $\II(f_n)$ as follows:
\begin{eqnarray}\label{e.IIandalpha4}
\sum_{e \in R_n} \Inf_e(f_n)^2  &\le &  O(1) \sum_{k=1}^{\log_2 n +O(1)} 2^{2k} \Bigl( \alpha_4(2^k) \frac {2^k}{n} \Bigr)^2 \nonumber\\ 
& \le & O(1)\;  \frac 1 {n^2} \sum_{k\le \log_2 n+O(1)} 2^{4k}\, \alpha_4(2^k)^2\,. 
\end{eqnarray}

It now remains to obtain a good upper bound on $\alpha_4(R)$, for all $R\geq 1$.

\subsection{An upper bound on the four-arm event in $\Z^2$}

This turns out to be a rather non-trivial problem. Recall that we obtained an easy lower 
bound on $\alpha_4$ using $\alpha_5$ (and Lemma \ref{l.alpha4} strengthens this 
lower bound). For an upper bound, completely different ideas are required. 
On $\Z^2$, the following estimate is available for the 4-arm event.

\begin{proposition}\label{pr.upperboundalpha4}
For critical percolation on $\Z^2$, there exists constants $\eps, C>0$ such 
that for any $R\ge 1$, one has 
\[
\alpha_4(1,R) \le C \Bigl( \frac 1 R \Bigr)^{1+\eps}\,.
\]
\end{proposition}

Before discussing where such an estimate comes from,
let us see that it indeed implies a polynomial decay for $\II(f_n)$.

Recall equation~\eqref{e.IIandalpha4}. Plugging in the above estimate, this 
gives us
\begin{eqnarray}
\sum_{e \in R_n} \Inf_e(f_n)^2 & \le & O(1) \;  \frac 1 {n^2} 
\sum_{k\le \log_2 n+O(1)} 2^{4k}\, (2^k)^{-2-2\eps} \nonumber \\
& \le & O(1) \frac 1 {n^2}  n^{2-2\eps} = O(1) n^{-2\eps}\,, \nonumber
\end{eqnarray}
which implies the desired polynomial decay and thus the fact that $\{ f_n \}$ is 
noise sensitive by Proposition \ref{pr.BKSquantitative}).
\vskip 0.3 cm

Let us now discuss different approaches which enable one to prove 
Proposition \ref{pr.upperboundalpha4}.

\bi
\item[(a)] Kesten proved implicitly this estimate in his celebrated paper 
\cite{\KestenScaling}. His main motivation for such an estimate
was to obtain bounds on the corresponding critical exponent which governs the so-called {\it critical length}.

\item[(b)] In \cite{\BKS}, in order to prove noise sensitivity of percolation 
using their criterion on $\II(f_n)$, the authors  referred to 
\cite{\KestenScaling}, but they also gave a completely different approach which 
also yields this estimate.

\ni Their alternative  approach is very nice: finding an upper bound for $\alpha_4(R)$ is related to finding an upper bound for the influences for crossings of
an $R\times R$ box.
For this, they noticed the following nice phenomenon: if a monotone function $f$ happens to be very little correlated with majority,
then its influences have to be small. The proof of this phenomenon uses for the first time in this context the concept of ``randomized algorithms''.
For more on this approach, see Chapter \ref{ch.RA}, 
which is devoted to these types of ideas.

\item[(c)] In \cite{\SS}, the concept of randomized algorithms is used in a more 
powerful way. See again Chapter \ref{ch.RA}. In this chapter, we provide
a proof of this estimate in Proposition \ref{prop.4armexponent}.

\ei

\begin{remark}
It turns out that that a multi-scale version of 
Proposition \ref{pr.upperboundalpha4} stating that
$\alpha_4(r,R) \le C \Bigl( \frac r R \Bigr)^{1+\eps}$
is also true. However, none of the three arguments given above
seem to prove this stronger version. A proof of this stronger version is
given in the appendix of \cite{\SchrammSmirnovNoise}. Since this multi-scale
version is not needed until Chapter \ref{ch.SNS}, 
we stated here only the weaker version.
\end{remark}

\section{Some other consequences of our study of influences}

In the previous sections, we handled the boundary effects in order to check that $\II(f_n)$
indeed decays polynomially fast. Let us list some related results implied by this analysis.

\subsection{Energy spectrum of $f_n$}
We start by a straightforward observation: since the $f_n$ are monotone, we have by
Proposition \ref{pr.inf&spectrum.monotone} that 
\[
\widehat{f_n}(\{x\}) = \frac 1 2 \Inf_x(f_n)\,,
\]
for any site $x$ (or edge $e$) in $R_n$.
Therefore, the bounds we obtained on $\II(f_n)$ imply the following control 
on the first layer of the energy spectrum of the crossing events $\{f_n\}$.

\begin{corollary}\label{cor.FirstLayerDecay}
Let $\{ f_n\}$ be the crossing events of the rectangles $R_n$.
\bi
\item If we are on the triangular lattice $\T$, then we have the bound
\[
E_{f_n}(1) = \sum_{|S|=1} \widehat{f_n}(S)^2 \le n^{-1/2+o(1)}\,.
\]
\item On the square lattice $\Z^2$, we end up with the weaker estimate
\[
E_{f_n}(1) \le C\, n^{-\eps}\,,
\]
for some $\eps, C >0$. 
\ei
\end{corollary}

\subsection{Sharp threshold of percolation}
The above analysis gave an upper bound on $\sum_k \Inf_k(f_n)^2$. 
As we have seen in the first chapters, the total influence $\Inf(f_n) = \sum_k \Inf_k(f_n)$
is also a very interesting quantity. Recall that, by Russo's formula, this is the quantity which shows ``how sharp'' the 
threshold is for $p \mapsto \P_p [ f_n =1 ]$. 

The above analysis allows us to prove the following.

\begin{proposition}\label{pr.UpperboundInf}
Both on $\T$ and $\Z^2$, one has 
\[
\Inf(f_n) \asymp \, n^2 \alpha_4(n)\,.
\]
In particular, this shows that on $\T$ that 
\[
\Inf(f_n) \asymp n^{3/4+o(1)}\,.
\]
\end{proposition}

\begin{remark}
Since $f_n$ is defined on $\{-1,1\}^{O(n^2)}$, note that the Majority function 
defined on the same hypercube has a much sharper threshold than the 
percolation crossings $f_n$.
\end{remark}

\proof
We first derive an upper bound on the total influence.
In the same vein (i.e., using dyadic annuli and quasi-multiplicativity) as we derived
(\ref{e.IIandalpha4}) and with the same notation one has
\begin{eqnarray}
\Inf(f_n) = \sum_e \Inf_e(f_n) & \le & \sum_e O(1) \alpha_4(n_1) \frac {n_1} {n} \nonumber  \\
& \le & O(1) \frac 1 n \sum_{k \le \log_2 n+O(1)} 2^{3k} \alpha_4(2^k)\,.  \nonumber
\end{eqnarray}
Now, and this is the main step here, using quasi-multiplicativity one has $\alpha_4(2^k) \le O(1) \frac {\alpha_4(n)}{\alpha_4(2^k, n)}$, which gives us

\begin{eqnarray}
\Inf(f_n) & \le & O(1) \frac {\alpha_4(n)} n \sum_{k \le \log_2 n+O(1)} 2^{3k} \frac 1 {\alpha_4(2^k, n)}  \nonumber \\
& \le & O(1) \frac {\alpha_4(n)} n \sum_{k \le \log_2 n +O(1)} 2^{3k} \frac{n^2}{2^{2k}}\; \text{ since }\alpha_4(r,R) \ge \alpha_5(r,R)\asymp (r/R)^{-2}\nonumber \\
& \le & O(1) n\, \alpha_4(n) \sum_{k \le \log_2 n+O(1)} 2^k \nonumber \\
&\le & O(1) n^2 \alpha_4(n)\; \nonumber
\end{eqnarray} 
as desired. 

For the lower bound on the total influence, we proceed as follows.
One obtains a lower bound by just summing over the influences of
points whose distance to the boundary is at least $n/4$. It would
suffice if we knew that for such edges or hexagons, the influence is
at least a constant times  $\alpha_4(n)$.
This is in fact known to be true. 
It is not very involved and is part of the {\it folklore} results in
percolation. However, it still would lead us too far from our topic.
The needed technique is known under the name of {\bf separation of
arms} and is clearly related to the statement of quasi-multiplicativity.
See \cite{\WWperc} for more details.
\QED

\section{Quantitative noise sensitivity}

In this chapter, we have proved that the sequence of 
crossing events $\{ f_n \}$ is noise sensitive. This can be roughly translated 
as follows: for any fixed level of noise $\eps>0$, as
$n\to\infty$, the large scale clusters of $\omega$ in the window $[0,n]^2$ are 
asymptotically independent of the large clusters of $\omega_\eps$.

\begin{remark}
Note that this picture is correct, but in order to make it rigorous, this would 
require some work, since so far 
we only worked with left-right crossing events. The non-trivial step here is to prove that in some sense, in the scaling limit $n\to \infty$,
any macroscopic property concerning percolation (e.g., diameter of clusters)
is measurable with respect to the $\sigma$-algebra generated by the crossing 
events.
This is a rather subtle problem since we need to make precise what kind of information we keep in what we call  the ``scaling limit'' of percolation
(or subsequential scaling limits in the case of $\Z^2$). 
An example of something which is not present in the scaling limit is whether
one has more open sites than closed ones since by noise sensitivity we know 
that this is asymptotically uncorrelated with crossing events.
We will not need to  discuss these notions of scaling limits more in these lecture notes, since the focus is mainly on the discrete model itself including 
the model of dynamical percolation which is presented at the end of these lecture
notes.
\end{remark}

At this stage, a natural question to ask is to what extent the percolation picture 
is sensitive to noise.
In other words, can we let the noise $\eps=\eps_n$ go to zero with the ``size of the system'' $n$, and yet keep this independence of large scale structures 
between $\omega$ and $\omega_{\epsilon_n}$ ? If yes, can we give 
{\em quantitative estimates} on how fast the noise $\eps=\eps_n$ may go to zero?
One can state this question more precisely as follows.

\begin{question}\label{q.QNS}
If $\{f_n\}$ denote our left-right crossing events,
for which sequences of noise-levels $\{ \eps_n \}$ do we have 
\[
\lim_{n\to \infty} \Cov [f_n(\omega), f_n(\omega_{\eps_n})] = 0 \; ?
\]
\end{question}

The purpose of this section is to briefly discuss this question based on the results we have obtained so far.

\subsection{Link with the energy spectrum of $\{ f_n \}$}

It is an exercise to show that Question \ref{q.QNS}
is essentially equivalent to the following one.

\begin{question}
For which sequences $\{k_n\}$ going to infinity do we have 
\[
\sum_{m=1}^{k_n} E_{f_n}(m) = \sum_{1\le |S| \le k_n} \widehat{f_n}(S)^2 \underset{n\to \infty}{\longrightarrow} 0\; ?
\]
\end{question}

Recall that we have already obtained some relevant information on this question. 
Indeed, we have proved in this chapter that 
$\II(f_n)=\sum_x \Inf_x(f_n)^2$ decays polynomially fast towards 0 (both on $\Z^2$ and $\T$). Therefore Proposition \ref{pr.BKSquantitative}
tells us that for some constant $c>0$, one has for both $\T$ and $\Z^2$ that

\begin{equation}\label{e.LowerTailLog}
\sum_{1\le |S| \le c \log n} \widehat{f_n}(S)^2 \to 0\,.
\end{equation}
Therefore, back to our original question \ref{q.QNS}, this gives us the following 
quantitative statement: if the noise $\eps_n$ satisfies
$\eps_n \gg \frac{1}{\log n}$, then 
$f_n(\omega)$ and $f_n(\omega_{\eps_n})$ are asymptotically independent.

\subsection{Noise stability regime}
Of course, one cannot be too demanding on the rate of decay of $\{\eps_n\}$. 
For example if $\eps_n \ll \frac 1 {n^2}$, then in the window $[0,n]^2$, with 
high probability, the configurations $\omega$ and $\omega_{\eps_n}$ are identical. 
This brings us to the next natural question concerning the 
{\em noise stability regime} of crossing events.

\begin{question}
Let $\{f_n\}$ be our sequence of crossing events. For which sequences $\{ \eps_n \}$
do we have 
\[
\Pb{ f_n(\omega) \neq f_n(\omega_{\eps_n})} \underset{n\to\infty}{\longrightarrow} 0\; ?
\]
\end{question}

It is an exercise to show that this question
is essentially equivalent to the following one.

For which sequences $\{k_n\}$ do we have 
\[
\sum_{|S| > k_n} \widehat{f_n}(S)^2 \to 0\; ?
\]

Using the estimates of the present chapter, one can give the following 
non-trivial bound on the noise stability regime of $\{f_n\}$.

\begin{proposition} \label{prop.uppertail}
Both on $\Z^2$ and $\T$, if 
\[
\eps_n = o\Bigl(\frac 1 {n^2 \alpha_4(n)}\Bigr)\,,
\]
then 
\[
\Pb{ f_n(\omega) \neq f_n(\omega_{\eps_n})} \underset{n\to\infty}{\longrightarrow} 0
\]
On the triangular grid, using the critical exponent, this gives us a bound 
of $n^{-3/4}$ on the noise stability regime of percolation.
\end{proposition}

\proof  
Let $\{\eps_n\}$ be a sequence satisfying the above assumption.
There are $O(n^2)$ bits concerned. For simplicity, assume
that there are exactly $n^2$ bits.
Let us order these in some arbitrary way: $\{x_1,\ldots, x_{n^2}\}$ (or on $\Z^2$, $\{e_1,\ldots, e_{n^2}\}$).

Let $\omega=\omega_0=(x_1,\ldots,x_{n^2})$ be sampled according to the uniform 
measure. Recall that
the noised configuration $\omega_{\eps_n}$ is produced as follows: for each $i\in[n^2]$, resample the bit $x_i$ with probability
$\eps_n$, independently of everything else, obtaining the bit $y_i$. (In particular $y_i \neq x_i$ with probability $\eps_n/2$).

Now for each $i\in[n^2]$ define the intermediate configuration 
\[
\omega_i:= (y_1,\ldots,y_i,x_{i+1},\ldots,x_{n^2})
\]

Notice that for each $i\in[n^2]$, $\omega_i$ is also sampled according to the 
uniform measure and one has for each $i\in \{1,\ldots, n^2\}$ that 
\[
\Pb{f_n(\omega_{i-1}) \neq f_n(\omega_{i})} = (\eps_n/2) \, \Inf_{x_{i}}(f_n)\,.
\]

Summing over all $i$, one obtains
\begin{eqnarray}
\Pb{f_n(\omega) \neq f_n(\omega_{\eps_n})} &=& \Pb{f_n(\omega_0)\neq f_n(\omega_{n^2})} \nonumber \\
&\le & \sum_{i=0}^{n^2-1} \Pb{f_n(\omega_i)\neq f_n(\omega_{i+1})} \nonumber \\
&= & (\eps_n/2) \,\sum_{i=1}^{n^2} \Inf_{x_i}(f_n) \nonumber \\
& = & (\eps_n/2) \, \Inf(f_n) \nonumber \\
&\le &\eps_n O(1) n^2 \alpha_4(n)\; \text{ by Proposition \ref{pr.UpperboundInf},} \nonumber
\end{eqnarray}
which concludes the proof.
\QED

\subsection{Where does the spectral mass lies?}
Proposition \ref{prop.uppertail} (together with 
Exercise \ref{exercise.stabilityandspectrum} in Chapter \ref{ch.SM})
implies that the Fourier coefficients of $\{f_n\}$ satisfy 
\begin{equation}\label{e.uppertailpolynom}
\sum_{|S| \gg n^2 \alpha_4(n)} \widehat{f_n}(S)^2 \underset{n\to\infty}{\longrightarrow} 0\,.
\end{equation}

\begin{figure}[!h]
\begin{center}
\includegraphics[width=0.9\textwidth]{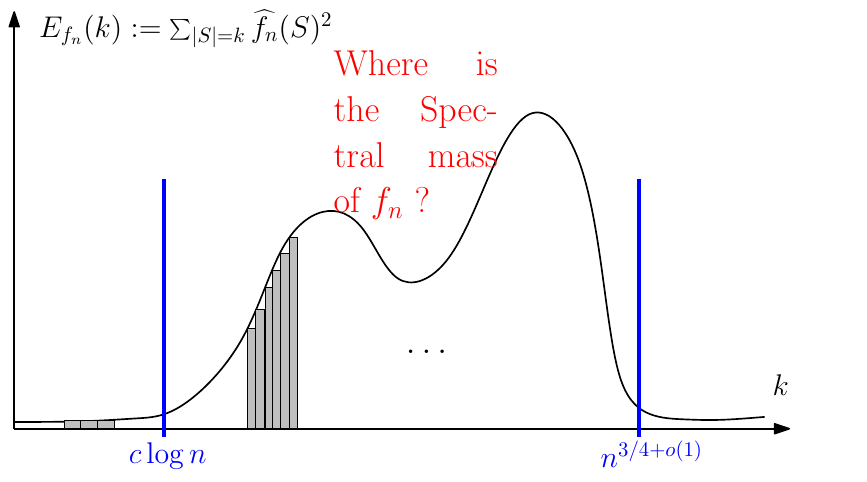}
\end{center}
\caption{
This picture summarizes our present knowledge of the energy spectrum of 
$\{f_n\}$ on the triangular lattice $\T$. Much remains to be understood to know 
where, in the range $[\Omega(\log n), n^{3/4+o(1)}]$, the spectral mass lies. 
This question will be analyzed in the following chapters.
}\label{f.quantitative}
\end{figure}

From Lemma \ref{l.alpha4}, we know that even on $\Z^2$, $n^2 \alpha_4(n)$ is 
larger than $n^\eps$ for some exponent $\eps>0$. Combining the estimates on 
the spectrum that we achieved so far (equations~\eqref{e.LowerTailLog} 
and~\eqref{e.uppertailpolynom}),
we see that in order to localize the spectral mass 
of $\{ f_n\}$, there is still a missing gap. See Figure \ref{f.quantitative}.

\vskip 0.3 cm

For our later applications to the model of dynamical percolation (in the
last chapter of these lecture notes), a better understanding of the noise 
sensitivity of percolation than the ``logarithmic'' control we achieved so far 
will be needed.

\note{Comment}{In \cite{\KestenScaling}, his statement which concerns us is the following: the exponent $\nu$ which dictates the 
correlation length satsifies the inequality $\nu>1$. And of course, this is possible only if the exponent corresponding to 
$\alpha_4$ is $>1$. Can one get a better bound from Kesten than just $\nu>1$ ? In section 5, he gives the lower bound 
\[
\nu \geq \frac {\delta + 1 }{\delta}\,.
\]

$\delta$ is critical exponent about the purely critical slice: it corresponds to 
\[
\P_{p_c} ( | \mathcal{C}(0) | \ge n) \approx n^{-1/\delta} \,.
\]
With a good estimate on $\alpha_1$ (which is easy using Aizenman), one can get I think a good control on $\delta$.
Note that his results seems to hold also for {\bf site percolation} on $\Z^2$ !

$\delta$ and $\xi_1$ follow the following relation
\[
\delta = \frac{2-\xi_1}{\xi_1}\,,
\]
If things work out well, this might give
\[
\nu \ge \frac {2}{2-\xi_1}
\]
So again, we need a non-trivial lower bound on $\xi_1$ which is equivalent to an upper bound on $\alpha_1(n)$.
\vskip 0.2 cm

Note also the interesting (but easy) phenomenon that $\nu=1/2<1$ for near-critical percolation on the Tree.
\vskip 0.3 cm

He also obtains bound on $\gamma$ (which has to do with the mean cluster size) at the very end, he proves that $\gamma\geq 8/5$ on $\Z^2$ and he also get the weaker $\gamma\geq 4/3$ for site
percolation on $\Z^2$ ! 
}

\chapter*{Exercise sheet on chapter \ref{ch.FE}}
\setcounter{exercise}{0}
\setcounter{section}{0}

Instead of being the usual exercise sheet, this page will be devoted to a single Problem whose goal will be to 
do ``hands-on'' computations of the first layers of the energy spectrum of the percolation crossing events $f_n$.
Recall from Proposition \ref{pr.NSeq} that a sequence of Boolean functions $\{ f_n \}$ is noise sensitive if and only if for any fixed $k\geq 1$,

\begin{equation}
\sum_{m=1}^k \, \sum_{|S|=m} \hat f_n(S)^2  = \sum_{m=1}^k E_{f_n}(m) \underset{n\to \infty}{\longrightarrow} 0\,. \nonumber
\end{equation}

In the present chapter, we obtained 
(using Proposition \ref{pr.inf&spectrum.monotone})
that this is indeed the case for $k=1$. The purpose here is to check by simple combinatorial arguments (without relying on hypercontractivity)
that it is still the case for $k=2$ and to convince ourselves that it works for all layers $k\geq 3$.

To start with, we will simplify our task by working on the torus $\Z^2 / n\Z^2$. This has the very nice advantage that there are no boundary issues here.
\vskip 0.5 cm

\section*{Energy spectrum of crossing
events on the torus (study of the first layers)}

Let $T_n$ be either the square grid torus $\Z^2/n\Z^2$ or the triangular grid torus $\T/ n\T$.
Let $f_n$ be the indicator of the event that there is an open circuit along the first coordinate of $T_n$.

\begin{enumerate}

\item Using RSW, prove that there is a constant $c>0$ such that for all $n\geq 1$, 
\[
c\le \Pb{f_n=1} \le 1-c\,.
\]
(In other words, $\{ f_n \}$ is non-degenerate.)

\item Show that for all edges $e$ (or sites $x$) in $T_n$
\[
\Inf_e(f_n) \le \alpha_4(\frac n 2)\,.
\]

\item Check that the BKS criterion (about $\II(f_n)$) is satisfied. Therefore $\{ f_n \}$ is noise-sensitive

\vskip 0.3 cm

From now on, one would like to forget about the BKS Theorem and try to do some 
hands-on computations in order to get a feeling why most frequencies should be 
large.

\item Show that if $x,y$ are two sites of $T_n$ (or similarly if $e,e'$ are two edges of $T_n$), then 
\[
 |\hat f(\{x,y\})| \le 2\Pb{\text{ $x$ and $y$ are pivotal points}}\,.
\]
Does this result hold for general Boolean functions?

\item Show that if $d:= |x-y|$, then 
\[
\Pb{\text{ $x$ and $y$ are pivotal points }} \le O(1) \frac {\alpha_4(n/2)^2} { \alpha_4(\frac d 2, \frac n 2)}\,.
\]
(Hint: use Proposition \ref{pr.quasi}.)

\item On the square lattice $\Z^2$, by carefully summing over all edges $e,e' \in T_n\times T_n$, show that 
\[
E_{f_n}(2) = \sum_{|S|=2} \widehat{f_n}(S)^2 \le \,  O(1) n^{-\eps}\,,
\]
for some exponent $\eps>0$.

\ni
Hint: you might decompose the sum in a dyadic way (as we did many 
times in the present section) depending on the mutual distance $d(e,e')$.

\item On the triangular grid, what exponent does it give for the decay of $\E_{f_n}(2)$? 
Compare with the decay we found in Corollary \ref{cor.FirstLayerDecay}
about the decay of the first layer $E_{f_n}(1)$ (i.e. $k=1$). See also Lemma
\ref{lem:talagrandext} in this regard. Discuss this.

\item For $\T$, what do you expect for higher (fixed) values of $k$? 
(I.e.\ for $E_{f_n}(k)$, $k\geq 3$) ?

\item {\it (Quite hard)} Try to obtain a nonrigorous
combinatorial argument similar to the one above in the particular case $k=2$, 
that for any fixed layer $k\geq 1$, 
\[
E_{f_n}(k) \underset{n\to \infty}{\longrightarrow} 0\,.
\]
This would give us an alternative proof of noise sensitivity of percolation
(at least in the case of the torus $T_n$) not relying on Theorem \ref{th:NSmainresult}.
\end{enumerate}

Observe that one can do similar things for rectangles but then one has to deal with
boundary issues.

\note{Comment}
{Other ideas of exercise :
\bi 
\item[-] Extend the BKS Theorem to other product measures $\mu_p$ {\bf (If one wants to do it, this should be a substantial 
Problem session since one needs to extend Bonami-Beckner and detail the Fourier basis outside of $p=1/2$. Unless there is 
an easier proof similar as the BKKKL extension ?)}. What does it say for site percolation ?
\ei

So it seems after discussion that it might not be very ``natural'' to use this idea of extending the problem to the continuous hypercube.
Yet, it is very likely that the spectral proof (with the appropriate basis) will carry easily with additional constants which will depend on $p$.
\vskip 0.3 cm

{\bf Problem: it seems RSW in a strong sense is more or less known only for $Z^2$, $\T$ and trivial extensions of these \ldots
Check with bond percolation on $\T$.} 
If that's the case, the motivation to extend BKS is smaller.
\vskip 0.3 cm

By the way, is RSW known for site percolation?  C: I believe not.
}

\chapter[\; Anomalous fluctuations]{Anomalous fluctuations}\label{ch.AF}

\note{Timing}{///////// 20 OR  40 minutes  OR  1 hour ///////// \vskip 1 cm}

In this lecture, our goal is to extend the technology we used to prove the KKL 
Theorems on influences and 
the BKS Theorem on noise sensitivity to a slightly different context: the study of fluctuations in {\bf first passage percolation}.

\section{The model of first passage percolation}
Let us first explain what the model is. Let $0<a<b$ be two positive numbers. 
We define a {\bf random metric} on the graph $\Z^d$, $d\geq 2$ as follows. 
Independently for each edge $e\in \E^d$, 
fix its length $\tau_e$ to be $a$ with probability 1/2 and $b$ with
probability $1/2$. This is represented by a uniform
configuration $\omega \in \{-1, 1\}^{\E^d}$.

This procedure induces a well-defined (random) metric $\mathrm{dist}_\omega$ on $\Z^d$ in the usual fashion.
For any vertices
$x,y \in \Z^d$, let
\[
\mathrm{dist}_\omega(x,y) := \inf_{ \begin{array}{ll} \text{paths }\gamma=\{e_1,\ldots, e_k\} \\ \text{ connecting }x \to y \end{array}}  \left\lbrace \sum \tau_{e_i}(\omega) \right\rbrace \,.
\]

\begin{remark}
In greater generality, the lengths of the edges are i.i.d.\ non-negative random variables, 
but here, following \cite{\BKSfpp}, we will restrict ourselves to the above uniform 
distribution on $\{a,b\}$ to simplify the exposition; see \cite{\BenaimRossignol} 
for an extension to more general laws.
\end{remark}

One of the main goals in first passage percolation is to understand the 
large-scale properties of this random metric space.
For example, for any $T\geq 1$, one may consider the (random) ball 
\[
B_\omega(x,T):= \{ y\in \Z^d: \mathrm{dist}_\omega(x,y) \le T \}.
\]

To understand the name {\it first passage percolation}, one can think of this model as follows. Imagine that water is pumped in at vertex $x$, and that 
for each edge $e$, it takes $\tau_e(\omega)$ units of time for the water 
to travel across the edge $e$. Then, $B_\omega(x,T)$ represents the 
region of space that has been wetted by time $T$.

\begin{figure}[!htp]
\begin{center}
\includegraphics[width=0.5 \textwidth]{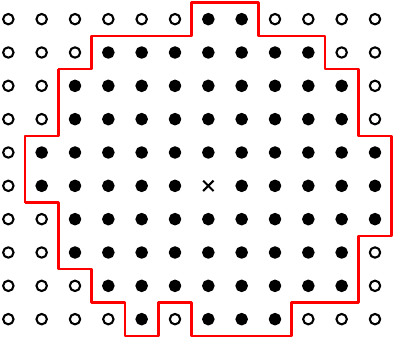}
\end{center}
\caption{A sample of a wetted region at time $T$, i.e.\ $B_\omega(x,T)$, in first passage percolation.}\label{f.fpp}
\end{figure}

An application of subadditivity shows that the renormalized ball 
$\frac 1 T B_\omega(0,T)$ converges as $T\to \infty$ towards a deterministic shape which 
can in certain cases be computed explicitly. This is a kind of ``geometric law of 
large numbers''. Whence the natural question:
\begin{question}
Describe the {\em fluctuations} of $B_\omega(0,T)$ around its asymptotic 
deterministic shape.
\end{question}

This question has received tremendous interest in the last 15 years or so. 
It is widely believed that these fluctuations should be in some sense 
``universal''. More precisely, the behavior of $B_\omega(0,T)$ around its 
limiting shape should not depend on the ``microscopic'' 
particularities of the model such as the law on the edges lengths
but only on the dimension $d$ of the underlying graph.
The shape itself depends on the other hand of course on the microscopic parameters, 
in the same way as the critical point depends on the graph in percolation.
\vskip 0.2 cm

In the two-dimensional case, using very beautiful combinatorial bijections with random matrices, certain cases of {\em directed} last passage percolation (where 
the law on the edges is taken to be geometric or exponential) have been understood 
very deeply. For example, it is known (see \cite{\JohanssonShapeFluct}) that the 
fluctuations of the ball of radius $n$ (i.e. the points whose last passage 
times are below $n$) around $n$ times its asymptotic deterministic shape are 
of order $n^{1/3}$ and the law of these fluctuations properly renormalized 
follows the Tracy-Widom distribution. Very interestingly, the fluctuations of the largest 
eigenvalue of GUE ensembles also follow this distribution.
\vskip 0.3 cm

\section{State of the art}

Returning to our initial model of (non-directed) first passage percolation, 
it is thus conjectured that, for dimension $d=2$,
fluctuations are of order $n^{1/3}$ following a Tracy-Widom Law.
Still, the current state of understanding of this model is far from this conjecture. 
\vskip 0.2 cm

Kesten first proved that the fluctuations of the ball of radius $n$ are at most $\sqrt{n}$ (this did not yet exclude a possible Gaussian behavior
with Gaussian scaling). 
Benjamini, Kalai and Schramm then strengthened this result by showing that the fluctuations are sub-Gaussian. This is still 
far from the conjectured $n^{1/3}$-fluctuations, but their approach has the great advantage of being very 
general; in particular their result holds in any dimension $d\geq 2$.

Let us now state their main theorem
concerning the fluctuations of the metric $\mathrm{dist}$.

\begin{theorem}[\cite{\BKSfpp}]\label{th.fpp}
For all $a,b,d$, there exists an absolute constant $C=C(a,b,d)$ such that in $\Z^d$,
\[
\Var(\mathrm{dist}_\omega(0,v)) \le C\frac {|v|}{\log |v|}
\]
for any $v\in \Z^d, |v|\geq 2$.
\end{theorem}

\note{Comment}{Keep in mind the fact that this is related (via Chatterjee) to the sensitivity of Geodesics for the lecture}

To keep things simple in these notes, we will only prove the analogous 
statement on the torus where one has more symmetries and invariance to 
play with.

\section{The case of the torus}
Let $\T_m^d$ be the $d$-dimensional torus $(\Z/m\Z)^d$.
As in the above lattice model, independently for each edge of
$\T_m^d$, we choose its length to be either $a$ or $b$ equally
likely. We are interested here
in the smallest length among all closed paths $\gamma$ ``winding'' 
around the torus along the first coordinate $\Z/ m\Z$ (i.e. those paths $\gamma$
which when
projected onto the first coordinate have winding number one). In $\cite{\BKSfpp}$, this is called the shortest {\it circumference}.
For any configuration $\omega\in \{a,b\}^{E(\T_m^d)}$, this shortest
circumference is denoted by $\Circ_m(\omega)$.

\begin{figure}[!htp]
\begin{center}
\includegraphics[width=0.8 \textwidth]{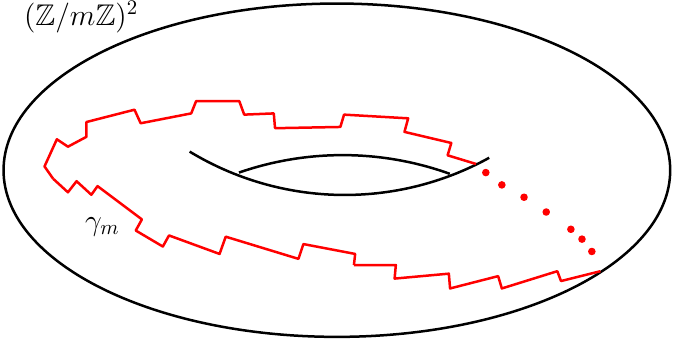}
\end{center}
\caption{The shortest geodesic along the first coordinate for the random metric $\mathrm{dist}_\omega$
on $(\Z/m\Z)^2$.}\label{f.frontpage}
\end{figure}

\begin{theorem}[\cite{\BKSfpp}]\label{th.torus}
There is a constant $C=C(a,b)$ (which does not depend on the dimension $d$), such that
\[
\var(\Circ_m(\omega)) \leq C \frac m {\log m}.
\]
\end{theorem}

\begin{remark}
A similar analysis as the one carried out below works in greater 
generality: if $G=(V,E)$ is some finite connected graph endowed with a random metric $d_\omega$ with $\omega\in \{a,b\}^{\otimes E}$,
then one can obtain bounds on the fluctuations of the random diameter $D=D_\omega$ of $(G,d_\omega)$. See \cite[Theorem 2]{\BKSfpp} for a precise statement in this more general context.
\end{remark}

\proof

For any edge $e$, let us consider the gradient along the edge 
$e$: $\nabla_e \Circ_m$. These gradient functions have values in $[-(b-a), b-a]$.
By dividing our distances by the constant factor $b-a$, we can even assume without 
loss of generality that our gradient functions have values in 
$[-1,1]$. Doing so, we end up being in a setup similar to the one we had 
in Chapter \ref{ch.hyper}. The {\bf influence} of an edge $e$ corresponds here to 
$\Inf_e(\Circ_m):=\Pb{\nabla_e \Circ_m(\omega) \neq 0}$.
We will prove later on that $\Circ_m$ has very small influences. 
In other words, we will show that the above gradient functions have small support, 
and hypercontractivity will imply the desired bound.
\vskip 0.2 cm

We have thus reduced the problem to the following general framework. Consider a real-valued 
function $f : \{-1,1\}^n \to \R$, such that for any variable $k$, 
$\nabla_k f \in [-1,1]$. We are interested in $\Var(f)$ and we
want to show that if ``influences are small'' then $\Var(f)$ is small.
It is easy to check that the variance can be written
\[
\Var(f) = \frac 1 4 \sum_k \sum_{\emptyset \neq S \subseteq [n]}  \frac 1 {|S|} \widehat{\nabla_k f} (S)^2\,.
\]
If all the variables have very small influence, then, as previously, 
$\nabla_k f$ should be of high frequency. Heuristically, this should then imply that 
\begin{eqnarray}
\Var(f) & \ll &   \sum_k \sum_{S\neq \emptyset} \widehat{\nabla_k f} (S)^2 \nonumber \\
& = & \sum_k \Inf_k(f) \,. \nonumber 
\end{eqnarray}

This intuition is quantified by the following lemma on the link between the fluctuations of a real-valued function $f$ on $\Omega_n$ and its influence vector. 

\begin{lemma}\label{l.torusfpp}
Let $f : \Omega_n \to \R$ be a (real-valued) function such that each
of its discrete derivatives $\nabla_k f, \, k\in [n]$ have their values in $[-1,1]$.
Let $\Inf_k(f):=\Pb{\nabla_k f \neq 0}$ be the influence of the $k$th bit.
Assume that the influences of $f$ are small in the sense that
there exists some $\alpha>0$ such that for any $k\in \{1, \ldots, n\}$, 
$\Inf_k(f) \le n^{-\alpha}$. Then there is some constant $C=C(\alpha)$ such that 
\[
\Var(f) \le \frac {C} {\log n} \sum_k \Inf_k(f)\,.
\]
\end{lemma}

\begin{remark}
If $f$ is Boolean, then this follows from Theorem \ref{th:KKL2}
with $C(\alpha)=c/\alpha$ with $c$ universal.
\end{remark}

The proof of this lemma is postponed to the next section. 
In the meantime, let us show that in our special case of first passage percolation on the torus, the assumption on small 
influences is indeed verified. Since the edge 
lengths are in $\{a,b\}$, the smallest contour $\Circ_m(\omega)$ in $\T_m^d$ around the first coordinate 
lies somewhere in $[am, bm]$. Hence, if $\gamma$ is a geodesic (a path in the torus with
the required winding number) satisfying $\mathrm{length}(\gamma)=\Circ_m(\omega)$,
then $\gamma$ uses at most $\frac b a m$ edges. There might be several different geodesics minimizing the circumference.
Let us choose randomly one of these in an ``invariant'' way and call it $\tilde \gamma$. For any edge $e\in E(\T_m^d)$, if,
by changing the length
of $e$, the circumference increases, then $e$ has to be contained in any geodesic $\gamma$, and in particular in $\tilde \gamma$. This implies that 
$\Pb{\nabla_e \Circ_m(\omega) > 0} \le \Pb{e \in \tilde \gamma }$. By symmetry we obtain that 
\[
\Inf_e(\Circ_m)=\Pb{\nabla_e \Circ_m(\omega) \neq 0} \le 2 \Pb{e \in  \tilde \gamma}\,.
\]

Now using the symmetries both of the torus $\T_m^d$ and of our observable $\Circ_m$, if $\tilde \gamma$ is chosen in an appropriate invariant way
(uniformly among all geodesics for instance), then it is clear that all 
the ``vertical'' edges (meaning those edges which, when projected onto the first 
coordinate, project onto a single vertex) have the same probability to 
lie in $\tilde\gamma$. The same is true for the ``horizontal'' edges. In particular
we have that
\[
\sum_{\text{``vertical'' edges }e} \Pb{e \in \tilde \gamma} \le \Eb{ |\tilde \gamma |} \le \frac b a m\,.
\]
Since there are at least order $m^d$ vertical edges, the influence of each of
these is bounded by $O(1)m^{1-d}$. The same is true for the horizontal edges.
All together this gives the desired assumption needed in 
Lemma \ref{l.torusfpp}. Applying this lemma, we indeed obtain that 
\[
\Var(\Circ_m(\omega)) \le O(1) \frac m {\log m}\,,
\]
where the constant does not depend on the dimension $d$; the dimension in fact 
helps us here, since it makes the influences smaller. \QED

\begin{remark}
At this point, we know that for any edge $e$, $\Inf_e(\Circ_m)=O(\frac m {m^d})$. 
Hence, at least in the case of the torus, one easily deduces from 
Poincar\'e's inequality the theorem by Kesten which says that 
$\Var(\Circ_m) = O(m)$.
\end{remark}

\section{Upper bounds on fluctuations in the spirit of KKL}

In this section, we prove Lemma \ref{l.torusfpp}.

\proof
Similarly as in the proofs of Chapter \ref{ch.hyper}, the proof relies on 
implementing hypercontractivity in the right way. We have that for any $c$,
\begin{eqnarray}
\var(f) & = & \frac 1 4 \sum_k \sum_{S\neq \emptyset} \frac 1 {|S|} \widehat{\nabla_k f} (S)^2 \nonumber\\
&\le& \frac 1 4 \sum_k \sum_{0 < |S| < c \log n} \widehat{\nabla_k f} (S)^2 + \frac {O(1)} {\log n} \sum_k \Inf_k(f) \nonumber
\end{eqnarray}
where the $O(1)$ term depends on the choice of $c$.

Hence it is enough to bound the contribution of small frequencies, $0<|S|< c \log n$, for some constant $c$ which will be chosen 
later. As previously we have for any $\rho\in (0,1)$ and using hypercontractivity,
\begin{eqnarray}
\sum_k \sum_{0 < |S| < c \log n} \widehat{\nabla_k f} (S)^2  &\le & \rho^{-2 c \log n} \sum_k \| T_\rho \nabla_k f \|_2^2 \nonumber \\
&\le  & \rho^{-2 c \log n} \sum_k \| \nabla_k f\|_{1+\rho^2} ^2 \nonumber \\
&\le & \rho^{-2 c \log n} \sum_k \Inf_k(f)^{2/(1+\rho^2)} \nonumber\\
&\le & \rho^{-2 c \log n} \bigl(\sup_{k} \Inf_k(f) \bigr)^{\frac{1-\rho^2}{1+\rho^2}} \sum_{k} \Inf_k(f) \nonumber \\ 
&\le& \rho^{-2 c \log n} n^{-\alpha \frac {1-\rho^2}{1+\rho^2} } \; \sum_k \Inf_k(f) \; \; \text{by our assumption}\,.\nonumber \\
\end{eqnarray}

Now fixing any $\rho\in(0,1)$, and then choosing the constant $c$ depending on 
$\rho$ and $\alpha$, the lemma follows. 
By optimizing on the choice of $\rho$, one could get better constants
if one wants. \qed

\section{Further discussion}

\ni
{\bf Some words on the proof of Theorem \ref{th.fpp}}
\vskip 0.3 cm

The main difficulty here is that the quantity of interest, $f(\omega):=\dist_\omega(0,v)$,
is no longer invariant under a large class of 
graph automorphisms. This lack of symmetry makes the study of influences more difficult. For example, edges 
near the endpoints 0 or $v$ have very high influence (of order one). To gain some more symmetry, the authors in \cite{\BKSfpp} rely on a very nice
``averaging'' procedure. We refer to this paper for more details.
\vskip 0.3 cm

\note{In detail}
{The idea is as follows: instead of looking at the (random) distance form 0 to $v$, they first pick a 
point $x$ randomly in the mesoscopic box $[-|v|^{1/4}, |v|^{1/4}]^d$ around the origin and then consider the distance from this
point $x$ towards $v+x$. Let $\tilde f$ denote this function ($\dist_\omega(x,v+x)$). $\tilde f$ uses extra randomness compared to $f$, but it is clear
that $\Eb{f}= \Eb{\tilde f}$ and it is not hard to see that when $|v|$ is large, $\var(f) \asymp \var(\tilde f)$. Therefore it is enough to study the fluctuations
of the more symmetric $\tilde f$ (we already see here that thanks to this avergaing procedure; the endpoints 0 and $v$ do not have anymore a high influence).
In some sense, along geodesics, this procedure {\bf ``spreads''} the influence on the $|v|^{1/4}$-neighborhood of the geodesics. More precisely, 
if $e$ is some edge, the influence of this edge is bounded by $2 \Pb{e \in x + \gamma}$, where $\gamma$ is chosen among geodesics from 0 to $v$.
Now, as we have seen in the case of the torus, geodesics are essentially one-dimensional (of length less than $O(1)|v|$); this is still true on the mesoscopic
scale: For any box $Q$ of radius $m:=|v|^{1/4}$, $|\gamma \cap Q| \le O(1) m$. Now by considering the mesoscopic box around $e$, it is like moving a ``line''
in a box of dimension $d$; the probability for an edge to be hit by that ``line'' is of order $m^{1-d}$. Therefore the influence of any edge $e$ for the 
``spread'' function $\tilde f$ is bounded by $O(1) |v|^{(1-d)/4} \le O(1)|v|^{-1/4}$.
This implies the needed assumption in Lemma \ref{l.torusfpp} and hence concludes 
the sketch of proof of Theorem \ref{th.torus}. See \cite{\BKSfpp} for a more detailed
proof.}

\ni
{\bf Known lower bounds on the fluctuations}
\vskip 0.3 cm

We discussed mainly here ways to obtain upper bounds on the
fluctuations of the shapes in first passage percolation.
It is worth pointing out that some non-trivial {\em lower} 
bounds on the fluctuations are
known for $\Z^2$. See \cite{\PemantlePeres, \NewmanPiza}.

\begin{remark}
We end by mentioning that the proof given in \cite{\BKSfpp} was based
on an inequality by Talagrand. The proof given here avoids this inequality.
\end{remark}

\vskip 0.4 cm
\note{Remark}{ Chr : Some polynomial lower bounds for the fluctuations are known in some 
specific directions (in $n^{1/8}$ it seems)!! }

\chapter*{Exercise sheet of Chapter \ref{ch.AF}}
\setcounter{exercise}{0}

\begin{problem}
Let $n\geq 1$ and $d\geq 2$.  Consider the random metric on the torus $\Z^d / n\Z^d$ as described in this chapter.
For any $k\geq 1$, let $\mathcal{A}_n^k$ be the event that the shortest ``horizontal'' circuit is $\leq k$. If $d\geq 3$, show that for any choice of $k_n=k(n)$, the family of events $\mathcal{A}_n^{k_n}$ is noise sensitive.  (Note that the situation here is similar to the Problem \ref{exer.clique} in Chapter \ref{ch.BF}.) Finally, discuss the 
two-dimensional case, $d=2$ (non-rigorously).
\end{problem}

\begin{exercise}
Show that Lemma \ref{l.torusfpp} is false if
$\Inf_k(f)$ is taken to be the square of the $L^2$ norm of $\nabla_k f$ 
rather than the probability of its support (i.e. find a counterexample). 
\end{exercise}

\note{Note that with usual definition of support, KKL not true for all fs since
variance scales with f and influences dont. So, clearly not true for large f but
true for small f.  If we take the L1 norm at the start of this exercise, still false?}


\note{Used to be an exercise}
{
Does the proof carries through for the model of {\bf last passage percolation} ? If not explain why. 
Can you modify the model of last passage-percolation so that something interesting can be said.
}

\chapter[\;\;\; Randomized algorithms and noise sensitivity]
{Randomized algorithms and noise sensitivity}\label{ch.RA}

\note{Timing}{//////////// 1.5 hour   ///////////// \vskip 0.5 cm}

In this chapter, we explain how the notion of {\bf revealment} for
so-called randomized algorithms can in some cases yield direct
information concerning the energy spectrum which may allow not only
noise sensitivity results but even quantitative noise sensitivity results.

\section{BKS and randomized algorithms}

In the previous chapter, we explained how Theorem \ref{th:NSmainresult}
together with bounds on the pivotal exponent for percolation yields
noise sensitivity for percolation crossings. However, 
in \cite{\BKS}, a different approach was in fact used for showing
noise sensitivity which, while still using Theorem \ref{th:NSmainresult},
did not use these bounds on the critical exponent. In that approach, one sees
the first appearance of randomized algorithms. In a nutshell,
the authors showed that (1) if a monotone function is very
uncorrelated with all majority functions, then it is noise sensitive
(in a precise quantitative sense) and (2) percolation crossings
are very uncorrelated with all majority functions. The latter
is shown by constructing a certain algorithm which,
due to the RSW Theorem \ref{th.RSW},
looks at very few bits but still looks at
enough bits in order to be able to determine the output of the function.

\section{The revealment theorem} \label{ss.reveal}

An {\bf algorithm}
for a Boolean function $f$ is an algorithm $A$ which  
queries (asks the values of) the bits one by one,
where the decision of which bit to ask can be based on the values
of the bits previously queried, and stops once $f$ is determined (being
determined means that $f$ takes the same value no matter how the remaining
bits are set).

A {\bf randomized algorithm}
for a Boolean function $f$ is the same as above but auxiliary
randomness may also be used to decide the next value queried
(including for the first bit).  [In computer science, the term
randomized decision tree would be used for our notion of
randomized algorithm, but we will not use this terminology.]

The following definition of {\it revealment} will be crucial.
Given a randomized algorithm $A$ for a Boolean function $f$,
we let $J_{A}$ denote the random set of bits queried by
$A$. (Note that this set
depends both on the randomness corresponding to the choice of $\omega$
and the randomness inherent in running the algorithm,
which are of course taken to be independent.)

\begin{defn}\label{d.revealment}
The {\bf revealment of a randomized algorithm} $A$ 
for a Boolean function $f$, denoted by $\delta_{A}$, is defined by
$$
\delta_{A}:=\max_{i\in \nn}\P(i\in J_{A}).
$$
The {\bf revealment of a Boolean function} $f$, 
denoted by $\delta_f$, is defined by
$$
\delta_f:=\inf_{A}\delta_{A}
$$
where the infimum is taken over all
randomized algorithms $A$ for $f$.
\end{defn}

This section presents a connection between noise sensitivity
and randomized algorithms. It will be used later to 
yield an alternative proof of noise sensitivity for 
percolation crossings which is not based upon Theorem
\ref{th:NSmainresult} (or Proposition \ref{pr.BKSquantitative}).
Two other advantages of the  algorithmic approach of the
present section over that mentioned in the previous section 
(besides the fact that it does not rest on Theorem \ref{th:NSmainresult})
is that it applies to
nonmonotone functions and yields a more ``quantitative'' version of
noise sensitivity.

We have only defined algorithms, randomized algorithms and revealment
for Boolean functions but the definitions immediately 
extend to functions $f:\Omega_n\to\R$. 

The main theorem of this section is the following.

\begin{theorem}[\cite{\SS}]\label{t.ss}
For any function $f:\Omega_n\to\R$ and 
for each $k=1,2,\dots$, we have that
\begin{equation}\label{e.noise}
E_f(k)=\sum_{S\subseteq\nn,\,|S|=k}\hat f(S)^2 \le \delta_f\,k\,\|f\|^2,
\end{equation}
where $\|f\|$ denotes the $L^2$ norm of $f$ with respect to the
uniform probability measure on $\Omega$ and 
$\delta_f$ is the revealment of $f$.
\end{theorem}

\medskip\noindent
Before giving the proof, we make some comments to help the reader see 
what is happening and suggest why a result like this might be true.
Our original function is
a sum of monomials with coefficients given by the Fourier coefficients.
Each time a bit is revealed by the algorithm, we obtain a new Boolean function
obtained by just substituting in the value of the bit we obtained into
the corresponding variable. On the algebraic side, those monomials which
contain this bit go down by 1 in degree while the other monomials
are unchanged. There might however be cancellation in the process which is 
what we hope for since when the algorithm stops, all the monomials (except the
constant) must have been killed. 
The way cancellation occurs is illustrated as follows.
The Boolean function at some stage might contain
$(1/3) x_2x_4x_5 +(1/3) x_2x_4$ and then the bit $x_5$ might be revealed and
take the value $-1$.  When we substitute this value into the variable,
the two terms cancel and disappear, thereby bringing us 1 step closer to a constant
(and hence determined) function. 

As far as why the result might be true, the intuition,
very roughly speaking, is as follows.  The theorem says that for
a Boolean function we cannot, for example, have $\delta=1/1000$ and 
$\sum_{i}\hat f(\{i\})^2 =1/2$. If the level 1 monomials of the function were
$$
a_1\omega_1 + a_2\omega_2 +\cdots + a_n\omega_n,
$$
then it is clear that after the algorithm is over, then with high
probability, the sum of the squares of the coefficients of the terms
which have not been reduced to a constant is still reasonably large. Therefore,
since the function at the end of the algorithm is constant, these remaining
terms must necessarily have been cancelled by higher degree monomials which,
after running the algorithm, have been ``reduced to'' degree 1 monomials.
If, for the sake of this heuristic
argument, we assume that each bit is revealed independently, then
the probability that a degree $k\ge 2$ monomial is brought down to
a degree 1 monomial (which is necessary for it to help to cancel
the degree 1 terms described above) is at most $\delta^{k-1}$
and hence the expected sum of the squares of the coefficients from the 
degree $k\ge 2$ monomials which are brought down to degree 1 is at
most $\delta^{k-1}$. The total such sum for levels 2 to $n$ is then at most
$$
\sum_{k=2}^n\delta^{k-1}\le 2\delta
$$
which won't be enough to cancel the (originally) degree 1 monomials
which remained degree 1 after running the algorithm if $\delta$ is much less
than $\sum_{i}\hat f(\{i\})^2$. A similar heuristic works for the other levels.

\proof
In the following, we let $\tilde{\Omega}$ denote the probability
space that includes the randomness in the input bits of $f$ and the 
randomness used to run the algorithm (which we assume to be
independent) and we let $\E$ denote the corresponding expectation.
Without loss of generality, elements of 
$\tilde{\Omega}$ can be represented as
$\tilde{\omega}=(\omega,\tau)$ where 
$\omega$ are the random bits and $\tau$ represents the
randomness necessary to run the algorithm.

Now, fix $k\ge 1$. Let 
$$
g(\omega):=\sum_{|S|=k} \hat f(S)\,\chi_S(\omega)\,,\qquad
\omega\in\Omega.
$$
The left hand side of (\ref{e.noise}) is equal to $\|g\|^2$.

Let $J\subseteq\nn$ be the random set of all bits examined by the algorithm.
Let $\ev A$ denote the 
minimal $\sigma$-field for which $J$ is measurable and every
$\omega_i$, $i\in J$, is measurable; this can be viewed as the
relevant information gathered by the algorithm.
For any function $h:\Omega\to\R$, let
$h_J:\Omega\to\R$ denote the random function obtained by substituting
the values of the bits in $J$. 
More precisely, if $\tilde{\omega}=(\omega,\tau)$
and $\omega'\in\Omega$, then $h_J(\tilde{\omega})(\omega')$ is
$h(\omega'')$ where $\omega''$ is $\omega$ on $J(\tilde{\omega})$ and is
$\omega'$ on $[n]\backslash J(\tilde{\omega})$.
In this way, $h_J$ is a random variable
on $\tilde{\Omega}$ taking values in the set of mappings from
$\Omega$ to $\R$ and it is immediate that this random variable is
$\ev A$-measurable.
When the algorithm terminates, the unexamined bits in $\Omega$
are unbiased and hence 
$\Eb{h\!\md\! \ev A}=\int h_J(={\hat h_J(\emptyset)})$
where $\int$ is defined, as usual, to be integration with
respect to uniform measure on $\Omega$.
It follows that $\Es{h}=\Es{\int h_J}$.

Similarly, for all $h$,
\begin{equation}\label{e.equiv}
\|h\|^2=\Eb{h^2}=\EB{\int h_J^2}=\Eb{\|h_J\|^2}.
\end{equation}

Since the algorithm determines  $f$, it is $\ev A$ measurable,
and we have
\begin{equation*}
\|g\|^2=\Es{g\,f}=\EB{\Eb{g\,f\!\md\! \ev A}}=
\EB{f\,\Eb{g\!\md\!\ev A}}.
\end{equation*}
Since $\Eb{g\!\md\! \ev A}={\hat g_J(\emptyset)}$,  Cauchy-Schwarz therefore gives
\begin{equation}\label{e.ept}
\|g\|^2 \le \sqrt{ \Es{\hat g_J(\emptyset)^2}}\,\|f\|\,.
\end{equation}

We now apply Parseval's formula to the (random) function $g_J$: 
this gives (for any $\tilde\omega=(\omega,\tau) \in \tilde\Omega$),
\[
\hat g_J(\emptyset)^2 = \|g_J\|_2^2 - \sum_{|S|>0} \hat g_J(S)^2 .
\]
Taking the expectation over $\tilde\omega\in \tilde \Omega$, this leads to 
\begin{eqnarray}
\Eb{\hat g_J(\emptyset)^2} & = & \Eb{\| g_J \|_2^2} - \sum_{|S|>0} \Eb{ \hat g_J(S)^2}  \nonumber \\
&=& \|g \|_2^2 - \sum_{|S|>0} \Eb{\hat g_J(S)^2} \text{  \; by
(\ref{e.equiv})} \nonumber \\
&=& \sum_{|S|=k} \hat g(S)^2 - \sum_{|S|>0} \Eb{\hat g_J(S)^2}  \left \lbrace \begin{array}{ll} \text{since $g$ is supported} \\ \text{on level-$k$ coefficients} \end{array} \right. \nonumber \\
&\le & \sum_{|S|=k} \Eb{\hat g(S)^2 - \hat g_J(S)^2} \left \lbrace \begin{array}{l} \text{by restricting to} \\ \text{level-$k$ coefficients} \end{array} \right. \nonumber
\end{eqnarray}

Now, since $g_J$ is built randomly from $g$ by fixing the variables in $J=J(\tilde\omega)$, and since $g$ by definition 
does not have frequencies larger than $k$, it is clear that for any
$S$ with $|S|=k$ we have
\begin{equation*}
\hat g_J(S) = \left\lbrace \begin{array}{ll}\hat g(S) = \hat f(S), & \text{if $S\cap J(\tilde \omega) =\emptyset$} \\ 0, & \text{otherwise.} \end{array} \right. 
\end{equation*}
Therefore, we obtain
\[
\| \Eb{g\!\bigm| \! J} \|_2^2 = \Eb{\hat g_J(\emptyset)^2} \le \sum_{|S|=k} \hat g(S)^2 \, \Pb{S\cap J \neq \emptyset} \le \|g\|_2^2 \, k \,\delta\,.
\]
Combining with (\ref{e.ept}) completes the proof.
\qed

\bigskip\noindent
Proposition \ref{pr.NSeq} and Theorem \ref{t.ss} immediately imply
the following corollary.

\begin{corollary}\label{cor.ss}
If the revealments satisfy 
$$
\lim_{n\to\infty}\delta_{f_n}= 0,
$$
then $\{f_n\}$ is noise sensitive.
\end{corollary}

\medskip\noindent
In the exercises, one is asked to show that certain sequences of
Boolean functions are noise sensitive by applying the above corollary.

\section{An application to noise sensitivity of percolation}

In this section, we apply Corollary \ref{cor.ss} to prove noise sensitivity
of percolation crossings. The following result gives the necessary assumption
that the revealments approach 0.

\begin{theorem}[\cite{\SS}]\label{t.square}
Let $f=f_n$ be the indicator function for the 
event that critical site percolation on the triangular grid
contains a left to right crossing of our $n\times n$ box.
Then $\delta_{f_n} \le n^{-1/4+o(1)}$
as $n\to\infty$.

For critical bond percolation on the square grid, this holds with
$1/4$ replaced by some positive constant $a>0$.
\end{theorem}

\proofoutline
We outline the argument only for the triangular lattice; the argument
for the square lattice is similar.
We first give a first attempt at a good algorithm. We consider
from Chapter \ref{ch.perc} the exploration path or 
interface from the bottom right of the square to
the top left used to detect a left right crossing.
This (deterministic) algorithm simply asks the bits
that it needs to know in order to continue the interface. Observe that
if a bit is queried, it is necessarily the case that there is both a
black and white path from next to the hexagon to the boundary. It follows,
from the exponent of $1/4$ for the 2-arm event in Chapter 
\ref{ch.perc}, that, for hexagons
far from the boundary, the probability that they are revealed is at
most $R^{-1/4+o(1)}$ as desired. However, one cannot conclude that
points near the boundary have small revealment and of course the right
bottom point is always revealed.

The way that we modify the above algorithm so that all points have
small revealment is as follows. We first choose a point $x$ at random from
the middle third of the right side. We then run two algorithms, the
first one which
checks whether there is a left right path from the right side {\it above}
$x$ to the left side and the second one which
checks whether there is a left right path from the right side {\it below}
$x$ to the left side. The first part is done by looking
at an interface from $x$ to the top left corner as above. 
The second part is done by looking
at an interface from $x$ to the bottom left corner as above (but
where the colors on the two sides of the interface need to be
swapped.)

It can then be shown with a little work (but no new
conceptual ideas) that this modified algorithm has the desired 
revealment of at most $R^{-1/4+o(1)}$ as desired. One of the things
that one needs to use in this analysis is the so-called one-arm half-plane
exponent, which has a known value of $1/3$. See \cite{\SS}
for details.  \qed

\ni
\begin{minipage}{0.5\textwidth}
\begin{center}
\includegraphics[width=0.7\textwidth]{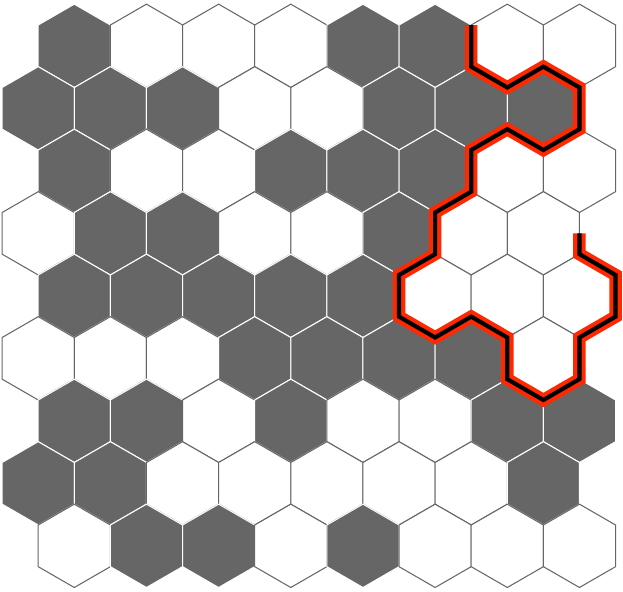}
\end{center}
\end{minipage}
\hskip 0.3 cm
\begin{minipage}{0.5 \textwidth}
\begin{center}
\includegraphics[width=0.7\textwidth]{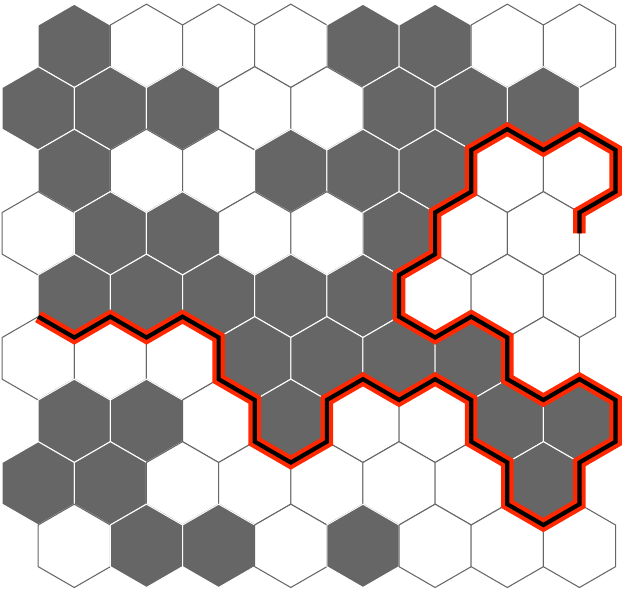}
\end{center}
\end{minipage}

\subsection {First quantitative noise sensitivity result}\label{ss.quant}

In this subsection, we give our first ``polynomial bound'' on the noise 
sensitivity of percolation. This is an important step in our understanding
of quantitative noise sensitivity of percolation initiated in 
Chapter \ref{ch.FE}.

Recall that in the definition of noise sensitivity, $\eps$ is held fixed. 
However, as we have seen in Chapter \ref{ch.FE}, it is of interest 
to ask if the correlations can still go to 0 when 
$\eps=\eps_n$ goes to 0 with $n$ but not so fast. The techniques of the 
present chapter imply the following result.

\begin{theorem}[\cite{\SS}] \label{th:crossingquant} Let $\{f_n\}$ be as in 
Theorem \ref{t.square}. Then, for the triangular lattice,
for all $\gamma < 1/8$,
\begin{equation}
\lim_{n\to\infty} \E[f_n(\omega)f_n(\omega_{1/n^\gamma})]- \E[f_n(\omega)]^2 =0.
\end{equation}
On the square lattice, there exists some $\gamma>0$ with the above property.
\end{theorem}

\proof
We prove only the first statement; the square lattice case is handled
similarly. First, (\ref{e.mainequation}) gives us that
every $n$ and $\gamma$,
\begin{equation} \label{e.again}
\E[f_n(\omega)f_n(\omega_{1/n^\gamma})]- \E[f_n(\omega)]^2 =
\sum_{k=1}E_{f_n}(k) (1-1/n^\gamma)^k.
\end{equation} 
Note that there are order $n^2$ terms in the sum.
Fix $\gamma < 1/8$. Choose $\epsilon >0$ so that 
$\gamma+\eps < 1/8$. For large $n$, we have
that $\delta_{f_n} \le 1/n^{1/4-\eps}$. 
The right hand side of (\ref{e.again}) is at most 
$$
\sum_{k=1}^{n^{\gamma+\eps/2}} k/n^{1/4-\eps}+
(1-1/n^{\gamma})^{n^{\gamma+\eps/2}}
$$
by breaking up the sum at $n^{\gamma+\eps/2}$ and applying
Theorems \ref{t.ss} and \ref{t.square}
to bound the $E_{f_n}(k)$ terms in the first part.
The second term clearly goes to 0 while the first part also goes to 0
by the way $\eps$ was chosen.
\qed

Observe that the {\it proof} of Theorem \ref{th:crossingquant} immediately
yields the following general result.

\begin{corollary} Let $\{f_n\}$ be a sequence of Boolean functions on $m_n$
bits with $\delta(f_n)\le O(1)/n^\beta$ for all $n$. Then 
for all $\gamma < \beta/2$, we have that
\begin{equation}
\lim_{n\to\infty} \E[f_n(\omega)f_n(\omega_{1/n^\gamma})]- \E[f_n(\omega)]^2 =0.
\end{equation}
\end{corollary}

\section{Lower bounds on revealments}
One of the goals of the present section is to show that
one cannot hope to reach the conjectured $3/4$-sensitivity exponent
with Theorem \ref{t.ss}.
Theorem \ref{th:crossingquant} told us that we obtain asymptotic
decorrelation if the noise is $1/n^\gamma$ for $\gamma< 1/8$. Note that this
differs from the conjectured ``critical exponent'' of $3/4$
by a factor of 6. In this section, we investigate the degree to which
the $1/8$ could potentially be improved and in the discussion, we
will bring up an interesting open problem. We will also derive an
interesting general theorem giving a nontrivial 
lower bound on the revealment for monotone functions. We start with
the following definition.

\begin{defn}\label{d.cost}
Given a randomized algorithm $A$ for a Boolean function $f$, let
{\bf $C(A)$} (the cost of $A$) be the expected number of queries that the 
algorithm $A$ makes. Let {\bf $C(f)$} (the cost of $f$) be the infimum of 
$C(A)$ over all randomized algorithms $A$ for $f$.
\end{defn}

\begin{remark}
\ni
(i). It is easy to see that $C(f)$ is unchanged if we take the infimum 
over deterministic algorithms. \\
(ii). Clearly $n\delta_A\ge C(A)$ and hence $n\delta_f\ge C(f)$. \\
(iii). $C(f)$ is at least the total influence $\Inf(f)$ since for any
algorithm $A$ and any $i$, the event that $i$ is
pivotal necessarily implies that the bit $i$ is queried by $A$.
\end{remark}

The following result due to O'Donnell and Servedio 
(\cite{\OdonnellServedio})is an essential improvement on the third
part of the last remark.

\begin{theorem} \label{t.OS} Let 
$f$ be a monotone Boolean function mapping $\Omega_n$ into $\{-1,1\}$.
Then $C(f)\ge \Inf(f)^2$ and hence $\delta_f\ge \Inf(f)^2/n$.
\end{theorem} 

\proof 
Fix any randomized algorithm $A$ for $f$.
Let $J=J_A$ be the random set of bits queried by $A$. We then have
$$
\Inf(f)=\E[\sum_i f(\omega)\omega_i]=
\E[f(\omega)\sum_i \omega_iI_{\{i\in J\}}]\le
\sqrt{\E[f(\omega)^2]}
\sqrt{\E[(\sum_i \omega_iI_{\{i\in J\}})^2]}
$$
where the first equality uses monotonicity
(recall Proposition \ref{pr.inf&spectrum.monotone}) and then the
Cauchy-Schwarz inequality is used.
We now bound the first term by 1. For the second moment inside the second
square root, the sum of the diagonal terms yields $\E[|J|]$ while the
cross terms are all 0 since for $i\neq j$, 
$\E[\omega_iI_{\{i\in J\}}\omega_jI_{\{j\in J\}}]=0$ as can be
seen by breaking up the sum depending on whether $i$ or $j$ is 
queried first. This yields the result.
\qed

Returning to our event $f_n$ of percolation crossings, 
since the sum of the influences is
$n^{3/4+o(1)}$, Theorem \ref{t.OS} tells us that
$\delta_{f_n} \ge n^{-1/2+o(1)}$. It follows from the method
of proof in Theorem \ref{th:crossingquant} that Theorem \ref{t.ss}
cannot improve the result of
Theorem \ref{th:crossingquant} past $\gamma=1/4$ which
is still a factor of 3 from the critical value $3/4$. Of course,
one could investigate the degree to which Theorem
\ref{t.ss} itself could be improved. 

Theorem \ref{t.square} tells us that there are algorithms $A_n$ for
$f_n$ such that $C(A_n)\le n^{7/4+o(1)}$. On the other hand,
Theorem \ref{t.OS} tell us that it is necessarily the case that
$C(A)\ge n^{6/4+o(1)}$.

\medskip\noindent
Open Question: 
Find the smallest $\sigma$ such that there are algorithms $A_n$ for
$f_n$ with $C(A_n)\le n^{\sigma}$. (We know $\sigma\in [6/4,7/4]$.)

\medskip
We mention another inequality relating revealment with influences 
which is a consequence of the results in \cite{\DecisionTrees}.

\begin{theorem} \label{t.OSSS} Let 
$f$ be a Boolean function mapping $\Omega_n$ into $\{-1,1\}$. Then
$\delta_f\ge \Var(f)/(n\max_i\Inf_i(f))$
\end{theorem} 

It is interesting to compare Theorems 
\ref{t.OS} and \ref{t.OSSS}. Assuming $\Var(f)$ is of order 1,
and all the influences are of order $1/n^\alpha$, then it is
easy to check that Theorem \ref{t.OS} gives a better bound
when $\alpha < 2/3$ and Theorem \ref{t.OSSS} gives a better bound
when $\alpha > 2/3$. For crossings of percolation, where $\alpha$
should be $5/8$, it is better to use Theorem
\ref{t.OS} rather than \ref{t.OSSS}. 

Finally, there are a number of interesting results concerning revealment
obtained in the paper \cite{\BalancedBoolean}. Four results are as follows.

\noindent
1. If $f$ is reasonably balanced on $n$ bits,
then the revealment is at least of order $1/n^{1/2}$. 

\noindent
2. There is a reasonably balanced function on $n$ bits
whose revealment is at most $O(1)(\log n)/n^{1/2}$. 

\noindent
3. If $f$ is reasonably balanced on $n$ bits and is monotone, 
then the revealment is at least of order $1/n^{1/3}$. 

\noindent
4. There is a reasonably balanced monotone function on $n$ bits 
whose revealment is at most $O(1)(\log n)/n^{1/3}$. 

We finally end this section by giving one more reference which gives an
interesting connection between percolation, algorithms and game theory;
see \cite{\RandomTurnHex}.

\section{An application to a critical exponent}

In this section, we show how Theorem \ref{t.ss} or in fact Theorem \ref{t.OS} 
can be used to show that the 4-arm exponent is strictly larger
than 1; recall that with SLE technology, this can be shown for
the triangular lattice.

\begin{proposition} \label{prop.4armexponent}
Both on the triangular lattice $\T$ and on $\Z^2$, there exists
$\epsilon_0 > 0$ such that
$$
\alpha_4(R) \le 1/R^{1+\epsilon_0}
$$
\end{proposition}

We will assume the separation of arms result mentioned earlier 
in Chapter \ref{ch.FE} which says that for the 
event $f_R$, the influence of any variable further than distance $R/10$ 
from the boundary, a set of variables that we will denote by $B$
for bulk, is $\asymp \alpha_4(R)$.

\proof 
Theorems \ref{t.square} and \ref{t.ss} imply that for some $a>0$,
$$
\sum_{i}\hat f_R(\{i\})^2 \le 1/R^a.
$$
Next, using the separation of arms as explained above, we have
\begin{equation}\label{e.armssep}
R^2\alpha^2_4(R) \le O(1) \sum_{i\in B} \Inf^2_i.
\end{equation}
Proposition \ref{pr.inf&spectrum.monotone} then yields
$$
R^2\alpha^2_4(R) \le O(1/R^a)
$$
and the result follows.
\qed

Observe that Theorem \ref{t.OS} could also be used as follows. 
Theorem \ref{t.square} implies that $C(f_R)\le R^{2-a}$ for some $a>0$
and then Theorem \ref{t.OS} yields $\Inf(f_R)^2\le R^{2-a}$.

Exactly as in (\ref{e.armssep}), one has, again using separation of arms, that
\begin{equation}
R^2\alpha_4(R) \le O(1) \sum_{i\in B} \Inf_i \le O(1) \Inf(f_R).
\end{equation}
Altogether this gives us
$$
R^4\alpha^2_4(R) \le O(1)R^{2-a},
$$
again yielding the result.

We finally mention that it is not so strange that either
of Theorems \ref{t.ss} or \ref{t.OS} can be used here since, as the reader
can easily verify, for the case of monotone functions all of whose variables
have the same influence, the case $k=1$ in Theorem \ref{t.ss} is
equivalent to Theorem \ref{t.OS}.

\begin{remark}
We now mention that the proof for the multi-scale version of
Proposition \ref{pr.upperboundalpha4} is an extension of
the approach of O'Donnell and Servedio above.
\end{remark}

\section{Does noise sensitivity imply low revealment?}

As far as these lectures are concerned, this subsection will not
connect to anything that follows and hence can be viewed as tangential.

It is natural to ask if the converse of Corollary \ref{cor.ss} might
be true. A moment's thought reveals that example \ref{ex.par}, 
Parity, provides a counterexample. However, it is more interesting
perhaps that there is a monotone counterexample to the converse which
is provided by example \ref{ex.clique}, Clique containment.

\begin{proposition}\label{pr.friedgutkahnwigderson}
Clique containment provides an example showing that the converse of 
Corollary \ref{cor.ss} is false for monotone functions.
\end{proposition}

\proofoutline
We first explain more precisely the size of the clique that we
are looking for. Given $n$ and $k$, let $f(n,k):=
\binom{n}{k}2^{-\binom{k}{2}}$, which is just the expected number of cliques
of size $k$ in a random graph. When $k$ is around $2\log_2(n)$, it is easy
to check that $f(n,k+1)/f(n,k)$ is $o(1)$ as $n\to\infty$. For such $k$,
clearly if  $f(n,k)$ is small, then with high probability there is no $k$-clique
while it can be shown, via a second moment type argument, that if 
$f(n,k)$ is large, then with high probability there is a $k$-clique.
One now takes $k_n$ to be around $2\log_2(n)$ such that
$f(n,k_n)\ge 1$ and $f(n,k_n+1)< 1$. Since 
$f(n,k+1)/f(n,k)$ is $o(1)$, it follows with some thought from the above that
the clique number is concentrated on at most 2 points. Furthermore, if 
$f(n,k_n)$ is very large and $f(n,k_n+1)$ very small, then it is concentrated
on one point. Again, see \cite{\AS} for details.

Finally, we denote the event that the random graph on $n$ vertices
contains a clique of size $k_n$ by $A_n$. 
We have already seen in one of the exercises that this example is
noise sensitive. We will only consider a sequence of $n$'s so that
$A_n$ is nondegenerate in the sense that the probabilities of this sequence
stay bounded away from 0 and 1. An interesting point is that there
is such a sequence. Again, see \cite{\AS} for this.
To show that the revealments do not go to
0, it suffices to show that the sequence of costs (see Definition 
\ref{d.cost} and the remarks afterwards) is $\Omega(n^2)$. We prove
something stronger but, to do this, we must first give a few more 
definitions.

\begin{definition}
For a given Boolean function $f$, a {\bf witness} for $\omega$
is any subset $W$ of the variables such that the elements of 
$\omega$ in $W$ determine $f$ in the sense
that for every $\omega'$ which agrees with $\omega$ on $W$, we have
that $f(\omega)=f(\omega')$. The 
{\bf witness size} of $\omega$, denoted $w(\omega)$, is the size of the
smallest witness for $\omega$. The
{\bf expected witness size}, denoted by $w(f)$, is $\E(w(\omega))$.
\end{definition}

Observe that, for any Boolean function $f$, the bits revealed by any
algorithm $A$ for $f$ and for any $\omega$ is always a witness for
$\omega$. It easily follows that the cost $C(f)$ satisfies
$C(f)\ge  w(f)$. Therefore, in order to
prove the proposition, it suffices to show that
\begin{equation}\label{e.kahn}
w(f_n) =\Omega(n^2).
\end{equation}

\begin{remark}
\ni
(i). The above also implies that with a fixed uniform probability,
$w(\omega)$ is $\Omega(n^2)$. \\
(ii). Of course when $f_n$ is $1$, there is always a (small) witness
of size $\binom{k_n}{2}\ll n$
and so the large average witness size comes from when
$f_n$ is $-1$.\\
(iii). However, it is not deterministically true that when
$f_n$ is $-1$, $w(\omega)$ is necessarily of size $\Omega(n^2)$.
For example, for $\omega\equiv -1$ (corresponding to the empty graph), 
the witness size is $o(n^2)$ as is easily checked. Clearly 
the empty graph has the smallest witness size among $\omega$ with $f_n=-1$.
\end{remark}

\begin{lemma}\label{lem.janson}
Let $E_n$ be the event that all sets of vertices of size
at least $.97n$ contains $C_{k_n-3}$. Then
$\lim_{n\to\infty}\P(E_n)=1$.
\end{lemma}

\proof
This follows, after some work, from the Janson inequalities.
See \cite{\AS} for details concerning these inequalities.
\qed

\begin{lemma}\label{lem.3goodvertices}
Let $U$ be any collection of at most $n^2/1000$ edges
in $C_n$. Then there exist distinct $v_1,v_2,v_3$ such that 
no edge in $U$ goes between any $v_i$ and $v_j$ and 
\begin{equation}\label{e.secondcondition}
|\{e\in U: e \mbox{ is an edge between $\{v_1,v_2,v_3\}$ and } 
\{v_1,v_2,v_3\}^c\}|\le  n/50.
\end{equation}
\end{lemma}

\proof 
We use the probabilistic method where we choose 
$\{v_1,v_2,v_3\}$ to be a uniformly chosen 3-set. It is immediate
that the probability that the first condition fails is at most 
$3|U|/\binom{n}{2}\le 1/100$. Letting $Y$ be the number of edges
in the set appearing in (\ref{e.secondcondition})
and $Y'$ be the number of $U$ edges touching $v_1$,
it is easy to see that
$$
\E(Y)\le 3\E(Y')=6|U|/n\le n/100
$$
where the equality follows from the fact that, for any graph,
the number of edges is half the total degree. By Markov's inequality,
the probability of the event in (\ref{e.secondcondition})
holds with probably at least 1/2.  This shows that the random 3-set 
$\{v_1,v_2,v_3\}$ satisfies the two stated conditions with positive probability
and hence such a 3-set exists. 
\qed

\medskip
By Lemma \ref{lem.janson}, we have $\P(A_n^c\cap E_n)\ge c>0$ for all large $n$.
To prove the theorem, it therefore suffices to show that if
$A_n^c\cap E_n$ occurs, there is no witness of size smaller than
$n^2/1000$.  Assume $U$ to be any set of edges of size smaller than
$n^2/1000$. Choose $\{v_1,v_2,v_3\}$ from Lemma \ref{lem.3goodvertices}.
By the second condition in this lemma,
there exists a set $S$ of size at least $.97n$ which is disjoint from
$\{v_1,v_2,v_3\}$ which has no $U$-edge to $\{v_1,v_2,v_3\}$.
Since $E_n$ occurs, $S$ contains a $C_{k_n-3}$, whose vertices
we denote by $T$.  Since there are no $U$-edges between
$T$ and $\{v_1,v_2,v_3\}$ or within $\{v_1,v_2,v_3\}$ (by the first
condition in Lemma \ref{lem.3goodvertices})
and $T$ is the complete graph, $U$ cannot be a witness since
$A_n^c$ occured.  \qed

\bigskip
The key step in the proof of Proposition \ref{pr.friedgutkahnwigderson} 
is (\ref{e.kahn}). This is stated without proof in \cite{\FKW}; 
however, E. Friedgut provided us with the above proof.

\chapter*{Exercise sheet of Chapter \ref{ch.RA}}

\setcounter{exercise}{0}

\begin{exercise}\label{exer.distcomplexity}
Compute the revealment for Majority function on 3 bits.
\end{exercise}

\begin{exercise}
Use Corollary \ref{cor.ss} to show that 
Examples \ref{ex.it3maj} and \ref{ex.tribes},
Iterated 3-Majority function and tribes, are noise sensitive.
\end{exercise}

\begin{exercise}
For transitive monotone functions, is there a relationship between
revealment and the minimal cost over all algorithms? 
\end{exercise}

\begin{exercise}\label{ex.SS=OS}
Show that for transitive monotone functions, Theorem \ref{t.OS}
yields the same result as Theorem \ref{t.ss} does for the case $k=1$.
\end{exercise}

\begin{exercise}
What can you say about the sequence of revealments for the
Iterated 3-Majority function?
[It can be shown that the sequence of revealments decays like $1/n^\sigma$
for some $\sigma$ but it is an open question what $\sigma$ is.]
\end{exercise}

\begin{exercise}
You are given a sequence of Boolean functions and told that it is not
noise sensitive using noise $\epsilon_n=1/n^{1/5}$. What, if anything,
can you conclude about the sequence of revealments $\delta_n$?
\end{exercise}

\begin{exercise}
Note that a consequence of Corollary \ref{cor.ss}
and the last line in Remark \ref{rem:BKScomments}
is that if $\{f_n\}$ is a sequence of monotone functions,
then, if the revealments of $\{f_n\}$ 
go to 0, the sums of the squared influences approach 0. 
Show that this implication is false without the monotonicity assumption.
\end{exercise}

\chapter{The spectral sample}\label{ch.SM}

\note{Timing}{///////// 45 minutes  ///////////}

It turns out that it is very useful to view the Fourier coefficients of a Boolean function 
as a random subset of the input bits where the ``weight'' or
``probability'' of a subset is its squared Fourier coefficient.
It is our understanding that it was Gil Kalai who suggested that
thinking of the spectrum as a random set could shed some light
on the types of questions we are looking at here.
The following is the crucial definition in this chapter.

\section{Definition of the spectral sample}

\begin{definition}
Given a Boolean function $f:\Omega_n \rightarrow \{\pm 1\}$ or
$\{0,1\}$, we let the {\bf spectral measure} $\SQ=\SQ_f$
of $f$ be the measure on subsets $\{1,\ldots,n\}$ given by
\[
\SQ_f(S):=\hat f(S)^2,\; S \subset \{1,\ldots, n \}\,.
\]

We let $\Spec_f=\Spec$ denote a subset of $\{1,\ldots,n\}$ chosen
according to this measure and call this the {\bf spectral sample}.
We let $\SQ$ also denote the corresponding expectation (even when
$\SQ$ is not a probability measure). 
\end{definition}
By Parseval's formula, the total mass of the so-defined spectral measure is 
\[
\sum_{S \subset \{ 1,\ldots, n\} } \hat f(S)^2 = \Eb{f^2} \,.
\]
This makes the following definition natural.
\begin{definition}
Given a Boolean function $f :\Omega_n \rightarrow \{\pm 1\}$ or $\{0 ,1 \}$,
we let the {\bf spectral probability measure} $\SP=\SP_f$ of $f$ be the probability measure 
on subsets of $\{1,\ldots, n\}$ given by 
\[
\SP_{f}(S):=\frac{\hat f(S)^2}{\E[f^2]}, \; S \subset \{1,\ldots, n \}\,.
\]
Since $\SP_f$ is just $\SQ_f$ up to a renormalization factor, the {\bf spectral sample}
$\Spec_{f} = \Spec$ will denote as well a random subset of $[n]$
sampled according to $\SP_f$. We let $\hat \E_f = \hat \E$ denote its corresponding expectation. 
\end{definition}


\begin{remark}\label{r.coco}
\leavevmode
\bi
\item[(i)] Note that if $f$ maps into $\{\pm 1\}$, then, by Parseval's
formula, $\SQ_f = \SP_f$  while
if it maps into $\{0,1\}$,  $\SQ_f$ will be a subprobability measure.

\item[(ii)] Observe that if $(f_n)_n$ is a sequence of non-degenerate  Boolean functions into $\{0,1\}$, then $\SP_{f_n}\asymp\SQ_{f_n}$. 
\item[(iii)] There is no statistical relationship between $\omega$ and $\Spec_f$ as
they are defined on different probability spaces. The spectral 
sample will just be a convenient point of view in order to understand
the questions we are studying. 
\ei
\end{remark}

Some of the formulas and results we have previously derived in these notes
have very simple formulations in terms of the spectral sample.
For example, it is immediate to check that
(\ref{e.correlationFourier}) simply becomes
\begin{equation} \label{eq:expansionrrandom}
\E[f(\omega)f(\omega_\epsilon)]= \SQ_f[(1-\epsilon)^{|\Spec|}]
\end{equation}
or
\begin{equation} \label{eq:expansionrrandomother}
\E[f(\omega)f(\omega_\epsilon)]-\E[f(\omega)]^2
= 
\SQ_f[(1-\epsilon)^{|\Spec|}I_{{\Spec}\neq \emptyset}].
\end{equation}

Next, in terms of the spectral sample, 
Propositions \ref{pr.NSeq} and \ref{pr.NStableeq} simply become the 
following proposition.

\begin{prop} \label{pr:NSNSAsample} 
If $\{f_n\}$ is a sequence of Boolean functions mapping into $\{\pm 1\}$, then we have
the following. \\
1. $\{f_n\}$ is noise sensitive if and only if 
$|\Spec_{f_n}|\to\infty$ in probability on the set
$\{|\Spec_{f_n}|\neq 0\}$.\\
2. $\{f_n\}$ is noise stable if and only if the random variables
$\{|\Spec_{f_n}|\}$ are tight.
\end{prop} 

There is also a nice relationship between the pivotal set 
$\calP$ and the spectral sample. The following result, which
is simply Proposition \ref{pr.inf&spectrum} (see also the remark
after this proposition), tells us that
the two random sets $\calP$ and $\Spec$ have 
the same 1-dimensional marginals.

\begin{prop} \label{pr:same1dmarg} 
If $f$ is a Boolean function mapping into $\{\pm 1\}$, then
for all $i\in [n]$ we have that
$$
\P(i\in \calP)= \SQ(i\in \Spec)
$$
and hence $\E(|\calP|)=\SQ(|\Spec|)$.
\end{prop} 

(This proposition is stated with $\SQ$ instead of $\SP$ since if $f$ maps into $\{0,1\}$ instead, 
then the reader can check that the above holds with an extra factor of 4 on the right hand side while
if $\SP$ were used instead, then this would not be true for any constant.)
Even though $\Spec$ and $\Piv$ have the same
``1-dimensional'' marginals, it is not however true that these two
random sets have the same distribution. For example, it is easily checked that
for $\MAJ_3$, these two distributions are different. Interestingly,
as we will see in the next section, $\Spec$ and $\Piv$ also always 
have the same ``2-dimensional'' marginals. This will prove useful when
applying second moment method arguments.

Before ending this section, let us give an alternative proof of
Proposition \ref{prop.uppertail}
using this point of view of thinking of $\Spec$ as a random set.

\medskip\ni
{\bf Alternative proof of Proposition \ref{prop.uppertail}}
The statement of the proposition when converted to the spectrum states
(see the exercises in this chapter if this is not clear)
that for any $a_n\to\infty$,
$$
\lim_{n\to\infty}\SP(|\Spec_n| \ge a_n n^2 \alpha_4(n))=0.
$$
However this immediately follows from Markov's inequality using 
Propositions \ref{pr.UpperboundInf} and \ref{pr:same1dmarg}.
\qed

\section{A way to sample the spectral sample in a sub-domain}

In this section, we describe a method of ``sampling'' the spectral
measure restricted to a subset of the bits. As an application
of this, we show that $\Spec$ and $\Piv$ in fact have the same
2-dimensional marginals, namely that for all $i$ and $j$,
$\P(i,j\in \calP)= \SQ(i,j\in \Spec)$.

In order to first get a little intuition about the spectral measure,
we start with an easy proposition.

\begin{proposition}[\cite{\GPS}]\label{pr.easy.spectral}
For a Boolean function $f$ and $A\subseteq \{1,2,\ldots,n\}$, we have
$$
\SQ(\Spec_f\subseteq A) =\E[|\E(f|A)|^2]
$$
where conditioning on $A$ means conditioning on the bits in $A$.
\end{proposition}

\proof
Noting that $\E(\chi_S|A)$ is $\chi_S$ if $S\subseteq A$ and 0 otherwise,
we obtain by expanding that
$$
\E(f|A)=\sum_{S\subseteq A } \hat f(S)\, \chi_S.
$$
Now apply Parseval's formula.
\qed

\medskip
If we have a subset $A\subseteq \{1,2,\ldots,n\}$, how do we ``sample''
from $A\cap \Spec$? A nice way to proceed is as follows: choose a random
configuration outside of $A$, then look at the induced function on $A$
and sample from the induced function's spectral measure. The
following proposition justifies in precise terms
this way of sampling. Its proof is
just an extension of the proof of Proposition \ref{pr.easy.spectral}.

\begin{proposition}[\cite{\GPS}]\label{pr.spectralsubdomain}
Fix a Boolean function $f$ on $\Omega_n$.
For $A\subseteq \{1,2,\ldots,n\}$ and $y\in\{\pm 1\}^{A^c}$, 
that is a configuration on $A^c$, let $g_y$
be the function defined on $\{\pm 1\}^{A}$ obtained by using $f$ but
fixing the configuration to be $y$ outside of $A$. Then for any
$S\subseteq A$, we have 
$$
\SQ(\Spec_f\cap A=S)=\E[\SQ(\Spec_{g_y}=S)]=\E[\hat g^2_y(S)].
$$
\end{proposition}

\proof
Using the first line of the proof of 
Proposition \ref{pr.easy.spectral}, it is easy to check that
for any $S\subseteq A$, we have that
$$
\Eb{f\,\chi_S\!\md\!\ev F_{A^c}}=\sum_{S'\subseteq A^c}\widehat f(S\cup S')\,\chi_{S'}\,.
$$
This gives 
$$
\EB{\Eb{f\,\chi_S\!\md\!\ev F_{A^c}}^2}=\sum_{S'\subseteq A^c}\widehat f(S\cup S')^2
=\SQ[\Spec\cap A=S]
$$
which is precisely the claim.
\qed

\begin{remark}
Observe that Proposition \ref{pr.easy.spectral} is a special case of
Proposition \ref{pr.spectralsubdomain} when $S$ is taken to be $\emptyset$
and $A$ is replaced by $A^c$.
\end{remark}

The following corollary was first observed by Gil Kalai.

\begin{corollary}[\cite{\GPS}]\label{cor.2dmarg}
If $f$ is a Boolean function mapping into $\{\pm 1\}$, then
for all $i$ and $j$,
$$
\P(i,j\in \calP)= \SQ(i,j\in \Spec).
$$
\end{corollary}
(The comment immediately following Proposition \ref{pr:same1dmarg} holds here as well.)

\proof
Although it has already been established that $\calP$ and $\Spec$
have the same 1-dimensional marginals, we first show how
Proposition \ref{pr.spectralsubdomain} can be used to establish this.
This latter proposition yields, with $A=S=\{i\}$, that
$$
\SQ(i\in \Spec)=\SQ(\Spec\cap \{i\} =\{i\})=
\E[\hat g^2_y(\{i\})].
$$
Note that $g_y$ is $\pm \omega_i$ if $i$ is pivotal and constant if
$i$ is not pivotal. Hence the last term is $\P(i\in \calP)$.

For the 2-dimensional marginals, one first checks this by hand when $n=2$.
For general $n$, taking $A=S=\{i,j\}$ in 
Proposition \ref{pr.spectralsubdomain}, we have 
$$
\SQ(i,j\in \Spec)=\P(\Spec\cap \{i,j\} =\{i,j\})=
\E[\hat g^2_y(\{i,j\})].
$$
For fixed $y$, the $n=2$ case tells us that 
$\hat g^2_y(\{i,j\})=\P(i,j\in \calP_{g_y})$. Finally, a little thought
shows that 
$\E[\P(i,j\in \calP_{g_y})]=\P(i,j\in \calP)$, completing the proof.
\qed

\section{Nontrivial spectrum near the upper bound for percolation}

We now return to our central event of percolation crossings
of the rectangle $R_n$ where $f_n$ denotes this event.
At this point, we know that for $\Z^2$, (most of) the spectrum
lies between $n^{\epsilon_0}$ (for some $\epsilon_0>0$)
and  $n^2 \alpha_4(n)$ while for 
$\T$ it sits between $n^{1/8+o(1)}$ and $n^{3/4+o(1)}$.
In this section, we show that there is a nontrivial amount of spectrum
near the upper bound $n^2 \alpha_4(n)$. For $\T$, in terms of quantitative noise
sensitivity, this tells us that if our noise sequence $\epsilon_n$ is
equal to $1/n^{3/4-\delta}$ for fixed $\delta>0$, then in the limit,
the two variables $f(\omega)$
and $f(\omega_{\epsilon_n})$ are not perfectly correlated; i.e.,
there is some degree of independence. (See the exercises for
understanding such arguments.)
However, we cannot conclude that
there is full independence since we don't know that ``all'' of the
spectrum is near $n^{3/4+o(1)}$ (yet!).

\begin{theorem}[\cite{\GPS}] \label{t.some.upper.spectrum}
Consider our percolation crossing functions $\{f_n\}$ (with values into  $\{ \pm 1 \}$) of the rectangles $R_n$
for $\Z^2$ or $\T$. There exists $c>0$ such that for all $n$,
$$
\SPb{|\Spec_n| \ge c n^2 \alpha_4(n)}\ge c.
$$
\end{theorem} 

The key lemma for proving this is the following second moment
bound on the number of pivotals which we prove afterwards.
It has a similar flavor to Exercise 6 in Chapter
\ref{ch.FE}.

\begin{lemma}[\cite{\GPS}]\label{l.pivotals.second.moment}
Consider our percolation crossing functions $\{f_n\}$ above and
let $R_n'$ be the box concentric with $R_n$ with half the radius.
If $X_n=|\calP_n\cap R_n'|$ is the cardinality of the set of pivotal points
in $R_n'$, then there exists a constant
$C$ such that for all $n$ we have that
$$
\E[|X_n|^2]\le C \E[|X_n|]^2.
$$
\end{lemma}

{\medbreak \noindent {\bf Proof of Theorem \ref{t.some.upper.spectrum}.}}
Since $\calP_n$ and $\Spec_n$ have the same 1 and 2-dimensional
marginals, it follows fairly straightforward from 
Lemma \ref{l.pivotals.second.moment} that we also have that for all $n$
$$
\SPb{ |\Spec_n\cap R_n'|^2} \le C \SPb{|\Spec_n\cap R_n'|}^2.
$$
Recall now the Paley-Zygmund inequality
which states that if $Z\ge 0$, then for all $\theta\in (0,1)$, 
\begin{equation}\label{e.PZ}
\P(Z\ge \theta \, \E[Z])\ge (1-\theta)^2 \frac{\E[Z]^2}{\E[Z^2]}.
\end{equation}
The two above inequalities (with 
$Z=|\Spec_n\cap R_n'|$ and $\theta=1/2$) imply that for all $n$,
$$
\SPb{ |\Spec_n\cap R_n'|\ge \frac {\SEb{ |\Spec_n\cap R_n'|} }2 } \ge \frac 1 {4C}.
$$

Now, by Proposition \ref{pr:same1dmarg}, one has that $\SEb{ |\Spec_n\cap R_n'|}=\E[X_n]$.
Furthermore (a trivial modification of) Proposition \ref{pr.UpperboundInf}
yields $\E[X_n] \asymp n^2 \alpha_4(n)$ which thus completes the proof.
\qed

\medskip\noindent
We are now left with 

{\medbreak \noindent {\bf Proof of Lemma \ref{l.pivotals.second.moment}.}}
As indicated at the end of the proof of 
Theorem \ref{t.some.upper.spectrum}, we have that
$\E(X_n)\asymp n^2 \alpha_4(n)$.
Next, for $x,y\in R_n'$, a picture shows that
$$
\P(x,y\in \calP_n)\le \alpha^2_4(|x-y|/2)\alpha_4(2|x-y|,n/2)
$$
since we need to have the 4-arm event around x to distance 
$|x-y|/2$, the same for $y$, and the 4-arm event in the annulus
centered at $(x+y)/2$ from distance $2|x-y|$ to distance $n/2$ and
finally these three events are independent. 
This is by quasi-multiplicity at most
$$
O(1)\alpha^2_4(n)/\alpha_4(|x-y|,n)
$$
and hence
$$
\E[|X_n|^2]\le 
O(1)\alpha^2_4(n)\sum_{x,y}\frac{1}{\alpha_4(|x-y|,n)}.
$$
Since, for a given $x$, there are at most $O(1)2^{2k}$ $y$'s with
$|x-y|\in [2^k,2^{k+1}]$, using quasi-multiplicity, the above sum is
at most
$$
O(1)n^2
\alpha^2_4(n)\sum_{k=0}^{\log_2(n)}\frac{2^{2k}}{\alpha_4(2^k,n)}.
$$
Using 
$$
\frac{1}{\alpha_4(r,R)}\le (R/r)^{2-\epsilon}
$$
(this is the fact that the four-arm exponent is strictly less than 2),
the sum becomes at most
$$
O(1)n^{4-\epsilon}\alpha^2_4(n)\sum_{k=0}^{\log_2(n)}2^{k\epsilon}.
$$
Since the last sum is at most $O(1)n^\epsilon$, we are done.
\qed

\medskip\noindent
In terms of the consequences for quantitative noise
sensitivity, Theorem \ref{t.some.upper.spectrum} implies the following 
corollary; see the exercises for similar implications. We state this only for the triangular lattice.
An analogous result holds for $\Z^2$.

\begin{corollary}\label{cor.not.full.correlation}
For $\T$, there exists $c>0$ so that if $\epsilon_n=1/(n^2\alpha_4(n))$, then
for all $n$,
$$
\P(f_n(\omega)\neq f_n(\omega_{\epsilon_n}))\ge c.
$$
\end{corollary}

\medskip\noindent
Note, importantly, this does {\it not} say that 
$f_n(\omega)$ and $f_n(\omega_{\epsilon_n})$ become asymptotically
uncorrelated, only that they are not asymptotically
completely correlated. 
To ensure that they are asymptotically uncorrelated is significantly
more difficult and requires showing that ``all'' of the spectrum is near
$n^{3/4}$. This much more difficult task is the subject of the next chapter.

\note{Comment}{It is possible that
we might want to put some general results (like the fact that  the prob. that
the spectrum is nonempty and is contained in a subset is at most the  
squared influence of that subset and that  the prob. that
the spectrum intersects a subset is at most the influence of that subset)
into Chapter \ref{ch.SM}.}

\chapter*{Exercise sheet on chapter \ref{ch.SM}}
\setcounter{exercise}{0}

\begin{exercise}\label{exercise.noiseandspectrum}
Let $\{f_n\}$ be an arbitrary sequence of Boolean functions mapping into $\{ \pm 1\}$
with corresponding spectral samples $\{ {\Spec}_n\}$. \\
(i). Show that $\SPb{0<|\Spec_n| \le A_n}\rightarrow 0$ implies that
$\SEb{(1-\epsilon_n)^{|\Spec_n|}I_{\Spec_n\neq \emptyset}}\rightarrow 0$ if $\epsilon_n A_n \rightarrow \infty$.  \\
(ii). Show that
$\SEb{(1-\epsilon_n)^{|\Spec_n|} I_{\Spec_n\neq \emptyset}}\rightarrow 0$ 
implies that
$\SPb{0<|\Spec_n| \le A_n} \rightarrow 0$ if $\epsilon_n A_n =O(1)$.
\end{exercise}

\begin{exercise}\label{exercise.stabilityandspectrum}
Let $\{f_n\}$ be an arbitrary sequence of Boolean functions mapping into
$\{ \pm 1\}$ with corresponding spectral samples $\{\Spec_n\}$. \\
(i). Show that 
$\Pb{f(\omega)\neq f(\omega_{\eps_n})}\rightarrow 0$ and 
$A_n\eps_n=\Omega(1)$ imply that
$\SPb{|\Spec_n| \ge A_n}\rightarrow 0$.
(ii). Show that
$\SPb{|\Spec_n| \ge A_n}\rightarrow 0$ and
$A_n\eps_n=o(1)$ imply that
$\Pb{f(\omega)\neq f(\omega_{\eps_n})}\rightarrow 0$.
\end{exercise}

\begin{exercise}
Prove Corollary \ref{cor.not.full.correlation}.
\end{exercise}

\begin{exercise}
For the iterated 3-Majority sequence, recall that the total influence
is $n^\alpha$ where $\alpha=1-\log 2/\log 3$. Show that for 
$\epsilon_{n}=1/n^\alpha$, 
$\P(f_n(\omega)\neq f_n(\omega_{\epsilon_n}))$ does not tend to 0.
\end{exercise}

\begin{exercise}\label{exercise4}
Assume that $\{f_n\}$ is a sequence of monotone Boolean functions on $n$
bits with total influence equal to $n^{1/2}$ up to constants. Show that the
sequence cannot be noise sensitive. Is it necessarily noise stable as the Majority function
is?
\end{exercise}

\begin{exercise}
Assume that $\{f_n\}$ is a sequence of monotone Boolean functions with mean 0 
on $n$ bits. Show that one cannot have noise sensitivity when using noise level
$\epsilon_n=1/n^{1/2}$.
\end{exercise}

\begin{exercise}
Show that $\calP$ and $\Spec$ have the same 2-dimensional marginals using
only Proposition \ref{pr.easy.spectral} rather than
Proposition \ref{pr.spectralsubdomain}. \\
Hint: It suffices to show that 
$\P(\{i,j\}\cap \Piv=\emptyset)=\SQ(\{i,j\}\cap \Spec=\emptyset)$.
\end{exercise}

\begin{exercise} (Challenging problem)
Do you expect that exercise \ref{exercise4} is sharp, meaning that, if
$1/2$ is replaced by $\alpha< 1/2$, then one can find noise sensitive examples?
\end{exercise}

\chapter[ Sharp noise sensitivity of percolation]{Sharp noise
  sensitivity of percolation }\label{ch.SNS}

\note{Timing}{///////////// 1.5 hour   ////////////// }

We will explain in this chapter the main ideas of the proof in \cite{\GPS} that most of the ``spectral mass'' 
lies near $n^2 \alpha_4(n) \approx n^{3/4+o(1)}$. This proof being rather long and involved, the content of this chapter will be far from a formal proof.
Rather it should be considered as a (hopefully convincing) heuristic 
explanation of the main results, and possibly for the interested readers as a 
``reading guide'' for the paper \cite{\GPS}. 

Very briefly speaking, the idea behind  the proof is to identify properties 
of the geometry of $\Spec_{f_n}$ which are reminiscent 
of a self-similar fractal structure. Ideally, $\Spec_{f_n}$ would behave like a spatial branching tree (or in other words a fractal percolation process), where 
distinct branches evolve independently of each other. This is conjecturally the case, but it turns out that it is very hard to control the 
dependency structure within $\Spec_{f_n}$.  
In \cite{\GPS}, only a tiny hint of spatial {\it independence} within $\Spec_{f_n}$ is proved. One of the main difficulties of the proof is to 
overcome the fact that one has very little independence to play with.

A substantial part of this chapter focuses on the much simpler case of fractal 
percolation. Indeed, this process can be seen as the simplest toy model
for the spectral sample $\Spec_{f_n}$. Explaining the simplified proof 
adapted to this setting already enables us to convey some of the main ideas
for handling $\Spec_{f_n}$.

\section{State of the art and main statement}

\begin{figure}
\begin{center}
\begin{tabular}[h]{| M{3 cm} | m{5 cm} | m{5 cm} |}

\cline{2-3}
\multicolumn{1}{ c }{} & \multicolumn{1}{|c|}{ on the square lattice $\Z^2$ } & \multicolumn{1}{c|} {on the triangular lattice $\T$} \\ \hline

The spectral mass diverges at polynomial speed &
There is a positive exponent $\eps>0$,  s.t. 
$\SPb{0< |\Spec_{f_n}| < n^\eps} \to 0$
& The same holds for all $\eps<1/8$ \\ \hline

  \multicolumn{1}{|m{ 3 cm}|}{Lower tail estimates}   &   \multicolumn{2}{ M{10 cm} | }{ On both lattices, Theorem \ref{t.ss} enables to obtain (non-sharp) lower tail estimates} \\ \hline
  
  A positive fraction of the spectral mass lies ``where it should'' & There is some universal $c>0$ s.t. 
  
  $\SPb{|\Spec_{f_n}|> c \, n^2 \alpha_4(n)} > c $ &  $\SPb{|\Spec_{f_n}|> c \, n^{3/4+o(1)}} > c $ \\ \hline
  
  \multicolumn{1}{|m{ 3 cm}|}{May be summarized by the following picture}  
  &   \multicolumn{2}{ M{10 cm} | }
  {
   \includegraphics[width=0.7\textwidth]{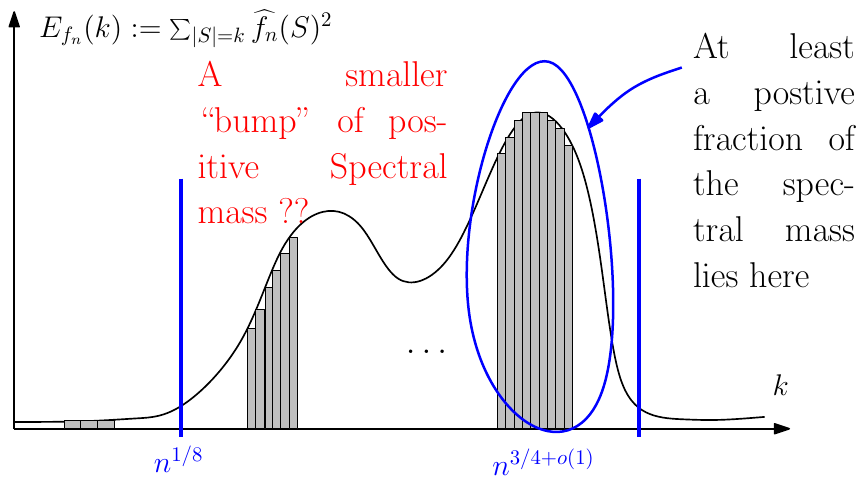}
  } \\ \hline
 
\end{tabular}
\end{center}
\caption{A summary of some of the results obtained so far for $\Spec_{f_n}$.}\label{t.tabular1}
\end{figure}

See Figure \ref{t.tabular1} where we summarize what we have learned so far 
about the spectral sample $\Spec_{f_n}$ of a left to 
right crossing event $f_n$.

From this table, we see that the main question now is to prove that all the 
spectral mass indeed diverges at speed $n^2 \alpha_4(n)$ which is $n^{3/4+o(1)}$
for the triangular lattice. This is the content of the following theorem.

\begin{theorem}[\cite{\GPS}]\label{th.GPS}

\begin{align*}
\limsup_{n \to \infty}  \SPb{ 0 < |\Spec_{f_n}| < \lambda\, n^2 \alpha_4(n)}  \underset{\lambda \to 0}{\longrightarrow} 0\, \,.
\end{align*}
\end{theorem}

On the triangular lattice $\T$, the rate of decay in $\lambda$ is known explicitly. Namely:

\begin{theorem}[\cite{\GPS}]\label{th.GPS2}
On the triangular grid $\T$, the lower tail of $|\Spec_{f_n}|$ satisfies
\begin{align*}
\limsup_{n \to \infty}  \SPb{ 0 < |\Spec_{f_n}| < \lambda\,  \SEb{|\Spec_{f_n}| })}  \underset{\lambda \to 0}{\asymp}  \lambda^{2/3}\, \,.
\end{align*}
\end{theorem}

This result deals with what one might call the ``macroscopic'' lower tail, 
i.e. with quantities which asymptotically are still of 
order $\SEb{|\Spec_{f_n}|}$
(since $\lambda$ remains fixed in the above $\limsup$).
It turns out that in our later study of dynamical percolation in
Chapter \ref{ch.DP}, we will need a sharp control on the full lower tail.
This is the content of the following stronger theorem:

\begin{theorem}[\cite{\GPS}]\label{th.GPS3}
On $\Z^2$ and on the triangular grid $\T$, for all $1\le r \le n$, one has 
\begin{align*}
\SPb{ 0 < |\Spec_{f_n}| < r^2 \alpha_4(r) }  \asymp  \frac{n^2}{r^2} \alpha_4(r,n)^2  \,.
\end{align*}

\ni
On the triangular grid, this translates into 
\begin{align*}
\SPb{ 0 < |\Spec_{f_n}| <  u  }  \approx  n^{- \frac 1 2} u^{\frac 2 3}  \,,
\end{align*}
where we write $\approx$ to avoid relying on $o(1)$ terms in the exponents.  

\end{theorem}

\section{Overall strategy}

In the above theorems, it is clear that we are mostly interested in the 
cardinality of $\Spec_{f_n}$. However, our strategy will consist in 
understanding as much as we can about the typical {\it geometry} of 
the random set $\Spec_{f_n}$ sampled according to the spectral probability measure 
$\hat \P_{f_n}$.

\begin{figure}[!h]
\centering
\includegraphics[width=1\textwidth]{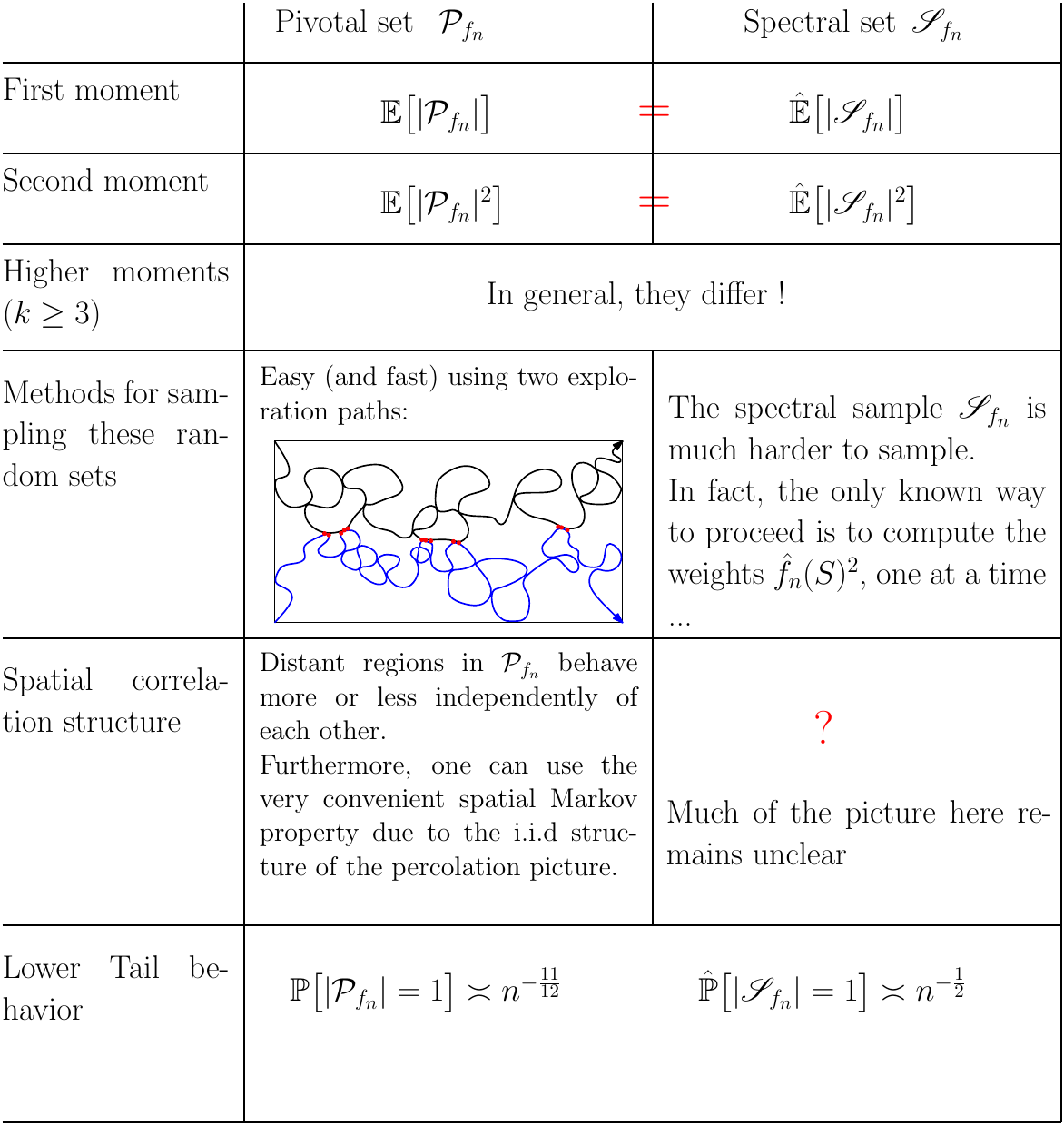}
\caption{Similarities and differences between $\Spec_{f_n}$ and $\Piv_{f_n}$.}\label{f.tableau}
\end{figure}

As we have seen so far, the random set $\Spec_{f_n}$ shares many properties with the set of pivotal points $\Piv_{f_n}$. A first possibility would be that 
they are asymptotically similar. After all, noise sensitivity is intimately 
related with pivotal points, so it is not unreasonable to hope for such a 
behavior. This scenario would be very convenient for us since the geometry of 
$\Piv_{f_n}$ is now well understood (at least on $\T$) thanks to the 
$\mathrm{SLE}$ processes. In particular, in the case of $\Piv_{f_n}$, one 
can ``explore'' $\Piv_{f_n}$ in a Markovian way by relying on 
exploration processes. Unfortunately, based on very convincing heuristics, it 
is conjectured that the scaling limits of $\frac 1 n \Spec_{f_n}$ and 
$\frac 1 n \Piv_{f_n}$ are singular random compact sets of the square. 
See Figure \ref{f.tableau} for a quick overview of the similarities and 
differences between these two random sets.

The conclusion of this table is that they indeed share many properties, but 
one cannot deduce lower tail estimates on $|\Spec_{f_n}|$ out of 
lower tail estimates on $|\Piv_{f_n}|$. Also, even worse, we will not be 
allowed to rely on spatial Markov properties for $\Spec_{f_n}$.

However, even though $\Piv_{f_n}$ and $\Spec_{f_n}$ differ in many ways, 
they share at least one essential property: a seemingly 
{\it self-similar fractal behavior}. The main strategy in \cite{\GPS} to 
control the lower-tail behavior of $|\Spec_{f_n}|$ is to prove that 
\underline{in some very weak sense}, $\Spec_{f_n}$ behaves like the 
simplest model among self-similar fractal processes in $[0,n]^2$: i.e. a super-critical spatial Galton-Watson tree 
embedded in $[0,n]^2$, also called a {\it fractal percolation process}. The lower tail of this very simple toy model will be investigated in detail in the next section
with a technique which will be suitable for $\Spec_{f_n}$. The main difficulty which arises in this program is the lack of knowledge of the independency
structure within $\Spec_{f_n}$.
In other words, when we try to compare $\Spec_{f_n}$ with a fractal percolation process, the self-similarity already requires some work, but the hardest part  is to deal with the fact that distinct ``branches'' (or rather their analogues) 
are not known to behave even slightly independently of each other. 
We will discuss these issues in Section \ref{s.tour} but will not give a complete proof.

\section{Toy model: the case of fractal percolation}\label{s.toy}

As we explained above, our main strategy is to exploit the fact that $\Spec_{f_n}$ has a certain self-similar fractal structure.
Along this section, we will consider the simplest case of such a self-similar fractal object: namely {\it fractal percolation}, and we will detail in this simple setting what our later strategy will be.
Deliberately, this strategy will not be optimal in this simplified case. In particular, we will not rely on the martingale techniques that one can use with 
fractal percolation or Galton-Watson trees,
since such methods would not be available for our spectral sample $\Spec_{f_n}$.

\subsection{Definition of the model and first properties}

To make the analogy with $\Spec_{f_n}$ easier let 
\[
n:=2^h\,, h\geq 1\, ,
\]
and let's fix a parameter $p\in(0,1)$.

\begin{figure}[!h]
\begin{center}
\includegraphics[width=0.6 \textwidth]{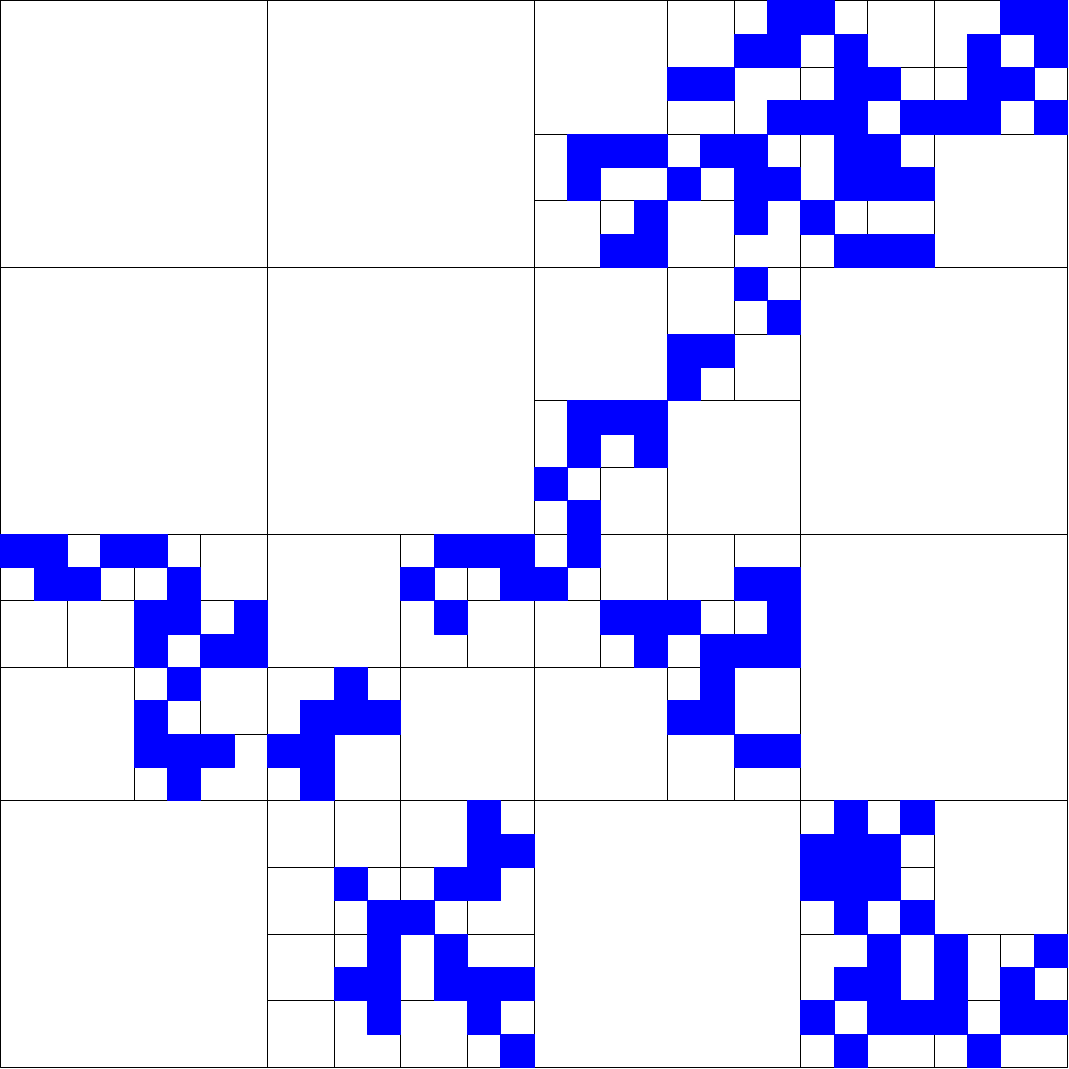}
\end{center}\caption{A realization of a fractal percolation $\Tree_{2^5}=\Tree^{5}$}\label{f.Tree5}
\end{figure}

Now, {\it fractal percolation} on $[0,n]^2$ is defined inductively as follows: 
divide $[0,2^h]^2$ into 4 squares and retain each of them 
independently with probability $p$. Let $\Tree^{1}$ be the union of the 
retained $2^{h-1}$-squares. The second-level tree $\Tree^{2}$ is obtained by reiterating 
the same procedure independently for each $2^{h-1}$-square in $\Tree^{1}$. Continuing in the same fashion all the way to the squares of unit size, one obtains 
$\Tree_n= \Tree:=\Tree^{h}$ which is a random subset of $[0,n]^2$. See \cite{\LyonsPeres} for more on the definition of {\it fractal percolation}. See also 
Figure \ref{f.Tree5} for an example of $\Tree^{5}$. 

\begin{remark}\label{r.2not}
We thus introduced two different notations for the same random set ($\Tree_{n=2^h} \equiv \Tree^h$). The reason for this is that on the one hand 
the notation $\Tree_n$ defined on $[0,n]^2=[0,2^h]^2$ makes the analogy with $\Spec_{f_n}$ (also defined on $[0,n]^2$) easier, while on the other hand inductive proofs 
will be more convenient with the notation $\Tree^h$.
\end{remark}

In order to have a supercritical Galton-Watson tree, one has to choose $p\in (1/4,1)$. Furthermore, one can easily check the following easy proposition.
\begin{proposition}\label{pr.tree}
Let $p\in (1/4, 1)$. Then

\[
\Eb{|\Tree_n|} = n^2 p^h = n^{2+ \log_2 p}\,,
\]
and 
\[
\Eb{|\Tree_n|^2} \le O(1) \Eb{|\Tree_n|}^2\,.
\]

In particular, by the second moment method
(e.g.\ the Paley-Zygmund inequality),
with positive probability, $\Tree_n$ is of order $n^{2+\log_2 p}$.
\end{proposition}

Let 
\[
\alpha := 2 + \log_2 p\,.
\]

This parameter $\alpha$ corresponds to the ``fractal dimension'' of $\Tree_n$.
To make the analogy with $\Spec_{f_n}$ even clearer, one could choose $p$ in such a way that $\alpha = 2+\log_2 p = 3/4$, but we will not need to.

The above proposition implies that on the event $\Tree_n \neq \emptyset$, with positive conditional probability $|\Tree_n|$ is large (of order $n^\alpha$).
This is the exact analogue of Theorem \ref{t.some.upper.spectrum} for the spectral sample $\Spec_{f_n}$.

Let us first analyze what would be the analogue of Theorem \ref{th.GPS} in the case of our toy model $\Tree_n$. We have the following.
\begin{proposition}\label{pr.01tree}
\begin{align*}
\limsup_{n \to \infty}  \Pb{ 0 < |\Tree_{n}| < \lambda\, n^\alpha)}  \underset{\lambda \to 0}{\longrightarrow} 0\, \,.
\end{align*}
\end{proposition}

\begin{remark}
If one could rely on martingale techniques, then this proposition is a 
corollary of standard results. Indeed, as is well-known
\[
M_i:=  \frac {|\Tree^{i}|} {(4p)^i}\,,
\]
is a positive martingale. Therefore it converges, as $n\to \infty$, 
to a non-negative random variable $W\geq 0$. Furthermore, the conditions of the 
Kesten-Stigum Theorem are fulfilled  (see for example Section 5.1 in 
\cite{\LyonsPeres}) and therefore $W$ is positive on the event 
that there is no extinction. This implies the above proposition.
\end{remark}


As we claimed above, we will intentionally follow a more hands-on approach in this section which will be more suitable to the random set $\Spec_{f_n}$
which we have in mind. Furthermore this approach will have the great 
advantage to provide the following much more precise result, 
which is the analogue of Theorem \ref{th.GPS3} for $\Tree_n$.

\begin{proposition}\label{pt.Tree3}
For any $1\le r \le n$,
\[
\Pb{0 < |\Tree_n| < r^\alpha  } \asymp (\frac r n)^{\log_2 1/\mu }\,,
\]
where $\mu$ is an explicit constant in $(0,1)$ computed in Exercise \ref{ex.induction}.
\end{proposition}

\subsection{Strategy and heuristics}

Letting $u \ll n^\alpha$, we wish to estimate $\Pb{0< |\Tree_n| < u}$. 
Even though we are only interested in the size of $\Tree_n$, we will try to estimate this quantity by understanding the 
{\it geometry} of the conditional set:
\[
\Tree_n^{|u}:= \mathcal{L} \Bigl( \Tree_n \Bigm| 0< |\Tree_n| < u \Bigr)\,.
\]

\begin{figure}
\begin{center}
\begin{tabular}[h]{| M{3 cm} | m{5 cm} | m{1 cm} | m{5 cm} |}

\hline
\multicolumn{1}{|c|}{} & \multicolumn{3}{c|}{} \\ 
\multicolumn{1}{|M{3 cm}|}{ How does  $\mathcal{L} \Bigl( \Tree_n \Bigm| 0< |\Tree_n| < u \Bigr)$ look ? } &
\multicolumn{1}{ M{5 cm}} { \includegraphics[width=0.33 \textwidth]{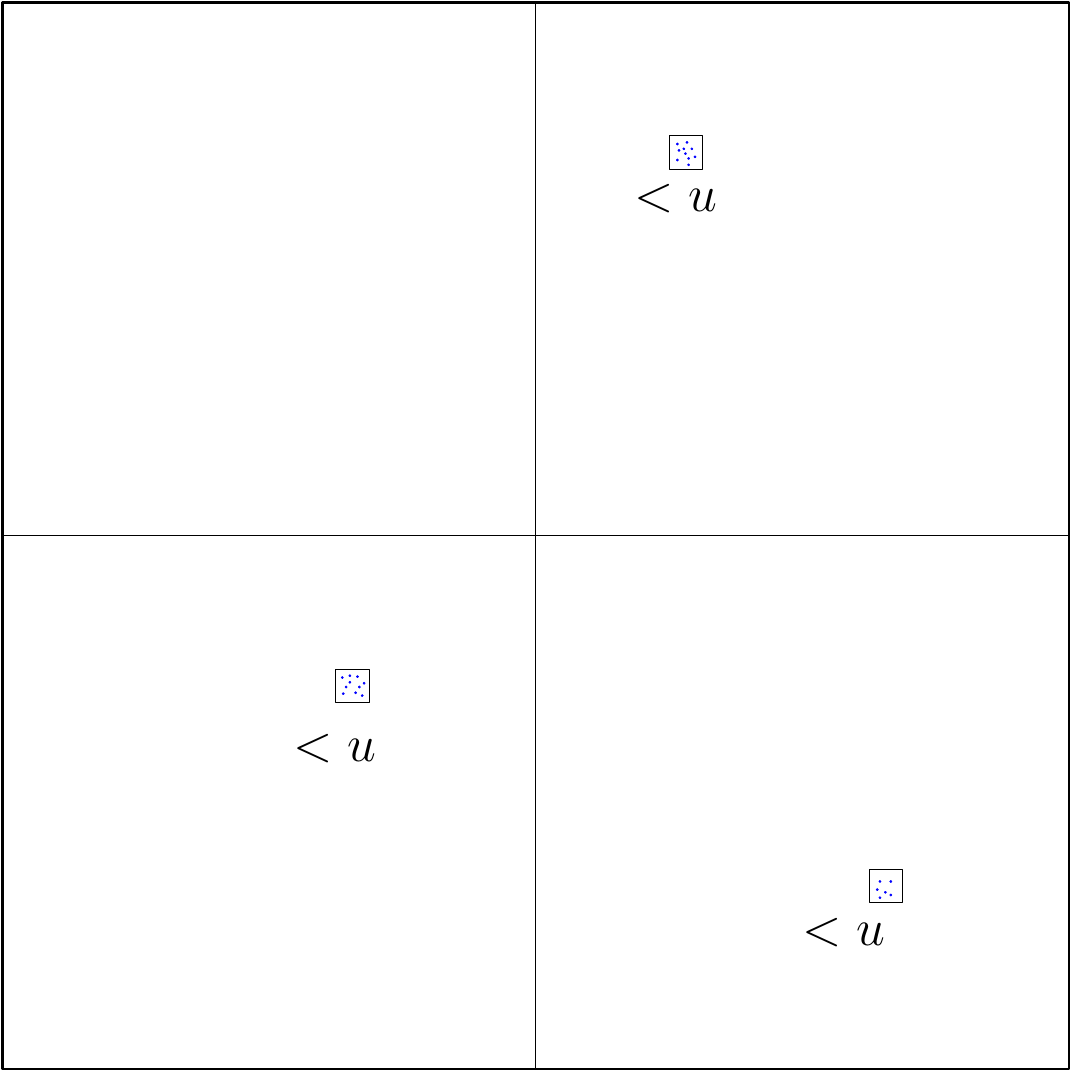}
More Entropy ( in $Vol^3$) but costs more to maintain these 3 ``islands'' alive.}
& 
\multicolumn{1} {c} {OR?} &
\multicolumn{1}{M{5 cm} |} { \includegraphics[width=0.33 \textwidth]{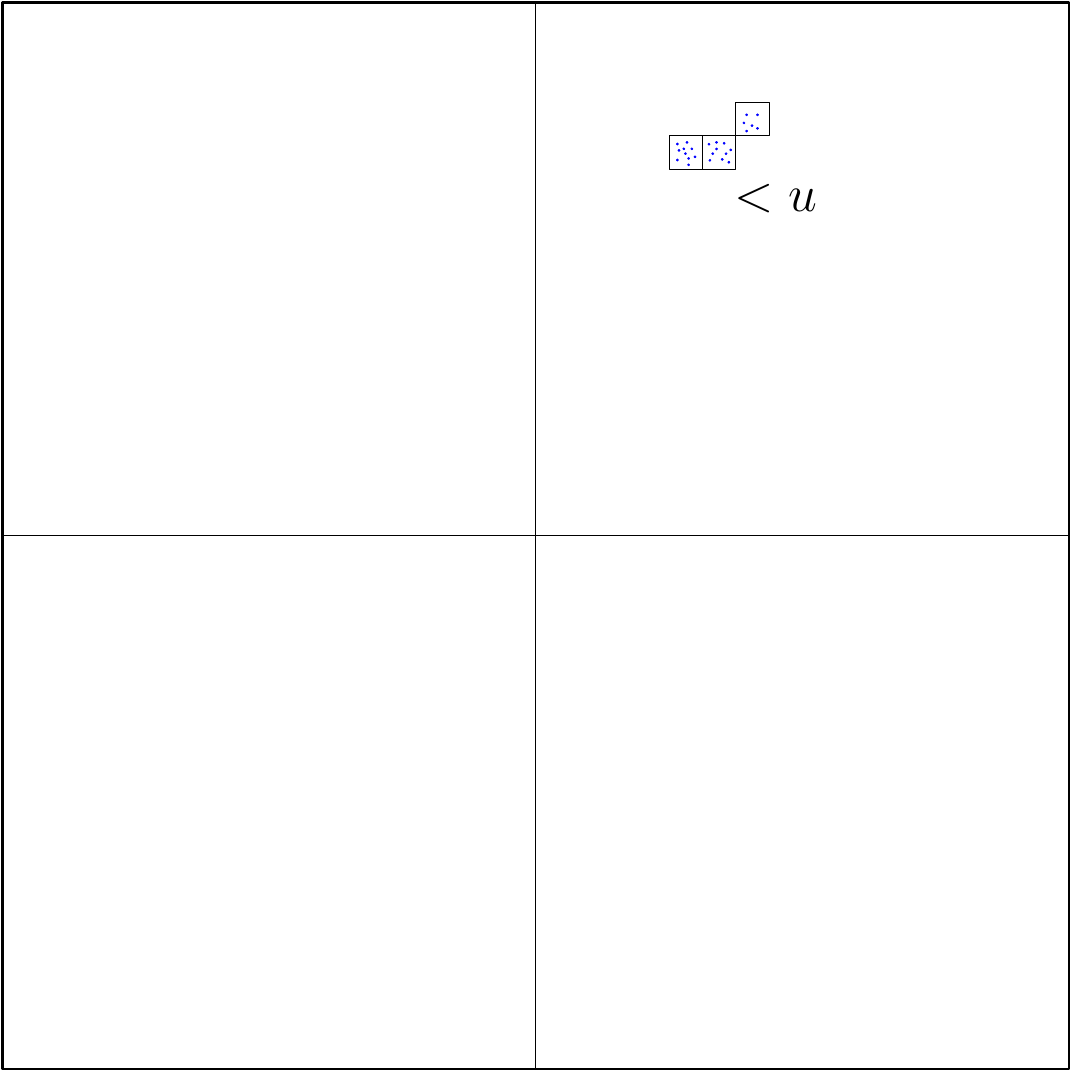}
Less Entropy ( in $Vol^1$) but only one island to maintain alive.}

 \\ \hline

\end{tabular}
\end{center}
\caption{Entropy v.s. Clustering effect}\label{t.entropy}
\end{figure}

The first natural question to ask is whether this conditional random set is 
typically {\it localized} or not. See Figure \ref{t.entropy}.

Intuitively, it is quite clear that the set $\Tree_n$ conditioned to be very 
small will tend to be localized. So it is the picture on the right in 
Figure \ref{t.entropy} which is  more likely. This would deserve a proof 
of course, but we will come back to this later.
The fact that it should look more and more localized tells us that
as one shrinks $u$, this should make our 
conditional $\Tree_n^{|u}$ more and more singular with respect to the 
unconditional one. But how much localization should we see?  
This is again fairly easy to answer, at least on the intuitive level. 
Indeed, $\Tree_n^{|u}$ should tend to localize until it reaches a certain 
mesoscopic scale $r$ such that 
$1\ll r \ll n$. One can compute how much it costs to maintain a single branch (or $O(1)$ branches) alive until scale $r$, but once this is achieved, one should let the system 
evolve in a ``natural'' way. In particular, once the tree survives all the way to a mesoscopic square of size $r$, it will (by the second moment method) produce $\Omega(r^\alpha)$ leaves 
there with positive probability.

To summarize, typically $\Tree_n^{|u}$ will maintain $O(1)$ many branches alive at scale $1\ll r \ll n$, 
and then it will let the branching structure evolve in a basically unconditional way. 
The intermediate scale $r$ is chosen so that $r^\alpha \asymp u$.

\begin{definition}\label{d.TreeR}
If $1\le r\le n=2^h$ is such that $r=2^l, 0\le l \le h$, let $\Tree_{(r)}$ denote the set of branches that were still alive at scale $r=2^l$ in the 
iterative construction of $\Tree_n$. In other words, $\Tree_{(r)}\equiv \Tree^{h-l}$ and $\Tree_n \subset \bigcup \Tree_{(r)}$. 
This random set $\Tree_{(r)}$ will be the analogue of the ``$r$-smoothing'' $\Spec_{(r)}$ of the spectral sample $\Spec_{f_n}$ defined 
later in Definition \ref{d.meso}.
\end{definition}

Returning to our problem, the above heuristics say that one expects to have for any $1\ll u \ll n^\alpha$.
\begin{align*}
\Pb{0< |\Tree_n| < u}  &\asymp \Pb{0< |\Tree_{(r)}| \le O(1)} \\
& \asymp \Pb{ |\Tree_{(r)}| = 1}\,,
\end{align*}
where $r$ is a dyadic integer chosen such that $r^\alpha \asymp u$. 
Or in other words, we expect that 

\begin{align}\label{e.Tr1}
\Pb{0< |\Tree_n| < r^\alpha} \asymp \Pb{|\Tree_{(r)}| = 1}\,.
\end{align}

In the next subsection, we briefly explain how this heuristic can be 
implemented into a proof in the case of the tree $\Tree_n$
in a way which will be suitable to the study of $\Spec_{f_n}$. We will only 
skim through the main ideas for this tree case.

\subsection{Setup of a proof for $\Tree_n$}

Motivated by the above heuristics, we divide our system into two scales: above and below the mesoscopic scale $r$.
One can write the lower tail event as follows (let $1 \ll r\ll n$):

\begin{align}\label{e.Midentity}
\Pb{0< |\Tree_n|<r^\alpha} 
& = \sum_{k\geq 1}\,  \Pb{ |\Tree_{(r)}| = k} \,  \Pb{ 0 < |\Tree_n| < r^\alpha \md |\Tree_{(r)}| = k} \,.
\end{align}

It is not hard to estimate the second term $\Pb{ 0 < |\Tree_n| < r^\alpha \md |\Tree_{(r)}| = k}$. 
Indeed, in this term we are conditioning on having exactly $k$ branches alive at scale $r$. Independently of where they are, ``below'' $r$, these $k$ branches 
evolve independently of each other. Furthermore, by the second moment method, there is a universal constant $c>0$ such that each of them exceeds the fatal amount of $r^\alpha$ leaves with probability at least $c$ (note that in the opposite direction, each branch could also go extinct with positive probability).
This implies that 
\[
\Pb{ 0 < |\Tree_n| < r^\alpha \md |\Tree_{(r)}| = k} \le (1-c)^k\,.
\]

\begin{remark}
Note that one makes heavy use of the independence structure within $\Tree_n$ 
here. This aspect is much more nontrivial for the spectral sample 
$\Spec_{f_n}$. Fortunately it turns out, and this is a key fact, that in 
\cite{\GPS} one can prove a weak independence statement which in some sense makes it possible to follow this route. 
\end{remark}

We are left with the following upper bound:
\begin{align}\label{e.factorTree}
\Pb{0< |\Tree_n|<r^\alpha} \le \sum_{k\geq 1}  \Pb{ |\Tree_{(r)}| = k} (1-c)^k\,.
\end{align}

In order to prove our goal of \eqref{e.Tr1}, by exploiting the exponential 
decay given by $(1-c)^k$ (which followed from independence), it is enough to 
prove the following bound on the mesoscopic behavior of $\Tree$:

\begin{lemma}\label{l.treeGk}
There is a {\bf sub-exponential} function $k\mapsto g(k)$ such that 
for all $1 \le r \le n$,
\[
\Pb{|\Tree_{(r)}| = k} \le g(k) \, \Pb{|\Tree_{(r)}| = 1}\,.
\]
\end{lemma}

Notice as we did in Definition \ref{d.TreeR} that since $\Tree_{(r)}$ has the same law as $\Tree^{h-l}$, this is a purely Galton-Watson tree type of question.

The big advantage of our strategy so far is that initially we were looking for 
a sharp control on $\Pb{0< |\Tree_n| < u}$ and now, using this ``two-scales'' 
argument, it only remains to prove a crude upper bound on the lower tail of 
$|\Tree_{(r)}|$. By scale invariance this is nothing else than obtaining a 
crude upper bound on the lower tail of $|\Tree_n|$. Hence this division into 
two scales greatly simplified our task.

\subsection{Sub-exponential estimate on the lower-tail (Lemma \ref{l.treeGk})}\label{ss.Treesubexpo}

The first step towards proving and understanding Lemma \ref{l.treeGk} is to understand the term $\Pb{|\Tree_{(r)}|=1}$. 
From now on, it will be easier to work with the ``dyadic'' notations instead, i.e. with $\Tree^i\equiv \Tree_{2^i}$ (see remark \ref{r.2not}).
With these notations, the first step is equivalent to understanding the probabilities $p_i:= \Pb{|\Tree^i|=1}$. This aspect of the problem is very specific to 
the case of Galton-Watson trees and gives very little insight into the 
later study of the spectrum $\Spec_{f_n}$.
Therefore we postpone the details to Exercise \ref{ex.induction}. 
The conclusion of this (straightforward) exercise is that $p_i$ behaves as $i\to \infty$ like

\[
p_i \sim c\,  \mu^i \,,
\]
for an explicit exponent $\mu \in (0,1)$ (see Exercise \ref{ex.induction}). 
In particular, in order to prove Proposition \ref{pt.Tree3},
it is now enough to find a sub-exponential function $k\mapsto g(k)$ such that for any $i,k \geq 1$,
\begin{align}\label{e.muk}
\Pb{|\Tree^i| = k } \le g(k) \mu^i\,.
\end{align}

More precisely, we will prove the following lemma.
\begin{lemma}\label{l.treeGk2}
Let $g(k):= 2^{\theta \log_2^2(k+2)}$, where $\theta$ is a fixed constant to be chosen later.
Then for all $i,k \geq 1$, one has 
\begin{align}\label{e.treeGk2}
\Pb{|\Tree^i|=k} \le g(k)\, \mu^i\,.
\end{align}
\end{lemma}

We provide the proof of this lemma here, since it can be seen as a ``toy proof'' of the corresponding sub-exponential estimate needed for the
$r$-smoothed spectral samples $\Spec_{(r)}$, stated in the coming Theorem \ref{th.GPSgk}. The proof of this latter theorem shares some similarities with the proof below
but is much more technical 
since in the case of $\Spec_{(r)}$ one has to deal with a more complex structure than the branching structure of a Galton-Watson tree.

\proof
We proceed by double induction.
Let $k\geq 2$ be fixed and assume that equation~\eqref{e.treeGk2} is already satisfied for all pair $(i',k')$ such that $k'<k$.
Based on this assumption, let us prove by induction on $i$ that all pairs $(i,k)$ satisfy equation~\eqref{e.treeGk2} as well. 

First of all, if $i$ is small enough, this is obvious by the definition of $g(k)$.
Let 
\[
J = J_k := \sup \{ i \geq 1: g(k) \mu^i >  10 \}\,.
\]

Then, it is clear that equation~\eqref{e.treeGk2} is satisfied for all $(i,k)$ with $i\le J_k$. 
Now let $i>J_k$.

If $\Tree^i$ is such that $|\Tree^i|=k\ge 1$, let $L=L(\Tree^i)\geq 0$ be the largest integer such that 
$\Tree^{i}$ intersects only one square of size $2^{i-L}$. This means that below scale $2^{i-L}$, the tree $\Tree^i$
splits into at least 2 live branches in distinct dyadic squares of size $2^{i-L-1}$. Let $d\in\{2,3,4\}$ be the number of such live branches.
By decomposing on the value of $L$, and using the above assumption, we get

$$
\Pb{|\Tree^i| = k}  \le  \Pb{ L(\Tree^i) > i-J_k} +
$$
$$
\frac{1}{1-q}\sum_{l=0}^{i-J_k} \Pb{L(\Tree^i)=l} \,  \sum_{d=2}^4  \binom{4}{d} \,  
(\mu^{i-l-1})^d 
\sum_{ \begin{array}{l} (k_j)_{1\le j \le d} \\k_j\ge 1,  \sum k_j = k \end{array} }
\prod_j g(k_j)
$$
where $q$ is the probability that our Galton-Watson tree goes extinct.

Let us first estimate what $\Pb{L(\Tree^i)\geq m}$ is for $m\geq 0$. 
If $m\geq 1$, this means that among the $2^{2m}$ dyadic squares of size $2^{i-m}$, only 
one will remain alive all the way to scale 1. Yet, it might be that some other such squares are still alive at scale $2^{i-m}$ but will go extinct by the time 
they reach scale 1. Let $p_{m,b}$ be the probability that the process $\Tree^{m+b}$\!\!, 
which lives in $[0,2^{m+b}]^2$, is entirely contained in a dyadic square of size $2^b$. 
With such notations, one has 
\[
\Pb{L(\Tree^i)\geq m}  = p_{m, i-m}\,. 
\]

Furthermore, if $i=m$, one has $p_{i,0}=p_i\sim c \mu^i$.
It is not hard to prove (see Exercise \ref{ex.induction}) the following lemma.
\begin{lemma}\label{l.pmh}
For any value of $m,b\geq 0$, one has
\[
p_{m,b}\le  \mu^m\,.
\]
In particular, one has a universal upper bound in $b\geq 0$.
\end{lemma}

It follows from the lemma that $\Pb{L(\Tree^i) = l} \le \Pb{ L(\Tree^i)\geq l} \le  \mu^l$  and 
\begin{align}
\Pb{ L(\Tree^i)> i-J_k} & \le  \mu^{i-J_k} \\ 
& \le  \, \frac 1 {10} \, g(k)\, \mu^i \,  \text{ by the definition of $J_k$}\,. 
\end{align}

This gives us that for some constant $C$
\begin{align*}
\Pb{|\Tree^i| = k}  &\le  \frac {\mu^{i}} {10} g(k)  + C\sum_{l=0}^{i-J_k} \mu^l  \sum_{d=2}^4 \, (\mu^{i-l})^d  
\sum_{ \begin{array}{l} (k_j)_{1\le j \le d} \\k_j\ge 1,  \sum k_j = k \end{array} }
\prod_j g(k_j)\\
&=   \frac {\mu^{i}} {10} g(k) + C\mu^i  \sum_{d=2}^4 \, \sum_{l=0}^{i-J_k}  (\mu^{i-l})^{d-1}  
\sum_{ \begin{array}{l} (k_j)_{1\le j \le d} \\k_j\ge 1,  \sum k_j = k \end{array} }
\prod_j g(k_j) \,. 
\end{align*}

Let us deal with the $d=2$ sum (the contributions coming from $d>2$ being even smaller). 
By concavity of $k\mapsto \theta \log_2^2(k+2)$, one obtains that for any $(k_1,k_2)$ such that $k_1+k_2=k$:
$g(k_1)g(k_2)\le g(k/2)^2$. Since there are at most $k^2$ such pairs, 
this gives us the following bound on the $d=2$ sum.

\begin{align*}
\sum_{l=0}^{i-J_k}  (\mu^{i-l})^{2-1}  
\sum_{ \begin{array}{l} (k_j)_{1\le j \le 2} \\k_j\ge 1,  \sum k_j = k \end{array} }
\prod_j g(k_j)
& \le  \sum_{l=0}^{i-J_k}  \mu^{i-l}  k^2 g(k/2)^2 \\
&\le \frac{1}{1-\mu}\, \mu^{J_k} \,k^2 \, g(k/2)^2 \\ 
&\le 10 \frac{1}{1-\mu} \, k^2 \, g(k/2)^2\,  (\mu g(k))^{-1}\,,
\end{align*}
by definition of $J_k$.

Now, some easy analysis implies that if one chooses the constant $\theta>0$ large enough, then 
for any $k\geq 2$, one has $C10 \frac{1}{1-\mu} k^2 \,  g(k/2)^2 \, (\mu g(k))^{-1} 
\le \frac 1 {10} g(k)$. 
Altogether (and taking into consideration the $d>2$ contributions), this implies that 

\[
\Pb{|\Tree^i| = k} \le \frac 2 5 g(k) \mu^i \le g(k) \mu^i \,,
\]
as desired.
\qed

\begin{remark}\label{r.entropysolved}
Recall the initial question from Figure \ref{t.entropy} which asked whether the clustering effect wins 
over the entropy effect or not. This question enabled us to motivate the setup of the proof
but in the end, we did not specifically address it. Notice that the above proof in fact solves the problem (see Exercise \ref{ex.entropy}).
\end{remark}


\section{Back to the spectrum: an exposition of the proof}\label{s.tour}

\subsection{Heuristic explanation}

Let us now apply the strategy we developed for $\Tree_n$ to the case of the 
spectral sample $\Spec_{f_n}$.
Our goal is to prove Theorem \ref{th.GPS3} (of which Theorems \ref{th.GPS} and \ref{th.GPS2} are straightforward corollaries). 
Let $\Spec_{f_n}\subset [0,n]^2$ be our spectral sample. We have seen (Theorem \ref{t.some.upper.spectrum}) that with positive 
probability $|\Spec_{f_n}|\asymp n^2 \alpha_4(n)$. For all $1< u < n^2 \alpha_4(n)$, we wish to understand the probability
$\SPb{0< |\Spec_{f_n}| <u}$. Following the notations we used for $\Tree_n$, let $\Spec_{f_n}^{|u}$ be the spectral sample conditioned on the event $\{ 0 <|\Spec_{f_n}|<u \}$.
\vskip 0.1 cm
\ni
{\bf Question:} How does $\Spec_{f_n}^{|u}$ typically look?

To answer this question, one has to understand whether $\Spec_{f_n}^{|u}$ tends to be {\it localized} or not. Recall from Figure \ref{t.entropy} the illustration
of the competition between entropy and clustering effects in the case of $\Tree_n$. The same 
figure applies to the spectral sample $\Spec_{f_n}$. 
We will later state 
a {\bf clustering lemma} (Lemma \ref{l.cluster}) which will strongly support 
the localized behavior described in the next proposition.

Therefore we are guessing that our conditional set $\Spec_{f_n}^{|u}$ will tend to localize into 
$O(1)$ many  squares of a certain scale
$r$ and will have a ``normal'' size within these $r$-squares. It remains to 
understand what this mesoscopic scale $r$ as a function of $u$ is.

By ``scale invariance'', one expects that if $\Spec_{f_n}$ is conditioned to live in a square of size $r$, then $|\Spec_{f_n}|$ will be of order $r^2 \alpha_4(r)$
with positive conditional probability. More precisely, the following lemma will be proved in 
Problem \ref{ex.scaleinv}.

\begin{lemma}\label{l.prob}
There is a universal $c\in(0,1)$ such that for any $n$ and
for any $r$-square $B\subset [n/4, 3n/4]^2$ in the ``bulk'' of $[0,n]^2$, one has
\begin{align}
\SPb{\frac{|\Spec_{f_n}|}{r^2 \alpha_4(r)} \in (c, 1/c)  \md \Spec_{f_n} \neq \emptyset \text{ and } \Spec_{f_n} \subset B}  > c\,.
\end{align}
\end{lemma}
In fact this lemma holds uniformly in the position of the $r$-square $B$ inside $[0,n]^2$, but we will not discuss this here.
\vskip 0.3 cm

What this lemma tells us is that for any $1<u< n^2 \alpha_4(n)$, if one  
chooses $r=r_u$ in such a way that  $ r^2 \alpha_4(r) \asymp u$, 
then we expect to have the following estimate:
\begin{align*}
\SPb{0<|\Spec_{f_n}|<u}  & \asymp \SPb{\Spec_{f_n} \text{ intersects $O(1)$ $r$-squares in }[0,n]^2 } \\
& \asymp \SPb{\Spec_{f_n} \text{ intersects a single $r$-square in }[0,n]^2 }
\end{align*}

At this point, let us introduce a concept which will be very helpful in what 
follows.
\begin{definition}[``$r$-smoothing'']\label{d.meso}
Let $1\le r \le n$.
Consider the domain $[0,n]^2$ and divide it into a grid of squares of edge-length $r$.
(If $1\ll r \ll n$, one can view this grid as a {\it mesoscopic} grid).

If $n$ is not divisible by $r$, write $n= m r + q$ and consider the grid of $r$-squares covering $[0,(m+1)r]^2$.

Now, for each subset $S\subset [0,n]^2$, define $S_{(r)}$ to be the set of  
$r\times r$ squares in the above grid
which intersect $S$. In particular $|S_{(r)}|$ will correspond to the number of such $r$-squares  which intersect $S$. With a slight abuse of notation, 
$S_{(r)}$ will sometimes also denote the actual subset of $[0,n]^2$ 
consisting of the union of these $r$-squares.

One can view the application $S \mapsto S_{(r)}$ as an {\bf $r$-smoothing} since all the details below the scale $r$ are lost.
\end{definition}

\begin{remark}
Note that in Definition \ref{d.TreeR}, we relied on a slightly different notion of ``$r$-smoothing'' since in that case, $\Tree_{(r)}$ could also include 
$r$-branches which might go extinct by the time they reached scale one. The advantage of this choice was that there was an exact scale-invariance 
from $\Tree$ to $\Tree_{(r)}$ while in the case of $\Spec_{f_n}$, there is no such exact scale-invariance from $\Spec$ to $\Spec_{(r)}$.
\end{remark}

With these notations, the above discussion leads us to believe that the following proposition should hold.
\begin{proposition}\label{pr.onebox}
For all $1\le r \le n$, one has 
\[
\SPb{ 0 <|\Spec_{f_n}| <r^2 \alpha_4(r) } \asymp \SPf{f_n}{|\Spec_{(r)}|=1}\,.
\]
\end{proposition}

Before explaining the setup used in \cite{\GPS} to prove such a result, let us check that it indeed implies Theorem \ref{th.GPS3}. 
By neglecting the boundary issues, one has 
\begin{align}
\SPf{f_n}{|\Spec_{(r)}| = 1} \asymp \sum_{\begin{array}{c} r\text{-squares } \\ B \subset [n/4, 3n/4]^2 \end{array}} \SPb{\Spec_{f_n} \neq \emptyset \text{ and } \Spec_{f_n} \subset B }\,.
\end{align}
There are $O(\frac {n^2}{r^2})$ such $B$ squares, and for each of these, one can check (see Exercise \ref{ex.onebox}) that 
\begin{equation}\label{e.4armsquared}
\SPb{\Spec_{f_n} \neq \emptyset \text{ and } \Spec_{f_n} \subset B } \asymp \alpha_4(r,n)^2\,.
\end{equation}
Therefore, Proposition \ref{pr.onebox} indeed implies Theorem \ref{th.GPS3}.

\subsection{Setup and organization of the proof of Proposition \ref{pr.onebox}}


To start with, assume we knew that disjoint regions in the spectral sample $\Spec_{f_n}$ behave more or less independently of each other
in the following (vague) sense. For any $k\geq 1$ and any mesoscopic scale $1\le r \le n$, if one conditions on $\Spec_{(r)}$ to be equal 
to $B_1 \cup \cdots \cup B_k$ for $k$ disjoint $r$-squares, then the conditional law of $\Spec_{\md \bigcup B_i}$  should be ``similar'' to an independent product 
of $\mathcal{L} \bigl[\Spec_{\md B_i} \bigm| \Spec\cap B_i \neq \emptyset \bigr],\, i\in \{1,\ldots, k\}$. Similarly as in the tree case (where the analogous property for $\Tree_n$ was an exact independence factorization), and assuming that the above comparison with an independent product could be made quantitative, 
this would potentially imply the following upper bound for a  certain absolute constant $c>0$:
\begin{align}\label{e.ifindependent}
\SPb{0 < |\Spec_{f_n}| < r^2 \alpha_4(r)} \le \, \sum_{k\geq 1} \,  \SPb{|\Spec_{(r)}| = k} \, (1-c)^k\,.
\end{align}

This means that even if one managed to obtain a good control on the dependency structure within $\Spec_{f_n}$ (in the above sense), one would still need to have a good 
estimate on $\SPb{|\Spec_{(r)}| = k}$ in order to deduce Proposition \ref{pr.onebox}. This part of the program is achieved in \cite{\GPS} without requiring any information 
on the dependency structure of $\Spec_{f_n}$. More precisely, the following result is proved:

\begin{theorem}[\cite{\GPS}]\label{th.subexpoS}
There is a sub-exponential function $g\mapsto g(k)$, such that for any $1\le r \le n$ and 
any $k\geq 1$,
\begin{align*}
\SPb{|\Spec_{(r)}| = k }\, \le \, g(k) \, \SPb{|\Spec_{(r)}| = 1}\,.
\end{align*}
\end{theorem}

The proof of this result will be described briefly in the next subsection.
\vskip 0.2 cm

One can now describe how the proof of Theorem \ref{th.GPS3} is organized in \cite{\GPS}. It is divided into 
three main parts:

\begin{enumerate}
\item The first part deals with proving the multi-scale sub-exponential bound on the lower-tail of $|\Spec_{(r)}|$ given by Theorem \ref{th.subexpoS}.
\item The second part consists in proving as much as we can on the dependency structure of $\Spec_{f_n}$. Unfortunately here,
it seems to be very challenging to achieve a good understanding of all the ``independence'' that should be present within $\Spec_{f_n}$.
The only hint of independence which was finally proved in \cite{\GPS} is a very {\it weak} one (see subsection \ref{ss.hint}).
In particular, it is too weak to readily imply a bound 
like ~\eqref{e.ifindependent}.

\item Since disjoint regions of the spectral sample $\Spec_{f_n}$ are not known to behave independently of each other, the third part of the proof consists 
in adapting the setup we used for the tree (where distinct branches evolve exactly independently of each other) into a setup  where the weak hint of independence 
obtained in the second part of the program turns out to be enough to imply the bound given by 
~\eqref{e.ifindependent} for an appropriate absolute constant $c>0$. This final part of the proof will be discussed in subsection \ref{ss.LDlemma}.

\end{enumerate}

The next three subsections will be devoted to each of these 3 parts of the program.

\subsection{Some words about the sub-exponential bound on the lower tail of $\Spec_{(r)}$}

In this subsection, we turn our attention to the proof of the first part of the program, i.e. on Theorem \ref{th.subexpoS}.
In fact, as in the case of $\Tree_n$, the following more explicit statement is proved in \cite{\GPS}.

\begin{theorem}[\cite{\GPS}]\label{th.GPSgk}
There exists an absolute constant $\theta>0$ such that for any $1\le r \le n$ and any $k\geq 1$,
\begin{align*}
\SPb{|\Spec_{(r)}| = k }\le \, 2^{\, \theta  \log_2^2(k+2)} \; \SPb{|\Spec_{(r)}| = 1}\,.
\end{align*}
\end{theorem}

\begin{remark}
Note that the theorems from \cite{\BKS} on the noise sensitivity of percolation are all particular cases ($r=1$) of this intermediate result in \cite{\GPS}. 
\end{remark}

\vskip 0.2 cm

The main idea in the proof of this theorem is in some sense to assign a {\it tree structure} to each possible set $\Spec_{(r)}$.
The advantage of working with a tree structure is that it is easier to work with inductive arguments. In fact, once
a mapping $\Spec_{(r)}\mapsto \text{``tree structure''}$ has been designed, the proof proceeds similarly as in the case of $\Tree_{(r)}$ by double induction on 
the depth of the tree as well as on $k\geq 1$. Of course, this mapping is a delicate affair:
it has to be designed in an ``efficient'' way so that it can compete against entropy effects caused by the exponential growth of the number of tree structures.

We will not give the details of how to define such a mapping, but let us describe informally how it works. 
More specifically than a tree structure, we will in fact  assign an {\it annulus structure} to each set $\Spec_{(r)}$. 

\begin{definition}\label{d.cluster}
Let $\Ann$ be a finite collection of disjoint (topological) annuli in the plane.
We call this an {\bf annulus structure}. 
Furthermore, we will say that a set $S \subset \R^2$ is {\bf compatible} with $\Ann$ (or vice versa) if it is contained in $\R^2\setminus\bigcup\Ann$ and intersects the inner disk of each annulus in $\Ann$. Note that it is allowed that
one annulus is ``inside'' of another annulus.
\end{definition}

\begin{figure}[!h]
\begin{center}
\includegraphics[width=0.6 \textwidth]{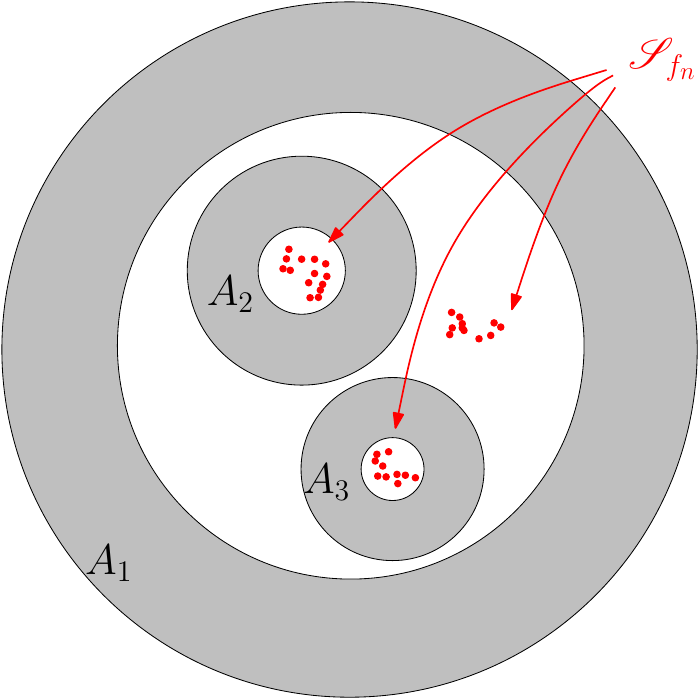}
\end{center}
\caption{An example of an {\bf annulus structure} $\Ann:= \{ A_1,A_2,A_3\}$ compatible with a spectral sample $\Spec_{f_n}$.}\label{f.annularstructure}
\end{figure}

The mapping procedure in \cite{\GPS} assigns to each $\Spec_{(r)}$ an annulus structure $\Ann \subset [0,n]^2$ in such a way that it is compatible 
with $\Spec_{(r)}$. See Figure \ref{f.annularstructure} for an example.
Again, we will not describe this procedure nor discuss the obvious boundary issues which arise here, but let us state a crucial
property satisfied by annulus structures.

\begin{lemma}[\bf clustering Lemma]\label{l.cluster}
If $\Ann$ is an annulus structure contained in $[0,n]^2$, then 
\begin{align*}
\SPb{\Spec_{(r)}\text{ is compatible with }\Ann} \le \prod_{A \in \Ann} \alpha_4(A)^2\,,
\end{align*}
where $\alpha_4(A)$ denotes the probability of having a {\it four-arm} event in the annulus $A$.
\end{lemma}

\begin{remark}
To deal with boundary issues, one would also need to incorporate within our annulus structures half-annuli centered on the boundaries as well as  
quarter disks centered at the corners of $[0,n]^2$. 
\end{remark}

Let us briefly comment on this lemma.
\bi
\item First of all, its proof is an elegant combination of linear algebra and percolation. It is a short and relatively elementary argument. See Lemma 4.3 in \cite{\GPS}.
\item It is very powerful in dealing with the possible non-injectivity of the mapping $\mathcal{S}_{(r)}\mapsto \Ann$. Indeed, 
while describing the setup above, one might have objected that if the mapping were not injective enough, then 
the cardinality of the ``fibers'' above each annulus structure would have to be taken into account as well. Fortunately, the above lemma reads as follows:
for any fixed annulus structure $\Ann$,
\begin{align*}
\sum_{\Spec_{(r)}:\, \Spec_{(r)} \mapsto \Ann} \SPb{\Spec_{(r)}} & \le \SPb{ \Spec_{(r)}\text{ is compatible with }\Ann} \le \prod_{A\in \Ann} \alpha_4(A)^2\,.
\end{align*}

\item Another essential feature of this lemma is that it quantifies very efficiently the fact that the clustering effect wins over the entropy effect in the sense 
of Figure \ref{t.entropy}. The mechanism responsible for this is that the probability of the {\it four-arm} event squared has an exponent (equal to $5/2$ on $\T$)
larger than the {\it volume} exponent equal to 2. To illustrate this, let us analyze the situation when $k=2$ (still neglecting boundary issues). The probability that 
the spectrum $\Spec_{f_n}$ intersects two and only two $r$-squares at macroscopic distance $\Omega(n)$ from each other can be easily estimated using the lemma. Indeed, in such a case, 
$\Spec_{(r)}$ would be compatible with an annulus structure consisting of two annuli, each being approximately of the type $A(r,n)$. There are $O(\frac {n^2}{r^2})\times O(\frac {n^2}{r^2})$ such possible annulus structures. Using the lemma each of them costs (on $\T$) $(\frac{r} n)^{5+o(1)}$. An easy exercise shows that this is much smaller than $\SPb{|\Spec_{(r)}|=2}$. In other words, if 
$|\Spec_{(r)}|$ is conditioned to be small, it tends to be localized. Also,
the way that the lemma is stated makes it very convenient to work with 
higher values of $k$.
\ei

The details of the proof of Theorem \ref{th.GPSgk} can be found in \cite{\GPS}. The double induction there is in some sense very close to the one we carried out in detail
in subsection \ref{ss.Treesubexpo} in the case of the tree; this is the reason why we included this latter proof.  For those who might read the proof in \cite{\GPS}, there is a 
notion of {\em overcrowded cluster} defined there; it exactly corresponds in the case of the tree to stopping the analysis above scale $J_k$ instead of going all the way to scale 1
(note that without stopping at this scale $J_k$, the double induction in subsection \ref{ss.Treesubexpo} would have failed).

\subsection{Some words on the weak independence property proved in \cite{\GPS}}\label{ss.hint}

This part of the program is in some sense the main one. 
To introduce it, let us start by a naive but tempting strategy. What the first part of the program (Theorem \ref{th.GPSgk}) tells us is that for any mesoscopic scale $1\le r \le n$, 
if $\Spec_{f_n}$ is non-empty, it is very unlikely that it will intersect few squares of size $r$. In other words, it is very unlikely that 
$|\Spec_{(r)}|$ will be small. Let $B_1, \ldots, B_m$ denote the set of $O(n^2/r^2)$ $r$-squares which tile $[0,n]^2$. One might try the following {\em scanning procedure}:
explore the spectral sample $\Spec_{f_n}$ inside the squares $B_i$ one at a time. More precisely, before starting the scanning procedure,
we consider our spectral sample $\Spec_{f_n}$ as a random subset of $[0,n]^2$ 
about which we do not know anything yet. Then, at step one, we 
reveal $\Spec_{| B_1}$. This gives us some partial information about
$\Spec_{f_n}$. What we still have to explore is a random set of $[0,n]^2 \setminus B_1$ which follows 
the law of a spectral sample conditioned on what was seen in $B_1$ and we keep going in this way. By Theorem \ref{th.GPSgk}, many of these squares will be non-empty. Now, it is not hard to
prove the following lemma (using similar methods as in Problem \ref{ex.scaleinv}).
\begin{lemma}\label{l.conditiononebox}
There is a universal constant $c>0$ such that for any $r$-square $B$ in the bulk $[n/4, 3n/4]^2$, one has 
\begin{align*}
\SPb{|\Spec_{f_n} \cap B| > c \,  r^2 \alpha_4(r) \, \bigm| \, \Spec_{f_n} \cap B \neq \emptyset} > c\,.
\end{align*}
\end{lemma}
This lemma in fact holds uniformly in the position of $B$ inside $[0,n]^2$. 

If one could prove the following (much) stronger result: there exists a universal constant $c>0$ such that 
uniformly on the sets $S \subset [0,n]^2 \setminus B$ one has 

\begin{align}\label{e.expectedconditionalbehavior}
\SPb{|\Spec_{f_n} \cap B| > c \,  r^2 \alpha_4(r) \, \bigm| \, \Spec_{f_n} \cap B \neq \emptyset \text{ and } \Spec_{| B^c} = S} > c\,,
\end{align}
then it would not be hard to make the above scanning strategy work together with Theorem \ref{th.GPSgk} in order to obtain Theorem \ref{th.GPS3}.
(Note that such a result would indeed give a strong hint of independence within $\Spec_{f_n}$.) 
However, as we discussed before, the current understanding of the independence within $\Spec_{f_n}$ is far from giving such a  
statement. Instead, the following result is proved in \cite{\GPS}. We provide here a slightly simplified version. 

\begin{theorem}[\cite{\GPS}]\label{th.WI}
There exists a uniform constant $c>0$ such that for any set $W\subset [0,n]^2$ and any $r$-square $B$ such that $B\cap W = \emptyset$,
one has 
\begin{align*}
\SPb{|\Spec_{f_n}\cap B| > c\, r^2 \alpha_4(r) \md \Spec_{f_n} \cap B \neq \emptyset \text{ and }\Spec_{f_n} \cap W = \emptyset}>c\,.
\end{align*}
\end{theorem}

Note that this theorem in some sense interpolates between part of Lemma \ref{l.prob} and Lemma \ref{l.conditiononebox}
which correspond respectively to the special cases $W=B^c$ and $W=\emptyset$. 
Yet it looks very weak compared to the expected
(\ref{e.expectedconditionalbehavior})
which is stated uniformly on the behavior of $\Spec_{f_n}$ outside of $B$.

Assuming this weak hint of independence (Theorem \ref{th.WI}), it seems we are in bad shape 
if we try to apply the above scanning procedure. Indeed, we face the following two obstacles:
\begin{enumerate}
\item  The first obstacle is that one would keep a good control only as far as one would not see any ``spectrum''.
Namely, while revealing $\Spec_{\md B_i}$ one at a time, the first time one finds a square $B_i$ such that $\Spec_{\md B_i}\neq \emptyset$, one would be forced to stop
the scanning procedure there. In particular,
if the size of the spectrum in this first non-trivial square does not exceed $r^2 \alpha_4(r)$, then we cannot conclude anything. 

\item The second obstacle is that, besides the conditioning $\Spec\cap W =\emptyset$, our estimate is also conditioned on the event that $\Spec\cap B \neq \emptyset$. In particular, in the above ``naive'' scanning
strategy where squares are revealed in a sequential way, at each step one would have to update the probability that $\Spec \cap B_{i+1} \neq \emptyset$
based on what was discovered so far. 
\end{enumerate}

It is the purpose of the third part of the program to adapt the above scanning strategy to these constraints. Before describing this third part in the next subsection,
let us say a few words on how to prove Theorem  \ref{th.WI}. 

A crucial step in the proof of this theorem 
is to understand the following ``one-point function'' for any $x\in B$ at distance at least 
$r/3$ from the boundary:
\[
\SPb{x\in \Spec_{f_n}\text{ and } \Spec_{f_n} \cap W = \emptyset} \,.
\]
A very useful observation is to rewrite this one-point function in terms of an explicit coupling of two 
i.i.d.\ percolation configurations.
It works as follows: let $(\omega_1,\omega_2)$ be a coupling of two i.i.d.\ percolations on $[0,n]^2$ which are such that
\[
\left\lbrace \begin{array}{ll} \omega_1 = \omega_2 & \text{ on } W^c \\ \omega_1, \omega_2 & \text{ are independent on } W \end{array}\right.
\]
One can check that the one-point function we are interested in is related to this coupling in the following simple way:
\[
\SPb{x\in \Spec_{f_n}\text{ and } \Spec_{f_n} \cap W = \emptyset} = \Pb{x \text{ is pivotal for }\omega_1 \text{ AND }\omega_2}\,.
\]

\begin{remark}
You may check this identity in the special cases where $W=\emptyset$ or $W=\{x\}^c$.
\end{remark}
Thanks to this observation, the proof of Theorem \ref{th.WI} proceeds by analyzing this $W$-coupling. See \cite{\GPS} for the complete details.

\subsection{Adapting the setup to the weak hint of independence}\label{ss.LDlemma}

As we discussed in the previous subsection, one faces two main obstacles if, on the basis of the weak independence given by Theorem \ref{th.WI},
one tries to apply the naive sequential scanning procedure described earlier. 

Let us start with the first obstacle. Assume that we scan the domain $[0,n]^2$ in a sequential way, i.e., we choose an increasing family of subsets $(W_l)_{l\geq 1} = (\{w_1,\ldots, w_l\})_{l\geq 1}$.
At each step, we reveal what $\Spec_{| \{w_{l+1}\}}$ is, conditioned on what was discovered so far (i.e., conditioned on $\Spec_{| W_l}$). 
From the weak independence Theorem \ref{th.WI}, it is clear that if we want this strategy to have any chance to be successful, we have to choose 
$(W_l)_{l\geq1}$ in such a way that $(\Spec_{f_n} \cap W_l)_{l\geq 1}$ will remain empty for some time (so that we can continue to rely on our weak independence result); of course 
this cannot remain empty forever, so the game is to choose the increasing family $(W_l)_{l\geq1}$
in such a way that the first time $\Spec_{f_n} \cap \{w_l\}$ will happen to be non-empty, it should give a strong indication that $\Spec_{f_n}$ is large in the $r$-neighborhood 
of $w_l$. 

As we have seen, revealing the entire mesoscopic boxes $B_i$ one at a time is not a successful idea. Here is a much better idea (which is not yet the right one due to the second obstacle, but we are getting close): in each $r$-square $B_i$, instead of revealing all the bits, let us reveal only a very small proportion $\delta_r$ of them.
Lemma \ref{l.conditiononebox} tells us that if $\Spec\cap B_i \neq \emptyset$, then each point $x\in B_i$ has probability of order $\alpha_4(r)$ to be in $\Spec_{f_n}$.
Therefore if we choose $\delta_r \ll (r^2 \alpha_4(r))^{-1}$, then with high probability, by revealing only a proportion $\delta_r$ of the points in $B_i$, we will ``miss'' the 
spectral sample $\Spec_{f_n}$. Hence, we have to choose $\delta_r \ge (r^2 \alpha_4(r))^{-1}$. In fact choosing $\delta \asymp (r^2 \alpha_4(r))^{-1}$ is exactly the right balance. 
Indeed, we know from Theorem \ref{th.GPSgk} that many $r$-squares $B_i$ will be touched by the spectral sample; now, in this more sophisticated scanning procedure, if the first such square encountered happens to contain few points (i.e. $\ll r^2 \alpha_4(r)$), then with the previous scanning strategy, we would 
``lose'',
but with the present one, due to our choice of $\delta_r$, most likely we will keep $\Spec_{f_n} \cap W_l = \emptyset$ so that we can continue further on until we reach 
a ``good'' square (i.e. a square containing of order $r^2 \alpha_4(r)$ points). 

Now, Theorems \ref{th.GPSgk} and \ref{th.WI} together tell us that with high probability, one will eventually reach such a good square. Indeed, suppose the $m$ first $r$-squares touched by the spectral sample happened to contain few points; then, most likely, if $W_{l_m}$ is the set of bits 
revealed so far, by our choice of $\delta_r$ we will still have $\Spec\cap W_{l_m} = \emptyset$. This allows us to still rely on Theorem \ref{th.WI}, which basically tells
us that there is a positive conditional probability for the next one to be a ``good'' square (we are neglecting the second obstacle here).
This says that the probability to visit $m$ consecutive bad squares seems to decrease exponentially fast. Since $m$ is typically very large (by Theorem \ref{th.GPSgk}), we conclude that, 
with high probability, we will finally reach good squares. In the first good square encountered, by our choice of $\delta_r$, there is now a positive probability to reveal a bit present in 
$\Spec_{f_n}$. In this case, the sequential scanning will have to stop, since we will not be able to use our weak independence result anymore, but this is not a big issue: indeed,
assume that you have some random set $S\subset B$. If by revealing each bit only with probability $\delta_r$, you end up finding a point in $S$, most likely your set $S$ 
is at least of size $\Omega(r^2 \alpha_4(r))$. This is exactly the size we are looking for in Theorem \ref{th.GPS3}.
\vskip 0.1 cm

Now, only the second obstacle remains. It can be rephrased as follows: assume you applied the above strategy in $B_1,\ldots, B_h$ (i.e. you revealed each point in $B_i, \, i\in\{1,\ldots, h\}$
only with probability $\delta_r$) and that you did not find any spectrum yet. In other words, if $W_l$ denotes the set of points visited so far, then $\Spec_{f_n} \cap W_l = \emptyset$. Now if $B_{h+1}$ is the next $r$-square to be scanned (still in a ``dilute'' way with intensity $\delta_r$), we seem to be in good shape since we know how to control the conditioning $\Spec_{f_n}\cap W_l = \emptyset$. However, if we want to rely on the uniform control given by Theorem \ref{th.WI}, we also need to further condition on $\Spec_{f_n} \cap B_{h+1} \neq \emptyset$. In other words, we need to control the following conditional expectation:
\[
\SPb{\Spec_{f_n} \cap B_{h+1} \neq \emptyset \bigm| \Spec_{f_n} \cap W_l = \emptyset}\,.
\]
It is quite involved to estimate such quantities. Fortunately, by changing our sequential scanning procedure into a slightly more ``abstract''  procedure, one can avoid 
dealing with such terms. More precisely, within each $r$-square $B$, we will still reveal only a $\delta_r$ proportion of the bits (so that the first obstacle is still taken care of),
but instead of operating in a sequential way (i.e. scanning $B_1$, then $B_2$ and so on), we will gain a lot by considering the combination of Theorem \ref{th.GPSgk} and Theorem \ref{th.WI} in a more abstract fashion. Namely, the following large deviation lemma from \cite{\GPS} captures exactly what we need in our present situation.

\begin{lemma}[\cite{\GPS}]\label{l.LD}
Let $X_i,Y_i \in \{0,1 \},\, i\in \{ 1,\ldots, m\}$  be random variables such that for each $i$
$Y_i \le X_i$ a.s.
If $\forall J\subset [m]$ and $\forall i\in [m]\setminus J$, we have
\begin{align}\label{e.LD}
\Pb{Y_i = 1 \md Y_j = 0,  \, \forall j\in J} \geq c \; \Pb{X_i =1 \md Y_j = 0, \,  \forall j \in J},
\end{align}
then if $X:= \sum X_i$ and $Y:= \sum Y_i$, one has that 
\begin{align*}
\Pb{Y=0 \, \md X>0} \le c^{-1} \Eb{ e^{-(c/e) X} \, \md X>0}\,.
\end{align*}
\end{lemma}

Recall that $B_1,\ldots, B_m$ denotes the set of $r$-squares which tile $[0,n]^2$.
For each $i\in [m]$, let $X_i := 1_{\Spec \cap B_i \neq \emptyset} $ and $Y_i:= 1_{\Spec \cap B_i \cap \mathcal{W} \neq \emptyset} $, where $\mathcal{W}$ is an independent uniform random 
subset of $[0,n]^2$ of intensity $\delta_r$.  

This lemma enables us to combine our two main results, Theorems \ref{th.WI} and 
\ref{th.GPSgk}, in a very nice way:
By our choice of the intensity $\delta_r$, Theorem \ref{th.WI} exactly states that the 
assumption~\eqref{e.LD} is satisfied for a certain constant $c>0$. 
Lemma \ref{l.LD} then implies that 

\begin{align*}
\SPb{ Y=0\, \md X>0 } \le c^{-1} \Eb{e^{- (c/e) X} \, \md X>0 }\,.
\end{align*}
Now, notice that $X=\sum X_i$ exactly corresponds to $|\Spec_{(r)}|$ while the event $\{X>0\}$ corresponds to $\{ \Spec_{f_n} \neq \emptyset \}$
and the event $\{Y=0\}$ corresponds to $\{ \Spec_{f_n} \cap \mathcal{W} = \emptyset \}$. Therefore Theorem \ref{th.GPSgk} leads us to

\begin{align}\label{eq.final}
\SPb{\Spec_{f_n} \cap \mathcal{W} = \emptyset\,,\,  \Spec_{f_n} \neq \emptyset } & \le c^{-1} \Eb{e^{- (c/e) |\Spec_{(r)}| } \,, \, \Spec_{f_n} \neq \emptyset} \nonumber \\ 
& \le c^{-1} \sum_{k \geq 1} \SPb{|\Spec_{(r)}| = k}  e^{-(c/e) k} \nonumber \\
&\le c^{-1} \Bigl( \sum_{k\geq 1} 2^{\theta \log_2^2(k+2)} e^{-(c/e) k)} \Bigr) \SPb{|\Spec_{(r)}| = 1} \nonumber  \\
&\le C(\theta) \, \SPb{|\Spec_{(r)}|= 1} \asymp \frac {n^2}{r^2} \alpha_4(r,n)^2\, ,
\end{align}
where (\ref{e.4armsquared}) is used in the last step.

This shows that on the event that $\Spec_{f_n}\neq \emptyset$, it is very unlikely that we do not detect the spectral sample on the $\delta_r$-dilute set $\mathcal{W}$.
This is enough for us to conclude using the following identity:
\[
\SPb{ \Spec_{f_n} \cap \mathcal{W} = \emptyset \, \md \Spec_{f_n}} = (1- \delta_r)^{|\Spec_{f_n}|} = (1- \frac 1 {r^2 \alpha_4(r)})^{|\Spec_{f_n}|}\,.
\] 
Indeed, by averaging this identity we obtain 
\begin{align*}
\SPb{\Spec_{f_n} \cap \mathcal{W} = \emptyset\,, \Spec_{f_n} \neq \emptyset} & = \SEb{ \SPb{\Spec_{f_n} \cap \mathcal{W} = \emptyset \, \md \Spec_{f_n}} \, 1_{\Spec_{f_n} \neq \emptyset}} \\
& = \SEb{(1- \frac 1 {r^2 \alpha_4(r)})^{|\Spec_{f_n}|}\;  1_{\Spec_{f_n} \neq \emptyset}  }\\ 
& \ge \Omega(1) \SPb{0< |\Spec_{f_n}| < r^2 \alpha_4(r)}\,,
\end{align*}
which, combined with~\eqref{eq.final} yields the desired upper bound in Theorem \ref{th.GPS3}. See Problem \ref{prob.sharp}
for the lower bound.

\section{The radial case}

The next chapter will focus on the existence of {\it exceptional times} in the model of dynamical percolation.
A main tool in the study of these exceptional times is the spectral measure $\SQ_{g_R}$ where $g_R$ is the Boolean function 
$g_R := \{ -1 ,1 \}^{O(R^2)} \to \{0,1\}$ defined to be the indicator function of the one-arm event $\{0 \longleftrightarrow \p B(0,R) \}$.
Note that by definition, $g_R$ is such that $\| g_R \|_2^2 = \alpha_1(R)$.

In \cite{\GPS}, the following ``sharp'' theorem on the lower tail of $\Spec_{g_R}$ is proved.

\begin{theorem}[\cite{\GPS}]\label{th.spectrumonearm}
Let $g_R$ be the one-arm event in $B(0,R)$.
Then for any $1\le r \le R$, one has 
\begin{align}
\SQf{g_R}{0< |\Spec_{g_R}|<r^2 \alpha_4(r)} \asymp \frac{\alpha_1(R)^2}{\alpha_1(r)}\,.
\end{align}
\end{theorem}

The proof of this theorem is in many ways similar to the chordal case (Theorem \ref{th.GPS3}).
An essential difference is that the ``clustering v.s. entropy'' mechanism is very different in this case.
Indeed in the chordal left to right case, when $\Spec_{f_n}$ is conditioned to be very small, the proof of Theorem \ref{th.GPS3}
shows that typically $\Spec_{f_n}$ localizes in some $r$-square whose location is ``uniform'' in the domain $[0,n]^2$. 
In the radial case, the situation is very different: $\Spec_{g_R}$ conditioned to be very small will in fact tend to localize in the $r$-square 
centered at the origin. This means that the analysis of the mesoscopic behavior (i.e. the analogue of Theorem \ref{th.GPSgk}) has to be adapted to the radial case.
In particular, in the definition of an annulus structure, the annuli containing the origin play a distinguished role. See \cite{\GPS} for complete details.

\chapter*{Exercise sheet on chapter \ref{ch.SNS}}

\setcounter{exercise}{0}

\begin{exercise}
Prove Proposition \ref{pr.tree}.
\end{exercise}

\begin{exercise}\label{ex.induction}
Consider the fractal percolation process $\Tree^i$, $i \geq 1$ introduced in this chapter.
(Recall that $\Tree_{2^i}\equiv \Tree^i$).
Recall that in Section \ref{s.toy}, it was important to estimate the quantity $\Pb{|\Tree^i|=1 }$. 
This is one of the purposes of the present exercise.

\bi

\item[(a)] Let $p_i:= \Pb{|\Tree^i | = 1}$. By recursion, show that there is a constant $c\in (0,1)$ so that, as $i\to \infty$ 
\[
p_i \sim c \mu^i\,,
\]
where $\mu:= 4p (1-p + pq)^3$ and $q$ is the probability of extinction for the Galton-Watson tree correponding to $(\Tree^i)_{i\geq 1}$.

\item[(b)] Using the generating function 
$s\mapsto f(s)(=E(s^{\mbox{ number of offspring}})$ 
of this Galton-Watson
tree, and by studying the behavior of its $i$-th iterates $f^{(i)}$, prove the 
same result with $\mu:= f'(q)$. Check that it gives the same formula.		

\item[(c)] Recall the definition of $p_{m,b}$ from Section \ref{s.toy}. 
Let $p_{m,\infty}$ be the probability that exactly 1 person at generation $m$
survives forever. Prove that 
\[
p_{m,\infty} = (1-q)\mu^m
\]
for the same exponent $\mu$. Prove Lemma \ref{l.pmh}. 
Finally, prove that 
$\lim_{b\to\infty}p_{m,b}=p_{m,\infty}$.
\ei

\end{exercise}

\begin{exercise}\label{ex.entropy}
Extract from the proof of Lemma \ref{l.treeGk2} the answer to the question asked in 
Figure \ref{t.entropy}.
\end{exercise}

\begin{exercise}\label{ex.implication}
Prove that 
\[
\text{Theorem \ref{th.GPS3}} \Rightarrow \text{Theorem \ref{th.GPS2}} \Rightarrow \text{Theorem \ref{th.GPS}} 
\]
\end{exercise}

\begin{exercise}\label{ex.onebox} 
Consider an $r$-square $B\subset [n/4, 3n/4]^2$ in the ``bulk'' of $[0,n]^2$.

\bi
\item[(a)] Prove using Proposition \ref{pr.easy.spectral} that 
\[
\SPb{\Spec_{f_n} \neq \emptyset \text{ and } \Spec_{f_n} \subset B} \asymp \alpha_4(r,n)^2
\]
\item[(b)] Check that the clustering Lemma \ref{l.cluster} is consistent with this estimate.
\ei

\end{exercise}

\begin{problem}\label{ex.scaleinv}
The purpose of this exercise is to prove Lemma  \ref{l.prob}.

\bi
\item[(a)] Using Proposition \ref{pr.easy.spectral}, prove that for any $x\in B$ at distance $r/3$ from the boundary, 
\[
\Pb{x\in \Spec_{f_n} \text{ and } \Spec_{f_n} \cap B^c =\emptyset} \asymp \alpha_4(r) \alpha_4(r,n)^2\,.
\]

\item[(b)] Recover the same result using Proposition \ref{pr.spectralsubdomain} instead.

\item[(c)] Conclude using Exercise \ref{ex.onebox} that $\SEb{|\Spec_{f_n}\cap \bar B| \md \Spec_{f_n} \neq \emptyset \text{ and } \Spec_{f_n} \subset B} \asymp r^2 \alpha_4(r)$, where $\bar B \subset B$ is the set of points $x\in B$ at distance at least $r/3$ from the boundary. 

\item[(d)] Study the second-moment $\SEb{|\Spec_{f_n}\cap \bar B|^2 \md \Spec_{f_n} \neq \emptyset \text{ and } \Spec_{f_n} \subset B}$.

\item[(e)] Deduce Lemma \ref{l.prob}.

\ei

\end{problem}

\begin{problem}\label{prob.sharp}

Most of this chapter was devoted to the explanation of the proof of Theorem \ref{th.GPS3}.
Note that we in fact only discussed how to prove the upper bound. This is because the lower bound is much easier to prove and this is the purpose of this problem.

\ni
\bi
\item[(a)] Deduce from Lemma \ref{l.prob} 
and Exercise \ref{ex.onebox}(a)
that the lower bound on $\SPb{0< |\Spec_{f_n}| < r^2 \alpha_4(r)}$ given in Theorem \ref{th.GPS3}
is correct. I.e., show that there exists a constant $c>0$ such that 
\begin{align*}
\SPb{0< |\Spec_{f_n}| < r^2 \alpha_4(r)} > c \frac{n^2}{r^2} \alpha_4(r,n)^2\,.
\end{align*}

\item[(b)] (Hard) In the same fashion, prove the lower bound part of Theorem \ref{th.spectrumonearm}.
\ei

\end{problem}

\chapter{\;\; Applications to dynamical percolation}\label{ch.DP}

\note{Timing}{/////////// 1.5 hour //////////// \vskip 1 cm}

\medbreak

In this section, we present a very natural model 
where percolation undergoes a time-evolution: this is the model of 
{\bf dynamical percolation} described below. The study of the 
``dynamical'' behavior of percolation as opposed to its ``static'' behavior
turns out to be very rich: interesting phenomena arise especially 
at the phase transition point. We will see that in some sense,
dynamical planar percolation at criticality is a 
very unstable or chaotic process. In order to understand this instability, 
sensitivity of percolation (and therefore 
its Fourier analysis) will play a key role.
In fact, the original motivation for the paper \cite{\BKS} on
noise sensitivity was to solve a particular problem 
in the subject of dynamical percolation. 
\cite{\SurveySteif} provides a recent survey on the subject of
dynamical percolation. 

We mention that one can read all but the last section of the present chapter 
without having read Chapter \ref{ch.SNS}.

\section{The model of dynamical percolation}

This model was introduced by 
H\"{a}ggstr\"{o}m, Peres and Steif~\cite{\HPS} 
inspired by a question that Paul Malliavin
asked at a lecture at the Mittag-Leffler Institute in 1995.
This model was invented independently by Itai Benjamini.

In the general version of this model as it was introduced, 
given an arbitrary graph $G$ and a parameter $p$,
the edges of $G$ switch back and forth according to independent 2-state
continuous time Markov chains where closed switches to open at rate $p$
and open switches to closed at rate $1-p$. Clearly, the product
measure with density $p$, denoted by $\pi_p$ in this chapter, is the 
unique stationary distribution for this Markov process. The general question 
studied in dynamical percolation is whether, when we start with the 
stationary distribution $\pi_p$, there exist atypical times at which the percolation
structure looks markedly different than that at a fixed time. In almost all
cases, the term  ``markedly different'' refers to the existence or
nonexistence of an infinite connected component. Dynamical percolation
on site percolation models, which includes our most important
case of the hexagonal lattice, is defined analogously.

We very briefly summarize a few early results in the area. 
It was shown in \cite{\HPS} that below criticality, there are no times
at which there is an infinite cluster and above criticality, there is
an infinite cluster at all times. See the exercises.
In \cite{\HPS}, examples of graphs which do not percolate at
criticality but for which there exist exceptional times where
percolation occurs were given. 
(Also given were examples of graphs which do percolate at
criticality but for which there exist exceptional times where
percolation does not occur.)
A fairly refined analysis of the case of so-called
{\em spherically symmetric} trees was given. See the exercises
for some of these.

Given the above results, it is natural to ask what 
happens on the standard graphs that we work with. Recall that for 
$\Z^2$, we have seen that there is no percolation at criticality.
It turns out that it is also known (see below) that for
$d\geq 19$, there is no percolation at criticality for $\Z^d$.
It is a major open question to prove that this is also the case for
intermediate dimensions; the consensus is that this should be the case.

\section{What's going on in high dimensions: $\Z^d, d\geq 19$?}

For the high dimensional case, $\Z^d, d\geq 19$, it was shown in
\cite{\HPS} that there are no exceptional times of percolation at criticality. 

\begin{thm}[\cite{\HPS}] \label{th:Zd}
For the integer lattice $\Z^d$ with $d\ge 19,$ dynamical critical percolation
has no exceptional times of percolation.
\end{thm}

The key reason for this is a highly nontrivial result due to work of
Hara and Slade (\cite{\HaraSladeLace}), 
using earlier work of Barsky and Aizenman (\cite{\Triangle}),
that says that if $\theta(p)$ is the probability that the origin percolates
when the parameter is $p$, then for $p\ge p_c$
\begin{equation} \label{eq:derive}
\theta(p)= O(p-p_c) \, .   
\end{equation}

(This implies in particular that there is no percolation at criticality.)
In fact, this is the only thing which is used in the proof and hence
the result holds whenever the percolation function satisfies this
``finite derivative condition'' at the critical point.

\proofoutline
By countable additivity, it suffices to show that there are no times
at which the origin percolates during $[0,1]$.
We use a first moment argument.
We break the time interval $[0,1]$ into $m$ intervals each of 
length $1/m$. If we fix one of these intervals,
the set of edges which are open at {\em some time} during this interval
is i.i.d.\ with density about $p_c+1/m$. Hence the probability that the
origin percolates with respect to these set of edges is by 
(\ref{eq:derive}) at most $O(1/m)$. It follows that 
if $N_m$ is the number of intervals where this occurs, then
$\E[N_m]$ is at most $O(1)$. 
It is not hard to check
that $N\le \liminf_m N_m$, where $N$ is the cardinality of the set 
of times during 
$[0,1]$ at which the origin percolates. Fatou's Lemma now yields
that $\E(N)<\infty$ and hence there are at most finitely many
exceptional times during $[0,1]$ at which the origin percolates. To go from
here to having no exceptional times can either be done by using some rather 
abstract Markov process theory or by a more hands on approach
as was done in \cite{\HPS} and which we refer to for details.
\qed

\begin{remark}
It is known that (\ref{eq:derive}) holds for any homogeneous tree
(see \cite{\Grimm} for the binary tree case) 
and hence there are no exceptional times of percolation in this case also.
\end{remark}

\begin{remark}
It is was proved by Kesten and Zhang \cite{\StrictInequalities},
that (\ref{eq:derive}) fails for $\Z^2$ and hence the proof method
above to show that there are no exceptional times
fails. This infinite derivative in this case might suggest that
there are in fact exceptional times for critical dynamical percolation on
$\Z^2$, an important question left open in \cite{\HPS}. 
\end{remark}

\section{$d=2$ and BKS}
One of the questions posed in \cite{\HPS} was whether there are
exceptional times of percolation for $\Z^2$.
It was this question which was one of the main motivations for the
paper \cite{\BKS}. While they did not prove the existence of
exceptional times of percolation, they did obtain the following very
interesting result which has a very similar flavor.

\begin{theorem}[\cite{\BKS}]\label{t.bksalmost}
Consider an $R\times R$ box on which we run critical dynamical
percolation. Let $S_R$ be the number of times during $[0,1]$
at which the configuration changes from having a percolation crossing
to not having one. Then 
$$
S_R\to\infty \mbox{ in probability as } R\to\infty.
$$
\end{theorem}

\medskip\noindent
Noise sensitivity of percolation as well as 
the above theorem tells us that certain large scale connectivity properties
decorrelate very quickly. This suggests that in some vague sense
$\omega_t^{p_c}$ ``changes'' very quickly as
time goes on and hence there might be some chance that
an infinite cluster appears since we are given many ``chances''.

In the next section, we begin our study of exceptional times
for  $\Z^2$ and the hexagonal lattice.

\section{The second moment method and the spectrum}

In this section, we reduce the question of exceptional times to a
``second moment method'' computation which in turn reduces to 
questions concerning the spectral behavior for specific 
Boolean functions involving percolation. Since $p=1/2$, our dynamics
can be equivalently defined by having each edge or hexagon be
rerandomized at rate 1.

The key random variable which one needs to look at is
$$
X=X_R:= \int_0^1 1_{0\overset{\omega_t}{\longleftrightarrow}R} \, dt 
$$
where $0\overset{\omega_t}{\longleftrightarrow}R$ is of course 
the event that at time $t$ there is an open 
path from the origin to distance $R$ away. Note that the above integral
is simply the Lebesgue measure of the set of times in $[0,1]$ 
at which this occurs.

We want to apply the second moment method here. We isolate the easy
part of the argument so that the reader who is not familiar
with this method understands it in a more general context.
However, the reader should keep in mind that the difficult part
is always to prove the needed bound on the second moments which in this
case is (\ref{eq:needed.inequality}).

\begin{prop} \label{pr:2nd.moment.bound}
If there exists a constant $C$ such that for all $R$
\begin{equation} \label{eq:needed.inequality}
\E(X_R^2)\le C\E(X_R)^2,
\end{equation}
then a.s.\ there are exceptional times of percolation.
\end{prop}

\proof
For any nonnegative random variable $Y$, the Cauchy-Schwarz inequality 
applied to $YI_{\{Y>0\}}$ yields
$$
\P(Y >0)\ge \E(Y)^2/\E(Y^2).
$$
Hence by (\ref{eq:needed.inequality}), we have that for all $R$,
$$
\P(X_R >0)\ge 1/C
$$
and hence by countable additivity (as we have a decreasing sequence of events)
$$
\P(\cap_R \{X_R >0\})\ge 1/C.
$$
Had the set of times that a fixed edge is on been a closed set, then
the above would have yielded by compactness that there is
an exceptional time of percolation with probability at least 
$1/C$. However, this is not a closed set. On the other hand, this
point is very easily fixed by modifying the process so that the times
each edge is on is a closed set and observing that a.s.\ no new
times of percolation are introduced by this modification. The details
are left to the reader. Once we have an exceptional time with positive probability,
ergodicity immediately implies that this occurs a.s. \qed

The first moment of $X_R$ is, due to Fubini's Theorem, simply
the probability of our one-arm event, namely $\alpha_1(R)$. The second moment of
$X_R$ is easily seen to be
\begin{equation}\label{e.2nd}
\E(X^2)=
\E(\int_0^1\int_0^1 1_{0\overset{\omega_s}{\longleftrightarrow}R}
\,1_{0\overset{\omega_{t}}{\longleftrightarrow}R}\,ds\,dt)
=
\int_0^1\int_0^1 \P(0\overset{\omega_s}{\longleftrightarrow}R,
0\overset{\omega_{t}}{\longleftrightarrow}R)\,ds\,dt
\end{equation}
which is, by time invariance, at most
\begin{equation}\label{e.2ndagain}
2\int_0^1 \P(0\overset{\omega_s}{\longleftrightarrow}R,
0\overset{\omega_{0}}{\longleftrightarrow}R)\,ds.
\end{equation}

The key observation now, which brings us
back to noise sensitivity, is that the integrand
$\P(0\overset{\omega_s}{\longleftrightarrow}R,
0\overset{\omega_{0}}{\longleftrightarrow}R)$
is precisely
$\E[f_R(\omega)f_R(\omega_\epsilon)]$ where
$f_R$ is the indicator of the event that there is an open path from
the origin to distance $R$ away and $\eps=1-e^{-s}$
since looking at our process at two different times is
exactly looking at a configuration and a noisy version.

What we have seen in this subsection is that proving the existence of
exceptional times comes down to proving a second moment estimate and
furthermore that the integrand in this second moment estimate 
concerns noise sensitivity, something for which
we have already developed a fair number of tools to handle.

\section[Proof of existence of exceptional times on $\T$]{Proof of existence of exceptional times for the hexagonal
lattice via randomized algorithms}

In \cite{\SS}, exceptional times were shown to exist for the
hexagonal lattice; this was the first transitive graph for which such a
result was obtained. However, the methods in this paper did not 
allow the authors to prove that $\Z^2$ had exceptional times.

\begin{theorem}[\cite{\SS}]\label{th.SS} 
For dynamical percolation on the hexagonal lattice $\T$ at the critical point 
$p_c=1/2$, there exist almost surely exceptional times $t\in[0,\infty)$ such 
that $\omega_t$ has an infinite cluster.
\end{theorem}

\vskip 0.3 cm

\proof
As we noted in the previous section,
two different times of our model can
be viewed as ``noising'' where the probability that a hexagon is
rerandomized within $t$ units of time is $1-e^{-t}$.
Hence, by (\ref{e.correlationFourier}), we have that
\begin{equation}
\Pb{0\overset{\omega_0}{\longleftrightarrow}R,\,0\overset{\omega_t}{\longleftrightarrow}R}
= \Eb{f_R}^2 + \sum_{\emptyset\neq S \subseteq B(0,R)} \hat{f_R}(S)^2 \exp(-t|S|)
\end{equation}
where $B(0,R)$ are the set of hexagons involved in the event $f_R$.
We see in this expression that, for small times $t$, the frequencies 
contributing in the correlation between 
$\{ 0\overset{\omega_0}{\longleftrightarrow}R \}$
and $\{ 0\overset{\omega_t}{\longleftrightarrow}R \}$ are of ``small''
size $|S| \lesssim 1/t$. Therefore, in order to detect the
existence of exceptional times, one needs to achieve good control on
the {\em lower tail} of the Fourier spectrum of $f_R$.

The approach of this section is to find an algorithm minimizing the revealment
as much as possible and to apply Theorem \ref{t.ss}.
However there is a difficulty here, since our algorithm
might have to look near the origin, in which case it is difficult to
keep the revealment small. There are other reasons for a potential problem.
If $R$ is very large and $t$ very small, then if one conditions on the event 
$\{  0\overset{\omega_0}{\longleftrightarrow}R \}$, since few sites are 
updated, the open path in $\omega_0$ from $0$ to distance $R$
will still be preserved in $\omega_t$ at least up to some distance $L(t)$ 
(further away, large scale connections start to decorrelate).
In some sense the geometry associated to the event 
$\{  0\overset{\omega}{\longleftrightarrow}R \}$ is ``frozen'' on a certain 
scale between time $0$ and time $t$. 
Therefore, it is natural to divide our correlation analysis into two 
scales: the ball of radius $r=r(t)$ and the annulus
from $r(t)$ to $R$. Obviously the ``frozen radius'' $r=r(t)$ 
increases as $t\to 0$. We therefore proceed as follows instead. For
any $r$, we have
\begin{eqnarray}\label{e.sorry.christophe}
\Pb{0\overset{\omega_0}{\longleftrightarrow}R,\,0\overset{\omega_t}{\longleftrightarrow}R}& \le  & \Pb{0\overset{\omega_0}{\longleftrightarrow} r} \, \Pb{r \overset{\omega_0}{\longleftrightarrow}R,\, r \overset{\omega_t}{\longleftrightarrow}R}\nonumber \\
&\le & \alpha_1(r)\,  \Eb{f_{r,R}(\omega_0) f_{r,R} (\omega_t)},
\end{eqnarray}
where $f_{r,R}$ is the indicator function of the event, denoted by
$r \overset{\omega}{\longleftrightarrow}R$, that there is an open path
from distance $r$ away to distance $R$ away. Now, as above, we have
\begin{equation}
\Eb{f_{r,R}(\omega_0) f_{r,R} (\omega_t)} \le 
\Eb{f_{r,R}}^2 + \sum_{k=1}^{\infty}\exp(-tk) \sum_{|S|=k}\hat f_{r,R}(S)^2.
\end{equation}

The Boolean function $f_{r,R}$ somehow avoids the singularity at the origin, 
and it is possible to find algorithms for this function 
with small revealments. In any case, letting $\delta=\delta_{r,R}$ be the
revealment of $f_{r,R}$, it follows from Theorem \ref{t.ss} and the fact
that $\sum_k k\exp(-tk)\le O(1)/t^2$ that
\begin{equation}\label{e.nextterm}
\Eb{f_{r,R}(\omega_0) f_{r,R} (\omega_t)} \le 
\alpha_1(r,R)^2 +O(1)\delta \alpha_1(r,R)/t^2.
\end{equation}

The following proposition gives a bound on $\delta$. We will sketch
why it is true afterwards.

\begin{prop}[\cite{\SS}]\label{p.annulus}
Let $2\le r< R$. Then
\begin{equation}\label{e.annulus.alg}
\delta_{r,R}\le O(1)\alpha_1(r,R)\,\alpha_2(r)\,.
\end{equation}
\end{prop}

Putting together (\ref{e.sorry.christophe}), (\ref{e.nextterm}),
Proposition \ref{p.annulus} and using quasi-multiplicativity of
$\alpha_1$ yields
$$
\Pb{0\overset{\omega_0}{\longleftrightarrow}R,\,0\overset{\omega_t}{\longleftrightarrow}R}\le
O(1)\frac{\alpha_1(R)^2}{\alpha_1(r)}\left(1+\frac{\alpha_2(r)}{t^2}\right).
$$
This is true for all $r$ and $t$. If we choose $r=r(t)=(1/t)^8$ and
ignore $o(1)$ terms in the critical exponents (which can easily be
handled rigorously), we obtain, using the explicit values for the one and two-arm
critical exponents, that
\begin{equation}\label{e.correlalgo}
\Pb{0\overset{\omega_0}{\longleftrightarrow}R,\,0\overset{\omega_t}{\longleftrightarrow}R} \le O(1)t^{-5/6} \alpha_1(R)^2\,.
\end{equation}

Now, since $\int_0^1 t^{-5/6} dt <\infty$, by integrating the above 
correlation bound over the unit interval, one obtains
that $\Eb{X_R^2} \le C \Eb{X_R}^2$ for some constant $C$ as
desired.
\qed

\medbreak \noindent 
{\bf Outline of proof of Proposition \ref{p.annulus}.} 

\ni
\begin{figure}[!h]\label{f.annular}
\begin{center}
\includegraphics[width=0.7 \textwidth]{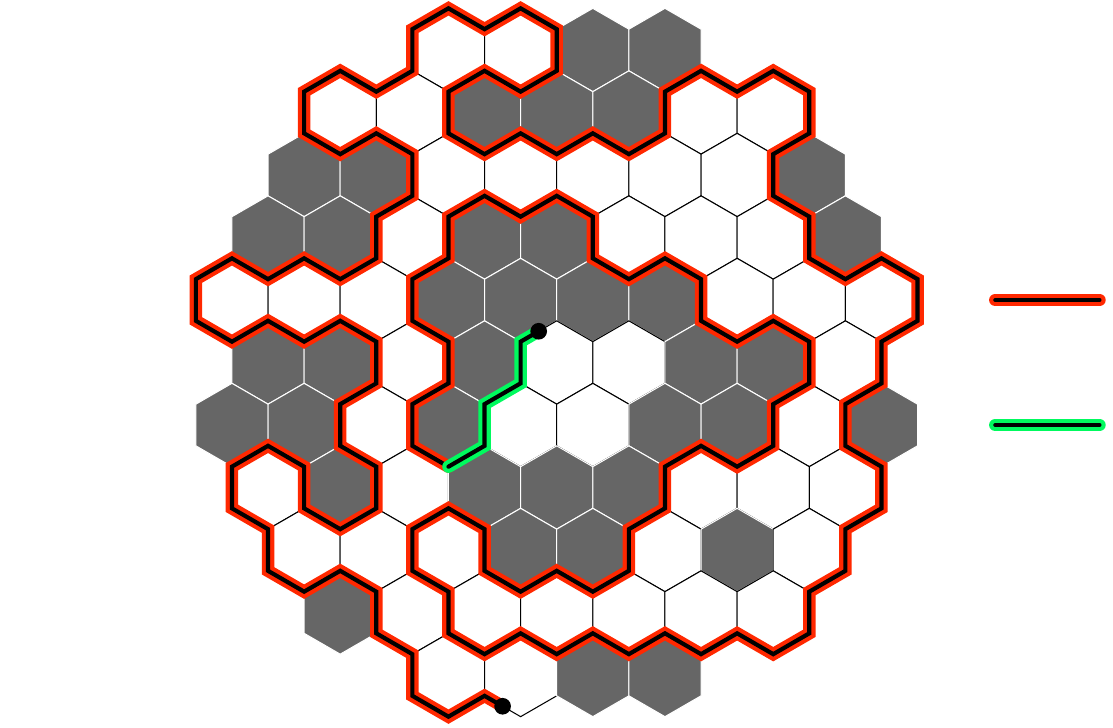}
\end{center}
\end{figure}
We use an algorithm that mimics the one we used for percolation crossings
except the present setup is ``radial''. As in the chordal case, 
we randomize the starting point of our exploration 
process by choosing a site uniformly on the `circle' of radius $R$. 
Then, we explore the picture with an exploration path $\gamma$ directed 
towards the origin; this means that as in the case of crossings, 
when the interface encounters an open (resp.\ closed) site, it turns say to the 
left (resp.\ right), the only difference being that when the exploration 
path closes a loop around the origin, it continues its exploration 
inside the connected component of the 
origin. (It is known that this discrete curve converges towards {\em radial} 
$\mathrm{SLE}_6$ on $\T$, when the mesh goes to zero.) 
It turns out that the so-defined exploration path gives all the information we 
need. Indeed, if the exploration path closes a clockwise loop around the 
origin, this means that there is a closed circuit around the origin making 
$f_{r,R}$ equal to zero. On the other hand, if the exploration path does 
not close any clockwise loop until it reaches radius $r$, it means 
that $f_{r,R}=1$. Hence, we run the exploration path until
either it closes a clockwise loop or it reaches radius $r$. This
is our algorithm. Neglecting boundary issues (points near radius $r$ or $R$), 
if $x$ is a point at distance $u$ from 0, with $2r < u < R/2$, in order 
for $x$ to be examined by the algorithm, it is needed that 
there is an open path from $2u$ to $R$ and the two-arm event holds in
the ball centered at $u$ with radius $u/2$. Hence for $|x|=u$,
$\Pb{x\in J}$ is at most
$O(1)\alpha_2(u) \alpha_1(u,R)$. Due to the explicit values of the one 
and two-arm exponents, this expression is decreasing in $u$. Hence, ignoring
the boundary, the revealment is at most
$O(1)\alpha_2(r) \alpha_1(r,R)$. See \cite{\SS} for more details.
\qed

We now assume that the reader is familiar with the notion of
Hausdorff dimension. We let $\Exc \subseteq [0,\infty]$ denote the 
(random) set of these exceptional times at which percolation occurs.
It is an immediate consequence of Fubini's Theorem that 
$\Exc$ has Lebesgue measure zero and hence we should look at its 
Hausdorff dimension if we want to measure its ``size''. The first result
is the following.

\begin{theorem}[\cite{\SS}]\label{th.SS.HD} 
The Hausdorff dimension of $\Exc$ is an almost sure constant in $[1/6, 31/36]$.
\end{theorem}

It was conjectured there that the dimension of the set of exceptional times is
a.s.\ $31/36$.

\proofoutline
The fact that the dimension is an almost sure constant follows
from easy 0-1 Laws. The lower bounds are obtained by placing a random
measure on $\Exc$ with finite so-called $\alpha$--energies for any $\alpha <1/6$
and using a result called Frostman's Theorem. (This is a standard technique once one 
has good control of the correlation structure.)
Basically, the $1/6$ comes from the fact
that for any  $\alpha <1/6$, one can multiply the integrand in
$\int_0^1 t^{-5/6} dt$ by $(1/t)^{\alpha}$ and still be integrable.
It is the amount of ``room to spare'' you have. If one could obtain
better estimates on the correlations, one could thereby improve the
lower bounds on the dimension. The upper bound is obtained via
a first moment argument similar to the proof of Theorem
\ref{th:Zd} but now using (\ref{e.5/36}). 
\qed

Before moving on to our final method of dealing with the spectrum, let us consider what we
might have lost in the above argument. Using the above argument, we optimized things
by taking $r(t)=(1/t)^8$. However, at time $t$ compared to time 0, we have noise which is
about $t$. Since we now know the exact noise sensitivity exponent, in order to obtain 
decorrelation, the noise level should be at least about the negative $3/4$th power of the 
radius of the region we are looking at. So, events in our annulus should
decorrelate if $r(t)>> (1/t)^{4/3}$. This suggests there might be potential
for improvement. Note we used an inner radius which is 6 times larger than
potentially necessary ($8=6\times 4/3$). This 6 is the same 6 by which
the result in Theorem \ref{th:crossingquant} differed by the true exponent 
($3/4=6\times 1/8$) and the same 6 explaining the gap in Theorem \ref{th.SS.HD}
($1-1/6)=6\times (1-31/36$). This last difference is also seen by comparing the exponents in
(\ref{e.correlalgo}) and the last term in (\ref{e.betterexponent}) below.

\section[Exceptional times via the geometric approach]{Proof of existence of exceptional times via the geometric approach of 
the spectrum}

Recall that our third approach for proving the noise sensitivity of
percolation crossings was based on a geometrical analysis of the
spectrum, viewing the spectrum as a random set. This approach
yielded the exact noise sensitivity exponent for percolation crossings
for the hexagonal lattice. This approach can also be used here as we
will now explain. Two big advantages of this approach are that it succeeded in
proving the existence of exceptional times for percolation crossings on
$\Z^2$, something which \cite{\SS} was not able to do, as well as obtaining
the exact Hausdorff dimension for the set of exceptional times, namely
the upper bound of $31/36$ in the previous result. 

\begin{theorem}[\cite{\GPS}]\label{th.HD} 
For the triangular lattice, the Hausdorff dimension of $\Exc$ is almost surely $31/36$.
\end{theorem}

\proof
As explained in the previous section, it suffices to lower the $5/6$ in
(\ref{e.correlalgo}) to $5/36$. (Note that (\ref{e.correlalgo}) was really only obtained
for numbers strictly larger than $5/6$, with the $O(1)$ depending on this number;
the same will be true for the $5/36$.)

Let $s(r)$ be the inverse of the map
$r\rightarrow r^2\alpha_4(r)\sim r^{3/4}$. So more or less, $s(r):=r^{4/3}$.
Using Theorem \ref{th.spectrumonearm}, we obtain the following: 
\begin{eqnarray} \label{e.betterexponent} 
\Eb{f_{R}(\omega_0) f_{R} (\omega_t)} &=& \sum_{S}\exp(-t|S|)\hat f_{R}(S)^2
\nonumber \\
&=&\sum_{k=1}^{\infty}\sum_{S:|S|\in [(k-1)/t,k/t)}\exp(-t|S|)\hat f_{R}(S)^2
\nonumber \\
&\le &\sum_{k=1}^{\infty}\exp(-k)\SQb{|\Spec_{f_R}|<k/t}
\nonumber \\
&\le &O(1)\sum_{k=1}^{\infty}\exp(-k)\frac{\alpha_1(R)^2}{\alpha_1(s(k/t))}
\nonumber \\
&\le &O(1)\alpha_1(R)^2  \sum_{k=1}^{\infty}\exp(-k) (\frac{k}{t})^{4/3\times 5/48}
\nonumber \\
&\le &O(1)\alpha_1(R)^2  (\frac{1}{t})^{5/36}.
\end{eqnarray}
This completes the proof. (Of course, there are $o(1)$ terms in these exponents which
we are ignoring.)
\qed

We have done a lot of the work for proving that there are exceptional times also on
$\Z^2$.

\begin{theorem}[\cite{\GPS}]\label{th.gpsz2} 
For dynamical percolation on $\Z^2$ at the critical point 
$p_c=1/2$, there exist almost surely exceptional times $t\in[0,\infty)$ such 
that $\omega_t$ has an infinite cluster.
\end{theorem}

\proof
$s(r)$ is defined as it was before but now we cannot say that $s(r)$ is about 
$r^{4/3}$. However, we can say that for some fixed  $\delta> 0$, we have
that for all  $r$,
\begin{equation}\label{e.powerlower}
s(r)\ge r^{\delta}
\end{equation}

From the previous proof, we still have
\begin{equation}
\frac{\Eb{f_{R}(\omega_0) f_{R} (\omega_t)}}{\alpha_1(R)^2}
\le O(1)\sum_{k=1}^{\infty}\exp(-k)\frac{1}{\alpha_1(s(k/t))}.
\end{equation}

Exactly as in the proof of Theorem \ref{th.SS}, we need to show that the right
hand side is integrable near 0 in order to carry out the second moment argument.

Quasi-multiplicativity can be used to show that
\begin{equation}
\alpha_1(s(1/t))\le k^{O(1)}\alpha_1(s(k/t)).
\end{equation}
(Note that if things behaved exactly as power laws, this would be clear.)

Therefore the above sum is at most
\begin{equation}
O(1)\sum_{k=1}^{\infty}\exp(-k)\frac{k^{O(1)}}{\alpha_1(s(1/t))}
\le 
O(1)\frac{1}{\alpha_1(s(1/t))}
\end{equation}

V. Beffara has shown that there exists $\epsilon_0 >0$ such that for all $r$,
\begin{equation}
\alpha_1(r)\alpha_4(r)\ge r^{\epsilon_0-2}.
\end{equation}
 
Note that Theorem \ref{t.5armexponent} and (\ref{e.reimerapplication}) tell us that
the left hand side is larger than $\Omega(1)r^{-2}$. The above tells us
that we get an (important) extra power of $r$ in (\ref{e.reimerapplication}).

It follows that 
\begin{equation}
\frac{1}{\alpha_1(s(1/t))}\le \alpha_4(s(1/t))s(1/t)^{2-\epsilon_0}
= (1/t) s(1/t)^{-\epsilon_0}.
\end{equation}

(\ref{e.powerlower}) tells us that the last factor is at most $t^\eta$ for some
$\eta>0$ and hence the relevant integral converges as desired. The rest of the
argument is the same.
\qed

One can also consider exceptional times for other events, such as for example
times at which there is an infinite cluster in the upper half-plane or 
times at which there are two infinite clusters in the whole plane,
and consider the corresponding Hausdorff dimension. A number of results of
this type, which are not sharp, are given in \cite{\SS} while various sharp results are given 
in \cite{\GPS}. 

\chapter*{Exercise sheet of Chapter \ref{ch.DP}}
\setcounter{exercise}{0}

\begin{exercise}
Prove that on any graph below criticality, there are no times
at which there is an infinite cluster while
above criticality,  there is an infinite cluster at all times.
\end{exercise}

\begin{exercise}
Consider critical dynamical percolation on a general graph 
satisfying $\theta(p_c)=0$. Show that a.s.\
$\{t: \omega_t \mbox{ percolates }\}$ has Lebesgue measure 0.
 \end{exercise}

\begin{exercise} (Somewhat hard).
A {\em spherically symmetric} tree is one where all vertices
at a given level have the same number of children, although
this number may depend on the given level. Let $T_n$ be the number
of vertices at the $n$th level. Show that there is percolation at $p$
if
$$
\sum_n \frac{1}{p^{-n}T_n} < \infty
$$

Hint: Let $X_n$ be the number of vertices in the $n$th level which are
connected to the root. Apply the second moment method to the sequence
of $X_n$'s.

The convergence of the sum is also necessary for percolation but
this is harder and you are not asked to show this. 
This theorem is due to Russell Lyons.
\end{exercise}

\begin{exercise}
Show that if $T_n$ is $n^2 2^n$ up to
multiplicative constants, then the critical value of the graph is 
$1/2$ and we percolate at the critical value. (This yields a graph
which percolates at the critical value.)
\end{exercise}

\begin{exercise}  (Quite a bit harder).
Consider dynamical percolation on a spherically symmetric tree.
Show that there for the parameter $p$, there are exceptional times
at which percolation occurs if
$$
\sum_n \frac{1}{np^{-n}T_n} < \infty.
$$
Hint: Find an appropriate random variable $X_n$ to which the second
moment method can be applied.
\end{exercise}

\begin{exercise}
Find a spherically symmetric tree which does not percolate
at criticality but for which there are exceptional times at which
percolation occurs.
\end{exercise}


\note{Old ideas of exercises here:}
{
\bi
\item[-] Problem session on the $d$-regular tree $\T^d$ case
\item[-] Problem session on the $\Z^d$ case, high $d$.
\item[-] Do the Half-plane case properly
\item[-] Conclude that there exist ``super exceptional'' times ! (either monochromatic or not)
\item[-] Exceptional times with 4 arms ? Disprove !
\ei
}

\bibliographystyle{alpha}
\bibliography{book}
\vskip 1 cm

\hskip 0.5 \textwidth
\begin{minipage}{0.5 \textwidth}
\noindent {\bf Christophe Garban}\\
CNRS and UMPA\\
Ecole Normale Sup\'erieure de Lyon, UMPA \\
46 all\'ee d'Italie \\
69364 Lyon Cedex 07 France \\
{christophe.garban@ens-lyon.fr}\\
\url{http://www.umpa.ens-lyon.fr/\string~cgarban/}.

\vskip 0.5 cm \noindent {\bf Jeffrey E. Steif}\\
Mathematical Sciences \\
Chalmers University of Technology \\
and \\
Mathematical Sciences \\
G\"{o}teborg University \\
SE-41296 Gothenburg, Sweden \\
{steif@math.chalmers.se} \\
\url{http://www.math.chalmers.se/\string~steif/}
\end{minipage}

\end{document}